\documentclass[12pt]{amsart}
\font\emailfont=cmtt10

\usepackage{amsmath,amsthm,amsfonts,amscd,flafter,epsf}

\newcommand\commentable[1]{#1}

\newcommand{\Tors}{\mathrm{Tors}}
\newcommand{\rk}{\mathrm{rk}}
\newcommand{\HF}{HF}

\newtheorem{theorem}{Theorem}[section]
\newtheorem{prop}[theorem]{Proposition}
\newtheorem{cor}[theorem]{Corollary}
\newtheorem{conj}[theorem]{Conjecture}
\newtheorem{lemma}[theorem]{Lemma}

\newtheorem{defn}[theorem]{Definition}

\newtheorem{remark}[theorem]{Remark}

\def\endproof{\relax\ifmmode\expandafter\endproofmath\else
  \unskip\nobreak\hfil\penalty50\hskip.75em\hbox{}\nobreak\hfil\bull
  {\parfillskip=0pt \finalhyphendemerits=0 \bigbreak}\fi}
\def\endproofmath$${\eqno\bull$$\bigbreak}
\def\bull{\vbox{\hrule\hbox{\vrule\kern3pt\vbox{\kern6pt}\kern3pt\vrule}\hrule}}

\newcommand{\Q}{\mathbb{Q}}
\newcommand{\R}{\mathbb{R}}

\newcommand{\C}{\mathbb{C}}

\newcommand{\Z}{\mathbb{Z}}

\newcommand{\OneHalf}{\frac{1}{2}}

\newcommand{\Zmod}[1]{\Z/{#1}\Z}

\newcommand{\Ker}{\mathrm{Ker}}

\newcommand{\cm}{\cdot}

\newcommand{\Nbd}[1]{{\mathrm{nd}}(#1)}
\newcommand{\nbd}[1]{\Nbd{#1}}
\newcommand{\CDisk}{D}

\newcommand{\ModSWfour}{\mathcal{M}}
\newcommand{\ModFlow}{\ModSWfour}

\newcommand{\SpinC}{{\mathrm{Spin}}^c}

\newcommand\Wedge{\Lambda}

\newcommand\abuts\Rightarrow
\newcommand\Sym{\mathrm{Sym}}

\newcommand\HFpRed{\HFp_{\red}}

\newcommand\RelSpinC{\underline{\SpinC}}
\newcommand\relspinc{\underline{\spinc}}

\newcommand\Filt{\mathcal F}

\newcommand\CFinfty{\CFinf}

\newcommand\x{\mathbf x}
\newcommand\w{\mathbf w}

\newcommand\y{\mathbf y}

\newcommand\ModSphere{\ModFlow\left({\mathbb S}\longrightarrow 
\Sym^{g-1}(\Sigma_{1})\times \Sym^2(\Sigma_{2})\right)}
\newcommand\ModSpheres\ModSphere
\newcommand\CF{CF}

\newcommand\CFa{\widehat{CF}}
\newcommand\CFp{\CFb}
\newcommand\CFm{\CF^-}

\newcommand\HFleq{\HF^{\leq 0}}

\newcommand{\red}{\mathrm{red}}

\newcommand\HFp{\HFb}

\newcommand\HFm{\HF^-}
\newcommand\CFinf{CF^\infty}
\newcommand\HFinf{HF^\infty}
\newcommand\CFb{CF^+}
\newcommand\HFa{\widehat{HF}}
\newcommand\HFb{HF^+}

\newcommand\Mas{\mu}
\newcommand\UnparModSp{\widehat \ModSp}
\newcommand\UnparModFlow\UnparModSp
\newcommand\Mod\ModSp

\newcommand{\cald}{{\mathcal D}}

\newcommand\PD{\mathrm{PD}}

\newcommand{\spinc}{\mathfrak s}
\newcommand{\spincfour}{\mathfrak r}

\newcommand{\spinct}{\mathfrak t}

\newcommand\ModMaps{\mathcal M}
\newcommand\ModSp\ModMaps

\newcommand\Ta{{\mathbb T}_{\alpha}}
\newcommand\Tb{{\mathbb T}_{\beta}}
\newcommand\Tc{{\mathbb T}_{\gamma}}

\newcommand\alphas{\mbox{\boldmath$\alpha$}}

\newcommand\betas{\mbox{\boldmath$\beta$}}
\newcommand\gammas{\mbox{\boldmath$\gamma$}}

\newcommand\PerDom{\mathcal P}

\newcommand\HFc{\HF^\circ}

\newcommand\Dual{\mathcal D}
\newcommand\Duality\Dual

\newcommand\phiCov{\widetilde \phi}
\newcommand\xCov{\widetilde x}
\newcommand\wCov{\widetilde w}
\newcommand\yCov{\widetilde y}
\newcommand\qCov{\widetilde q}
\newcommand\alphaCov{\widetilde \alpha}
\newcommand\betaCov{\widetilde \beta}
\newcommand\gammaCov{\widetilde \beta}
\newcommand\piCov{{\widetilde \pi}_2}

\newcommand\spincrel\relspinc
\newcommand\relspinct{\underline{\mathfrak t}}

\newcommand\CFKsa{CFK^0}
\newcommand\bCFKsa{^bCFK^0}
\newcommand\spincf{\mathfrak r}

\newcommand\Mass{\overline n}

\newcommand\bCFKp{{}^b\CFKp}
\newcommand\bCFKm{{}^b\CFKm}

\newcommand\bPhim{{}^b\Phi^-}
\newcommand\Phip{\Phi^+}

\newcommand\CFK{CFK}
\newcommand\HFK{HFK}

\newcommand\CFKp{\CFK^+}
\newcommand\CFKa{\widehat\CFK}
\newcommand\CFKm{\CFK^-}
\newcommand\CFKinf{\CFK^{\infty}}

\newcommand\HFKa{\widehat\HFK}

\newcommand\HFKinf{\HFK^{\infty}}

\newcommand\Mark{m}

\newcommand\BasePt{w}
\newcommand\FiltPt{z}

\newcommand\StartPt{w}
\newcommand\EndPt{z}

\newcommand\bPsip{{}^b\Psi^+}
\newcommand\bPsia{{}^b{\widehat\Psi}}

\newcommand\Psim{\Psi^-}

\commentable{

\title[{Holomorphic disks and knot invariants}] 
{Holomorphic disks and knot invariants}

\author[Peter Ozsv{\'a}th]{Peter Ozsv\'ath}
\address{Department of
Mathematics, Columbia University, New York 10027 \newline
\indent{\emailfont{petero@math.columbia.edu}}}
\thanks{PSO was supported by NSF grant number DMS 9971950 and a Sloan 
Research Fellowship}

\author[Zolt{\'a}n Szab{\'o}]{Zolt{\'a}n Szab{\'o}} 
\address{Department of
Mathematics, Princeton University, New Jersey 08540 \newline
\indent{\emailfont{szabo@math.princeton.edu}}}}
\thanks{ZSz was supported by NSF grant number DMS 0107792
and a Packard Fellowship}

\begin{document}

\begin{abstract}  
  We define a Floer-homology invariant for knots in an oriented
  three-manifold, closely related to the Heegaard Floer homologies for
  three-manifolds defined in an earlier paper. We set up basic
  properties of these invariants, including an Euler characteristic
  calculation, and a description of the
  behaviour under connected sums. Then, we establish a
  relationship with $\HFp$ for surgeries along the knot. Applications
  include calculation of $\HFp$ of three-manifolds obtained by
  surgeries on some special knots in $S^3$, and also calculation of
  $\HFp$ for certain simple three-manifolds which fiber over the
  circle.
\end{abstract}

\maketitle
\section{Introduction}

The purpose of this paper is to define and study an invariant for
null-homologous knots $K$ in an oriented three-manifold $Y$.  These
naturally give rise to invariants for oriented, null-homologous links
in $Y$, since oriented, null-homologous links in $Y$ correspond to
certain oriented, null-homologous knots in $Y\#^{n-1}(S^2\times S^1)$
(where $n$ denotes the number of components of the original link).  In
the interest of exposition, we describe first some special cases of
the construction in the case where the ambient three-manifold is $S^3$
-- the case of ``classical links.''

\subsection{Classical links.}
\label{subsec:ClassLinks}
Suppose $Y\cong S^3$, and suppose $L$ is an oriented link. In its simplest
form, our construction gives a sequence of graded Abelian groups
$\HFKa(L,i)$, where here $i\in \Z$. If $L$ has an odd number of components,
the grading is a $\Z$-grading, while if it has an even number of components,
the grading function takes values in $\OneHalf+\Z$. 

These homology groups satisfy a number of basic properties, which we
outline presently. Sometimes, it is simplest to state these properties 
for $\HFKa(L,i,\Q)$, the homology with rational coefficients:
$\HFKa(L,i,\Q)\cong \HFKa(L,i)\otimes_{\Z}\Q$.

First, the Euler characteristic is related to the Alexander-Conway
polynomial of $L$, $\Delta_L(T)$ by the following formula:
\begin{equation}
\label{eq:EulerChar}
\sum \chi(\HFKa(L,i,\Q))\cm T^{i} = (T^{-1/2}-T^{1/2})^{n-1} \cm \Delta_L(T),
\end{equation}
where $n$ denotes the number of components of $L$
(it is interesting to compare this with~\cite{Casson}, \cite{MengTaubes},
and~\cite{FSknots}).
The sign conventions on the Euler characteristic here are given by
\begin{eqnarray*}
\chi(\HFKa(L,i,\Q))
&=& \sum_{d \in  \left(\frac{n-1}{2}\right) + \Z}
(-1)^{d+\frac{n}{2}-\OneHalf} \rk \left(\HFKa_d(L,i,\Q)\right).
\end{eqnarray*}

Our grading conventions are justified by the following property.  If
${\overline L}$ denotes the mirror of $L$ (i.e. switch over- and
under-crossings in a projection for $L$), then
\begin{equation}
\label{eq:Reflection}
\HFKa_d(L,i,\Q)\cong \HFKa_{-d}({\overline L},-i,\Q).
\end{equation}
Another symmetry these invariants enjoy is the following conjugation symmetry:
\begin{equation}
\label{eq:KnotConjInv}
\HFKa_d(L,i,\Q)\cong \HFKa_{d-2i}(L,-i,\Q),
\end{equation}
refining the symmetry of the Alexander polynomial.
Finally, the invariants remain unchanged after an overall orientation reversal
of the link,
\begin{equation}
\label{eq:LinkOrientationReversal}
\HFKa_d(L,i,\Q)\cong \HFKa_d(-L,i,\Q).
\end{equation}

These groups also satisfy a K\"unneth principle for connected
sums. Specifically, let $L_1$ and $L_2$ be a pair of disjoint,
oriented links which can be separated from one another by a
two-sphere. Choose a component of $L_1$ and $L_2$, and let $L_1\# L_2$
denote the connected sum performed at these components. Then,
\begin{equation}
\label{eq:KunnethFormula}
\HFKa(L_1\#L_2,i,\Q)\cong
\bigoplus_{i_1+i_2=i} \HFKa(L_1,i_1,\Q)\otimes_\Q
\HFKa(L_2,i_2,\Q)\end{equation}
(see Corollary~\ref{cor:Kunneth} for a more general statement).  Of
course, this can be seen as a refinement of the fact that the
Alexander polynomial of links is multiplicative under connected sums.

Let $L_1$ and $L_2$ be two links as above. Then, we can also form
their disjoint union $L_1\coprod L_2$. We have that 
\begin{equation}
\label{eq:DisjointUnion}
\HFKa(L_1\coprod
L_2,i,\Q)\cong \HFKa(L_1\# L_2,i,\Q)\otimes W
\end{equation} where $W$ is a
two-dimensional graded vector space splitting into two one-dimensional
pieces $$W\cong W_{-1/2}\oplus W_{1/2},$$
where $W_{\pm 1/2}$ has
grading $\pm 1/2$.  The above is a manifestation of the advantage of
$\HFKa$ over the Alexander polynomial: the Alexander polynomial of any
split link vanishes.

\begin{figure}
\mbox{\vbox{\epsfbox{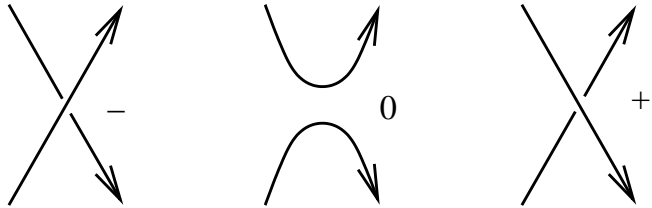}}}
\caption{\label{fig:Skein}
Skein moves at a double-point.}
\end{figure}

These invariants also satisfy a ``skein exact sequence''
(compare~\cite{FloerKnots}, \cite{BraamDonaldson}, \cite{AbsGraded},
\cite{Khovanov}). Suppose that $L$ is a link, and suppose that $p$ is
a positive crossing of some projection of $L$. Following the usual
conventions from skein theory, there are two other associated links,
$L_0$ and $L_-$, where here $L_-$ agrees with $L_+$, except that the
crossing at $p$ is changed, while $L_0$ agrees with $L_+$, except that
here the crossing $p$ is resolved. These three cases are illustrated
in Figure~\ref{fig:Skein}. There are two cases of the skein exact
sequence, according to whether or not the two strands of $L_+$ which
project to $p$ belong to the same component of $L_+$.

Suppose first that the two strands which project to $p$ belong to the
same component of $L_+$. In this case, the skein exact sequence reads:
\begin{equation}
\label{eq:SkeinExact}
\begin{CD}
...@>>>\HFKa(L_{-})@>>>\HFKa(L_0)@>>>\HFKa(L_{+})@>>>...,
\end{CD}
\end{equation}
where all the maps above respect the splitting of $\HFKa(L)$ into
summands (e.g. $\HFKa(L_-,i)$ is mapped to $\HFKa(L_0,i)$).
Furthermore, the maps to and from $\HFKa(L_0)$ drop degree by
$\OneHalf$. The remaining map  from
$\HFKa(L_+)$ to $\HFKa(L_-)$ does not necessarily respect the 
absolute grading; however, it can be expressed as a sum of
homogeneous maps\footnote{A graded group is one which splits as
$A=\bigoplus_{k\in\Q} A_k$. A homogeneous element is one which 
is supported in $A_k$ for some $k$. A map between graded groups
$\phi\colon A\longrightarrow B$
is said to be
{\em homogeneous of degree $d$} if $\phi$ maps all the $A_k$ into $B_{k+d}$. 
The map $\phi$ is said to be homogeneous if it is homogeneous
of degree $d$, for some $d$.},
none of which increases absolute grading.
When the two strands belong to
different components, we obtain the following:
\begin{equation}
\label{eq:SkeinExact2}
\begin{CD}
...@>>>\HFKa(L_{-})@>>>\HFKa(L_0)\otimes V @>>>\HFKa(L_{+})@>>>...,
\end{CD}
\end{equation}
where $V$ denotes the four-dimensional vector space $$V=V_{-1}\oplus
V_0\oplus V_1,$$ where here $V_{\pm 1}$ are one-dimensional pieces
supported in degree $\pm 1$, while $V_0$ is a two-dimensional piece
supported in degree $0$.  Moreover, the maps respect the decomposition
into summands, where the $i^{th}$ summand of the middle piece
$\HFKa(L_0)\otimes V$ is given by 
$$\left(\HFKa(L_0,i-1)\otimes V_1\right) 
\oplus \left(\HFKa(L_0,i)\otimes V_0\right) 
\oplus \left(\HFKa(L_0,i+1)\otimes V_{-1}\right).
$$
The shifts in the 
absolute gradings work just as they did in the previous case.

In addition to establishing the above properties, we calculate this
invariant for some examples. For a more detailed discussion of the
link invariants in the classical case, we refer the reader to the
sequel,~\cite{ClassicalKnots}. In particular, in that article, we
determine the invariants for alternating links. In a different
direction, if $L$ is a fibered link and $d$ denotes the degree of
$\Delta_L(T)$, then we have that
$\HFKa(L,d+\left(\frac{n-1}{2}\right))\cong \Z$. This statement is a
generalization of the fact that the Alexander-Conway polynomial of a
fibered link is monic.  We do not prove this fact in the present
paper, but rather give the proof in~\cite{Contact}.

Suppose now that our link consists of a single component, i.e. it is a
knot $K$.  In this case,
the homology groups also
 give bounds on the genus $g$ of the knot $K$:
if $\HFKa(K,s)\neq 0$, then $$|s|\leq g.$$

Since the knot homology groups can be viewed as a refinement of the
Alexander polynomial, it is natural to ask whether they contain more
information (as they did in the case of links). In fact, it is easy to
see that they do -- as graded Abelian groups, $\HFKa$ distinguishes
the right-handed from the left-handed trefoil.  More interestingly,
building on the material in~\cite{ClassicalKnots},
we show in~\cite{calcKT} that $\HFKa$
distinguishes all non-trivial pretzel knots of the form $P(2a+1,2b+1,2c+1)$ from
the unknot (note that these include knots with trivial Alexander
polynomial).  Indeed, in that paper, it is also shown that $\HFKa$ can 
distinguish knots which differ by a Conway mutation.

Moreover, we conjecture that the genus bounds stated earlier are sharp:

\begin{conj}
If $K\subset S^3$ is a knot with genus $g$, then
$$\HFKa(K,g)\neq 0.$$
\end{conj}

Our evidence for this conjecture is based on three ingredients.  The
first is a relationship between the knot Floer homology groups and
Heegaard Floer homology, which is established in
Section~\ref{sec:Relationship}. Next, we appeal to the conjectured
relationship between Heegaard Floer homology and Seiberg-Witten Floer
homology (c.f.~\cite{HolDisk}). And finally, non-triviality of the 
corresponding Seiberg-Witten objects  is established by work of
Kronheimer-Mrowka~\cite{KMthurston}, which in turn rests on work of
Gabai~\cite{Gabai} and Eliashberg-Thurston~\cite{EliashbergThurston}.

\subsection{General case}

In the above discussion, we restricted attention to a rather special
case of our constructions, in the interest of exposition. As we
mentioned, the constructions we give here actually generalize to the
case where the ambient three-manifold $Y$ is an arbitrary closed, oriented
three-manifold and $K$ is a null-homologous knot; and indeed, the algebra
also gives more information than simply the knot homology groups. In
their more general form, the constructions give a
$\Z\oplus\Z$-filtration on the chain complex $\CFinf(Y,\spinc)$ used
for calculating the Heegaard homology $\HFinf(Y,\spinc)$ defined
in~\cite{HolDisk} (where here $\spinc$ is any $\SpinC$ structure over
$Y$). This filtration also induces a $\Z$ filtration on
$\CFa(Y,\spinc)$; and the knot homology groups we described above are
the homology groups of the associated graded complex (in the case
where $Y\cong S^3$).

The paper is organized as follows. In Section~\ref{sec:TopPre}, we set
up the relevant topological preliminaries: the passage from links to
knots, Heegaard diagrams and knots, and some of the algebra of
filtered chain complexes. In Section~\ref{sec:DefKnot} we give the
definition of the knot filtration, and prove its topological
invariance. In that section, we also establish some of the symmetries
of the invariants.  In Section~\ref{sec:Relationship}, we explain how
the filtration can be used to calculate the Floer homology groups of
three-manifolds obtained by ``large surgeries'' along the knot $K$. In
Section~\ref{sec:Adj}, we derive an ``adjunction inequality'' which
relates $\HFKa$ and the genus of the knot. In
Section~\ref{sec:Examples}, we give some sample calculations for the
knot homology groups of some special classes of classical knots. With
this done, we return to some general properties: the K\"unneth
principle for connected sums, and the surgery long exact sequences in
Sections~\ref{sec:ConnectedSums} and
\ref{sec:Sequences} respectively. As an application of these general results,
we conclude with some calculations of $\HFp(Y,\spinc)$ for certain
simple three-manifolds which fiber over the circle (in the case where
the first Chern class of the $\SpinC$ structure $\spinc$ evaluates
non-trivially on the fibers). The three-manifolds considered here
include $Y\cong S^1\times
\Sigma_g$ for arbitrary $g$, and also mapping tori of a single,
non-separating Dehn twist.  Finally in Section~\ref{sec:Links} we
point out how the results from the paper specialize to the properties
and formulae stated in Subsection~\ref{subsec:ClassLinks}.

We will return to other applications of these constructions in future
papers.  For example, the constructions play a central role
in~\cite{Contact}, where they are used to define invariants of contact
structures for three-manifolds.

\subsection{Further remarks.}
The fact that a knot in $Y$ induces a filtration on $\CFinf(Y)$ has
been discovered independently by Rasmussen in~\cite{Rasmussen}, where
he uses the induced filtration to compute the Floer homologies of
three-manifolds obtained as surgeries on two-bridge knots. 
In fact, we have learned that many of the constructions presented here
have been independently discovered by Rasmussen in his thesis~\cite{RasmussenThesis}.

In another direction, it worth pointing out the striking similarity
between the knot homology groups described here and Khovanov's
homology groups~\cite{Khovanov} (see also~\cite{BarNatan}). Although
our constructions here are quite different in spirit from Khovanov's,
the final result is analogous: we have here a homology theory whose
Euler characteristic is the Alexander polynomial, while Khovanov
constructs a homology theory whose Euler characteristic is the Jones
polynomial. It is natural to ask whether there is a simultaneous
generalization of these two homology theories.

\subsection{Acknowledgments.}
The authors wish to warmly thank Paolo Lisca, Tomasz Mrowka, Jacob
Rasmussen, and Andr{\'a}s Stipsicz for some interesting
discussions. We would also like to thank the referee for some helpful
comments.

\section{Topological preliminaries}
\label{sec:TopPre}

\subsection{From knots  to links}
\label{sec:FromKnotsToLinks}
Although we discussed links in the introduction, we will focus mainly
on the case of  knots in the paper. Our justification for the
apparent loss of generality is the observation that oriented
$n$-component links in $Y$ correspond to oriented knots in
$Y\#^{n-1}(S^2\times S^1)$. For related constructions,
see~\cite{Hoste}; compare also~\cite{BraamDonaldson} and~\cite{FSknots}.

More precisely, suppose that $L$ is an $n$-component oriented
three-manifold $Y$. Then, we can construct an oriented knot in
$Y\#^{n-1}(S^2\times S^1)$ as follows. Fix $2n-2$ points
$\{p_i,q_i\}_{i=1}^{n-1}$ in $L$, which are then pairwise grouped
together (i.e., $n-1$ embedded zero-spheres in $L$), in such a manner
that if we formally identify each $p_i$ with $q_i$ in $L$, we obtain a 
connected graph.  We view $p_i$ $q_i$ as the feet of a
one-handle to attach to $Y$. Let $\kappa(Y,\{p_i,q_i\})$ denote the
new three-manifold. Of course, $\kappa(Y,\{p_i,q_i\})\cong
Y\#^{n-1}(S^2\times S^1)$. Now, inside each one-handle, we can find a
band along which to perform a connected sum of
the component of $L$ containing $p_i$ with the component containing
$q_i$. We choose the band so that the induced orientation of its boundary
is compatible with the orientation of $L_1$ and $L_2$.
Our hypotheses on the number and distribution of the
distinguished points ensures that the newly-constructed link (after
performing all $n-1$ of the connected sums) is a single-component knot,
which we denote $\kappa(L,\{p_i,q_i\})$, inside $\kappa(Y,\{p_i,q_i\})$

\begin{prop}
\label{prop:LinksToKnots}
The above construction induces a well-defined map from the set of
isotopy classes of oriented, $n$-component links in $Y$ to the set of
isotopy classes of oriented knots in $Y\#^{n-1}(S^2\times S^1)$. More
precisely, if $L$ and $L'$ are isotopic links, an
$\{p_i,q_i\}_{i=1}^{n-1}$ and $\{p_i',q_i'\}_{i=1}^{n-1}$ are
auxiliary choices of points as above, then there is an
orientation-preserving diffeomorphism
$$ \kappa(Y,\{p_i,q_i\})\longrightarrow \kappa(Y,\{p_i',q_i'\})$$
which carries the oriented knot $\kappa(L,\{p_i,q_i\})$ to a knot
which is isotopic to $\kappa(L,\{p_i',q_i'\})$.
\end{prop}

\begin{proof}
First, we argue that for a fixed set of points
$\{p_i,q_i\}_{i=1}^{n-1}$ in $L$ as above, the isotopy class of the
induced knot in $Y\#^{n-1}(S^2\times S^1)$ is independent of the
choices of bands connecting the $p_i$ to $q_i$ (inside the
one-handles).  This can be seen using an isotopy supported inside the
one-handle.

Next, we argue that the isotopy class of $\kappa(L)$ is independent of
the choice of points $\{p_i,q_i\}_{i=1}^{n-1}$. 

The key observation is that if we have three components $K_1$, $K_2$,
$K_3$, and distinguished points $p_1\in K_1$ with $q_1\in K_2$, and
$p_2\in K_1$ and $q_2\in K_3$, then we can handleslide the one-handle
specified by $\{p_2,q_2\}$ over $\{p_1,q_1\}$ to obtain a new
one-handle $\{p_3,q_3\}$, where $p_3\in K_2$, $q_3\in K_3$
(c.f. Figure~\ref{fig:HandleSlide}).

\begin{figure}
\mbox{\vbox{\epsfbox{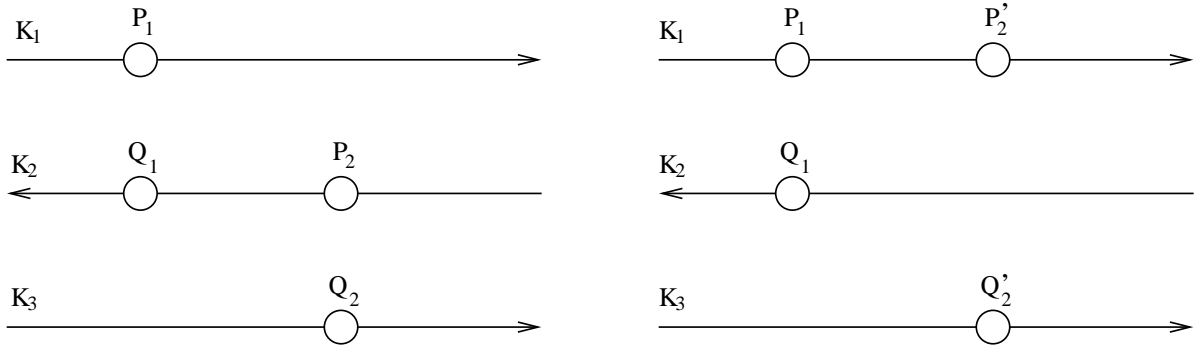}}}
\caption{\label{fig:HandleSlide}
{\bf{Handleslides.}} 
An illustration of the handleslide used in Proposition~\ref{prop:LinksToKnots}.
We slide the handle specified by $\{p_2,q_2\}$ over that specified by
$\{p_1,q_1\}$ to obtain the new one-handle specified by $\{p_1',q_1'\}$}
\end{figure}

We say that the points are arranged in a linear chain if for some
ordering on the components of $L$, $\{K_i\}_{i=1}^n$, if for each
$i=1,...,n-1$, we have that $p_i\in K_i$ and $q_i\in K_{i+1}$. First,
we claim that given any initial collection of points
$\{p_i',q_i'\}_{i=1}^n$, we can pass to a linear chain. To see this,
we let $K_1$ denote a circle with the minimal number of distinguished
points on it. Clearly  $K_1$ has only
one distinguished point, which we denote $p_1$. By
induction, we can find a sequence $K_1,...,K_i$ so that $K_1$ contains
only one distinguished point, $K_2,..,K_{i-1}$ all contain only two
distinguished points; indeed, $p_i\in K_i$ and $q_i\in K_{i+1}$ for
$i=1,...,i-1$.  If $K_{i}$ contains no other distinguished $p_i$, then
it follows from the connectivity hypothesis that $i=n$, and we have
arranged all our points in a linear chain. Otherwise, we fix some
$p_i\in K_i$. Its corresponding point $q_i$ lies on some other
component of $L$ which, according to the hypothesis, is not among
$K_1,...,K_{i}$. We denote that circle by $K_{i+1}$. Again, by a
sequence of handleslides over this one-handle, we can push all other
marked points from on $K_i$ to $K_{i+1}$.

Finally, any two linear chains can be connected by handleslides as
above. Specifically, choose any linear ordering of the components of
$L$, and consider the induced linear chain.  In the case where $K_1,
K_2,K_3 $ are the first three components in this linear ordering, the
handleslide described illustrated in Figure~\ref{fig:HandleSlide}
induces the permutation of $K_1$ and $K_2$. More generally, it is easy
to realize a permutation of two consecutive components as a
composition of two such handleslides.
\end{proof}

The above construction is more than a curiosity: it plays a central
role in the ``skein exact sequence'' described in Section~\ref{sec:Sequences}
below.

\subsection{Heegaard diagrams and knots}

\begin{defn}
For us, a {\em knot} will consist of a pair $(Y,K)$, where $Y$ is an
oriented three-manifold, and $K\subset Y$ is an embedded, oriented,
null-homologous circle.
\end{defn}

A knot $(Y,K)$ has a Heegaard diagram
$$(\Sigma,\alphas,\betas_0,\mu),$$ where here $\alphas$ is an
unordered $g$-tuple of pairwise disjoint attaching circles
$\alphas=\{\alpha_1,...,\alpha_g\}$, $\betas_0$ is a $(g-1)$-tuple of
pairwise disjoint attaching circles $\{\beta_2,...,\beta_g\}$, 
$\mu$ is an embedded, oriented circle in $\Sigma$ which is disjoint
from the $\betas_0$, and, of course $g$ is the genus of $\Sigma$.  This data is chosen so that
$(\Sigma,\alphas,\betas_0)$ specifies  the
knot-complement $Y-\nbd{K}$, i.e. if we attach disks along the
$\alphas$ and $\betas_0$, and then add a three-ball, we obtain the
knot-complement. Moreover, $\mu$ represents the ``meridian'' for the
knot in $Y$; thus, $(\Sigma,\alphas,\{\mu\}\cup\betas_0)$ is a
Heegaard diagram for $Y$.

\begin{defn}
A {\em marked Heegaard diagram for a knot $(Y,K)$} is a quintuple
$(\Sigma,\alphas,\betas_0,\mu,\Mark)$, where here $\Mark\in\mu\cap
(\Sigma-\alpha_1-...-\alpha_g)$.  
\end{defn}

The distinction between pointed
Heegaard diagrams (in the sense of~\cite{HolDisk}) and marked Heegaard
diagrams is that in the former, the distinguished point lies in the
complement of all the $g$-tuples, while in the latter, the
distinguished point lies on one of the curves.

A marked Heegaard diagram for a knot $K$ 
can also be used to construct 
a doubly-pointed Heegaard diagram for $Y$:

\begin{defn}
A {\em doubly-pointed Heegaard diagram for a three-manifold $Y$} is a
tuple $(\Sigma,\alphas,\betas,\StartPt,\EndPt)$, where
$(\Sigma,\alphas,\betas)$ is a Heegaard diagram for $Y$, and
$\StartPt$ and $\EndPt$ are a pair of distinct basepoints in $\Sigma$
which do not lie on any of the $\alphas$ or $\betas$.
\end{defn}

We extract a doubly-pointed Heegaard diagram for $Y$ from a Heegaard
diagram for $(Y,K)$ as follows.  Let
$(\Sigma,\alphas,\betas_0,\mu,\Mark)$ be a marked Heegaard diagram
for a knot complement $(Y,K)$. We write $\betas=\betas_0\cup\mu$.
Fix an arc $\delta$ which meets $\mu$
transversely in a single intersection point, which is the
basepoint $\Mark$, and which is disjoint from all the $\alphas$ and
$\betas_0$. Then, let $\FiltPt$ be the initial point of $\delta$, and
$\BasePt$ be its final point. 

Note that the ordering of the two points $\BasePt$ and $\FiltPt$ is
specified by the orientation of $K$. Specifically, choose a longitude
$\lambda$ for $K$, which we think of this as a curve in the Heegaard
surface. Clearly, the orientation on $K$ gives an orientation for
$\lambda$.  Now choose $\BasePt$ so that if $\delta$ is oriented as a path
from $\FiltPt$ to $\BasePt$, then we have an equality of algebraic
intersection numbers: $\#(\delta\cap\mu)=\#(\lambda\cap\mu)$ (for
either orientation of $\mu$).

Let $(\Sigma,\alphas,\betas,\StartPt,\EndPt)$ be a doubly-pointed
Heegaard diagram of genus $g$.  As in~\cite{HolDisk}, we consider the
$g$-fold symmetric product $\Sym^g(\Sigma)$, equipped with the pair of
tori
\begin{eqnarray*}
\Ta=\alpha_1\times...\times\alpha_g &{\text{and}}&
\Tb=\beta_1\times...\times\beta_g.
\end{eqnarray*}
Following notation from~\cite{HolDisk}, if $\x, \y\in\Ta\cap\Tb$, 
we consider spaces of homotopy classes of Whitney disks $\pi_2(\x,\y)$.
Letting 
$v\in\Sigma-\alpha_1-...-\alpha_g-\beta_1-...-\beta_g$, and choosing 
$\phi\in\pi_2(\x,\y)$, we let
$n_v(\phi)$ denote the oriented intersection number
$$n_v(\phi)=\#\phi^{-1}(\{v\}\times\Sym^{g-1}(\Sigma)).$$ A complex
structure on $\Sigma$ induces a complex structure on $\Sym^g(\Sigma)$
with the property that any Whitney disks $\phi$ which admits a
holomorphic representative has $n_v(\phi)\geq 0$. For a suitable small
perturbation of the holomorphic curve condition (which can be encoded
as a one-parameter family $J$ of almost-complex structures over
$\Sym^g(\Sigma)$), we obtain a
notion of pseudo-holomorphic disks, which form moduli spaces
satisfying suitably transverse Fredholm theory (along with the above
non-negativity hypothesis). We refer the reader
to~\cite{HolDisk} for a detailed discussion.

\subsection{Intersection points and $\SpinC$ structures}
\label{subsec:RelSpinC}

Recall that in~\cite{HolDisk}, for a Heegaard diagram for $Y$,
$(\Sigma,\alphas,\betas,\BasePt)$, we give a map
$$\spinc_\BasePt\colon \Ta\cap\Tb\longrightarrow \SpinC(Y).$$ When $Y$
is equipped with an oriented, null-homologous knot $K$, and
corresponding marked Heegaard diagram, this map has a refinement as
follows. Note first that there is a canonical zero-surgery $Y_0(K)$.
We define a map $$\relspinc_\Mark\colon \Ta\cap\Tb \longrightarrow
\SpinC(Y_0(K))$$ where here in the definition of $\Tb$ we use
$\betas=\mu\cup\betas_0$.  Sometimes, we denote the set
$\SpinC(Y_0(K))$ by $\RelSpinC(Y,K)$, and call it a set of {\em
relative $\SpinC$} structures for $(Y,K)$.  Note that
$$\RelSpinC(Y,K)\cong
\SpinC(Y)\times \Z.$$ The map from $\RelSpinC(Y,K)\longrightarrow
\SpinC(Y)$ is obtained by first restricting $\relspinct$ to $Y-K$, and
then uniquely extending it to $Y$ (note that the composition of these
two maps is the map $\spinc_\BasePt$ described earlier). Projection to
the second factor comes from evaluation $\OneHalf
\langle c_1(\relspinct),[{\widehat F}]\rangle$, where here ${\widehat
F}$ denotes a surface in $Y_0(K)$ obtained by capping off some fixed
Seifert surface $F$ for $K$. If $\relspinct\in\RelSpinC(Y,K)$ projects to 
$\spinc\in\SpinC(Y)$, we say that $\relspinct$ {\em extends} $\spinc$.

\begin{figure}
\mbox{\vbox{\epsfbox{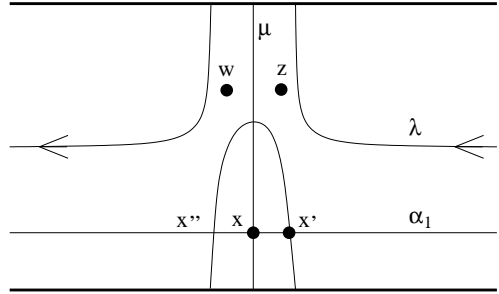}}}
\caption{\label{fig:WindConvention}
{\bf{Winding once along the meridian.}}  We have pictured the winding
procedure to illustrate the relative $\SpinC$ structure corresponding
to an intersection point. Note that the pictured region (all of which
takes place in a cylinder -- i.e. the two dark curves are to be
identified) is only a portion of the Heegaard diagram: the tuples $\x$
and $\x'$ in principle contain additional intersection points not
pictured here. The orientation of $\lambda$ induced from $K$
is illustrated  by the arrow.}
\end{figure}

To define $\relspinc_\Mark$, we replace the meridian with a
longitude $\lambda$ for the knot, chosen to wind once along the
meridian, never crossing the marked point $\Mark$,
so that each
intersection point $\x$ has a pair of closest points $\x'$ and $\x''$.
Let $(\Sigma,\alphas,\gammas,w)$
denote the corresponding Heegaard diagram for $Y_0(K)$
(i.e. $\gammas=\lambda\cup\betas_0$).  We then define
$\relspinc_{\Mark}(\x)$ to be the $\SpinC$ structure over $Y_0(K)$ 
given by $\spinc'_\BasePt(\x')=\spinc'_\BasePt(\x'')$, where here
$\spinc'_\BasePt\colon \Ta\cap\Tc \longrightarrow \SpinC(Y_0(K))$
denotes the usual map from intersection points to $\SpinC$ structures
(which we distinguish here from the analogous map $\spinc_\BasePt$ for $Y$).
Of course, the points  $\BasePt$ and $\FiltPt$ lie in the same
component of $\Sigma-\alpha_1-...-\alpha_g-\gamma_1-...-\gamma_g$, so
$\spinc'_\BasePt=\spinc'_\FiltPt$. See Figure~\ref{fig:WindConvention}
for an illustration.

Recall (c.f.~\cite{HolDisk}) that we can express the evaluation of
$\langle c_1(\relspinc_m(\x)),[{\widehat F}]\rangle$ in terms of data
on the Heegaard diagram. More precisely, each two-dimensional homology
class for $Y_0$ has a corresponding ``periodic domain'' $P$ in
$\Sigma$, i.e. a chain in $\Sigma$ whose local multiplicity at $w$ is
zero, and which bounds curves among the $\alphas$ and $\gammas$. Let
$\y\in\Ta\cap\Tc$ be an intersection point, let $\chi(P)$ denote the
``Euler measure'' of $P$ (which is simply the Euler characteristic of
$P$, if all its local multiplicities are zero or one), and let
${\overline n}_{\y}(P)$ denote the sum of the local multiplicities of
$P$ at the points $y_i$ comprising the $g$-tuple $\y$ (taken with a
suitable fraction if it lies on the boundary of $P$, as defined in
Section~\ref{HolDiskTwo:sec:Adjunction} of~\cite{HolDiskTwo};
compare also Equation~\eqref{eq:DefMult} below). Then,
it is shown in Proposition~\ref{HolDiskTwo:prop:MeasureCalc}
of~\cite{HolDiskTwo} that
\begin{equation}
\label{eq:MeasureCalc}
\langle c_1(\spinc'_\BasePt(\y)),[{\widehat F}]\rangle = 
\chi(P)+2{\Mass}_{\y}(P),
\end{equation}
provided that $P$ is the periodic domain representing ${\widehat F}$.

\begin{lemma}
\label{lemma:RelSpinCThree}
Let $(\Sigma,\alphas,\betas_0,\mu,\Mark)$ be a marked Heegaard diagram
for an oriented knot. Given $\x,\y\in\Ta\cap \Tb$, for any
$\phi\in\pi_2(\x,\y)$, we have that
$$\relspinc_{\Mark}(\x)-\relspinc_{\Mark}(\y)=(n_{\FiltPt}(\phi)-n_{\BasePt}(\phi))\cm\PD[\mu],$$
where here $[\mu]\in H_1(Y_0(K);\Z)$ is the homology class
obtained by thinking of 
the meridian of $K$ as a loop in $Y-\nbd{K}\subset Y_0(K)$, 
and oriented so that $\#(\mu\cap F)=1$. 
\end{lemma}

\begin{proof}
We begin with some preliminary observations.  Consider the natural
Heegaard triple for the the cobordism from $Y$ to $Y_0$
$(\Sigma,\alphas,\betas,\gammas,\BasePt)$, where $\gammas$ is obtained
as before by replacing the meridian $\mu$ in the $\betas$ by the
longitude $\lambda$. Note also that there is a second basepoint
$\FiltPt\in\Sigma$, which lies in the same component of 
$\Sigma-\alpha_1-...-\alpha_g-\gamma_1-...-\gamma_g$ as
$\BasePt$.
 
We claim that more generally if $\psi\in\pi_2(\x,\Theta,\y)$ is a Whitney
triangle (where, as usual, $\Theta$ is an
intersection point to $\Tb\cap\Tc$ representing the canonical
top-dimensional homology generator for $\HFleq(\#^{g-1}(S^2\#
S^1))$), then we have that
\begin{equation}
\label{eq:TriangleSpinC}
\spinc'_\BasePt(\y)
=\relspinc_\Mark(\x)+(n_\BasePt(\psi)-n_\FiltPt(\psi))\cm[\PD(\mu)].
\end{equation}

We see this as follows. First, fix $\x\in\Ta\cap\Tb$ and an
integer $k\in\Z$, and consider the set $S(\x,k)$ of all
$\y\in\Ta\cap\Tc$ for which there is a Whitney triangle
$\psi\in\pi_2(\x,\Theta,\y)$ with
$n_\BasePt(\psi)-n_\FiltPt(\psi)=k$. Clearly, this set is a single
$\SpinC$-equivalence class of intersection points on $Y_0(K)$.
Now it is an easy consequence of the definition of $\relspinc$ that
Equation~\eqref{eq:TriangleSpinC} holds in the case where $k=0$.

In view of the above remarks, to verify
Equation~\eqref{eq:TriangleSpinC}, it suffices to verify
it in the case where the triangle is
supported in the winding region. More precisely, fix an an integer
$k$, and consider the ``small'' triangle $\psi_k$ supported in the
winding region (after winding $\lambda$ sufficiently many times
in a neighborhood of $\mu$). For example,
when $k\leq 0$, after winding sufficiently many times, we can find a
triangle $\psi_k$ which lies on one side of $\mu$, so that
$n_\BasePt(\psi)=0$ and $n_\FiltPt(\psi)=k$. Now
$\psi_k\in\pi_2(\x,\Theta,\x_k')$, where the $\x_i'$ are the various
points for $\Ta\cap\Tc$ closest to $\x$ representing various $\SpinC$
structures $[S(\x,k)]$. It is an easy consequence now of
Equation~\eqref{eq:MeasureCalc} and our orientation conventions that
$$\langle c_1(\spinc'_\BasePt(\x_i)),[{\widehat F}]
\rangle -\langle c_1(\spinc'_\BasePt(\x_j)),[{\widehat F}] \rangle =
2(i-j).$$ This completes the verification of
Equation~\eqref{eq:TriangleSpinC}.

The equation stated in the lemma now follows from
Equation~\eqref{eq:TriangleSpinC} at once: let
$\psi\in\pi_2(\y,\Theta,\y')$ be a triangle with
$n_\FiltPt=n_\BasePt=0$, and then apply the formula for the juxtaposed
triangle $\psi'=\phi*\psi\in\pi_2(\x,\Theta,\y')$ with
$n_\BasePt(\psi')-n_{\FiltPt}(\psi')=n_\BasePt(\phi)-n_{\FiltPt}(\phi)$.
\end{proof}

\subsection{Filtered complexes}

Fix a partially ordered set $S$. An $S$-filtered group is a free
Abelian group $C$ generated freely by a distinguished set of
generators ${\mathfrak S}$ which admit a map $$\Filt\colon {\mathfrak
S}\longrightarrow S.$$ We write elements of $C$ as sums
$$\sum_{\sigma\in{\mathfrak S}} a_\sigma\cm \sigma,$$ where
$a_\sigma\in \Z$.  If
\begin{eqnarray*}
a=\sum_{\sigma\in{\mathfrak S}} a_\sigma\cm \sigma 
&{\text{and}}&
b= \sum_{\sigma\in{\mathfrak S}'} b_\sigma\cm\sigma 
\end{eqnarray*}
are elements of $S$-filtered groups $(C,\Filt,{\mathfrak S})$ and 
$(C',\Filt',{\mathfrak S}')$
respectively, then we write
$a\leq b$
if 
$$\max_{\{\sigma\in{\mathfrak S}\big| a_\sigma\neq 0\}} \Filt(\sigma)
\leq 
\min_{\{\sigma\in{\mathfrak S}'\big| b_\sigma\neq 0\}}\Filt'(\sigma).$$
A morphism of $S$-filtered groups 
$$\phi\colon (C,\Filt,{\mathfrak S})\longrightarrow 
(C',\Filt',{\mathfrak S}')$$
is a group homomorphism with the property that
$$\phi(a)\leq a$$
for all $a\in C$.

An {\em $S$-filtered chain complex} is an $S$-filtered group equipped
with a differential which is an $S$-filtered morphism; a morphism of
$S$-filtered chain complexes is a chain map which is also an
$S$-filtered morphism.  Let $T\subset S$ be a subset of $S$ with the
property that if $b\in T$, then all elements $a\in S$ with $a\leq b$
also are contained in $T$. If $T\subset S$ is such a subset, and if
$(C_*,\partial,\Filt)$ is an $S$-filtered complex, then $T$ gives rise
to a subcomplex of $C_*$, which transforms naturally under morphisms.

For example, recall that $\CFinf(Y,\spinct)$ defined in~\cite{HolDisk}
is generated by pairs $$[\x,i]\in (\Ta\cap\Tb)\times \Z$$ (where the
tori $\Ta$ and $\Tb$ lie in the $g$-fold symmetric product of a genus
$g$ Heegaard surface for $Y$). Indeed, the map $$\Filt[\x,i]=i$$
induces a natural filtration on this complex, and hence endows it with
a canonical subcomplex corresponding to the negative integers,
$\CFm(Y,\spinct)$, and a quotient complex $\CFp(Y,\spinct)$.

In the present paper, we will typically encounter complexes filtered by $\Z\times
\Z$, with the partial ordering $(i,j)\leq (i',j')$ if $i\leq i'$ and
$j\leq j'$.  These complexes have a number of naturally associated
subcomplexes. For instance, there is a subcomplex corresponding to the
quadrant $\{(i,j)| i\leq 0 {\text{~and~}} j\leq 0\}$; and there is
also a subcomplex corresponding to the union of three quadrants
$\{(i,j)|i\leq 0 {\text{~or~}} j\leq 0\}$.

\section{The knot filtration: definitions and basic properties}
\label{sec:DefKnot}

In this section, we begin by giving the definition of the knot
filtration.  In Subsection~\ref{subsec:TopInv} we prove its
topological invariance.  In Subsection~\ref{subsec:AbsGrading}, we
describe an absolute grading used on the knot filtration.  In
Subsection~\ref{subsec:Symm} we describe several symmetries of the
filtration: first under orientation reversal of the ambient
space, then under orientation reversal of the knot, and finally under
conjugation invariance of the underlying relative $\SpinC$ structures.

\subsection{Definition of the knot filtration}
\label{subsec:DefKnots}

Let $(\Sigma,\alphas,\betas,\BasePt,\FiltPt)$ be a doubly-pointed
Heegaard diagram, and $J$ be an allowed one-parameter family of
almost-complex structures over $\Sym^g(\Sigma)$. We can associate to
this data a $\Z^2$-filtered chain complex, following the constructions
of~\cite{HolDisk}. Specifically, we let
$\CFinf(\Sigma,\alphas,\betas,\BasePt,\FiltPt)$ be the free Abelian
group generated by triples $[\x,i,j]$ with $\x\in\Ta\cap\Tb$,
$i,j\in\Z$.  We endow this with the differential:
$$\partial^\infty[\x,i,j]=\sum_{\y\in\Ta\cap\Tb}
\sum_{\{\phi\in\pi_2(\x,\y)\big|\Mas(\phi)=1\}}
\#\left(\UnparModFlow(\phi)\right)[\y,i-n_{\BasePt}(\phi),j-n_{\FiltPt}(\phi)],$$
where here $\UnparModFlow(\phi)$ denotes the quotient of the moduli
space of $J$-holomorphic disks representing the homotopy type of
$\phi$, $\ModFlow(\phi)$, divided out by the natural action of $\R$
on this moduli space, and $\Mas(\phi)$ denotes the formal
dimension of
$\ModFlow(\phi)$. The signed count here
uses a coherent choice of orientations, as described in~\cite{HolDisk}.
Moreover, we can endow the chain complex with
the structure of a $\Z[U]$-module, by defining:
$$U\cm[\x,i,j]=[\x,i-1,j-1].$$
Of course, the filtration on the chain complex 
$$\Filt\colon (\Ta\cap\Tb)\times \Z^2\longrightarrow \Z^2$$
is given by
$$\Filt[\x,i,j]=(i,j),$$
where $\Z\times \Z$ is given the partial ordering $(i,j)\leq (i',j')$
when $i\leq i'$ and $j\leq j'$. 

The complex $\CFKinf$ naturally splits
into a sum of complexes. Specifically, generators $[\x,i,j]$ and
$[\y,\ell,m]$ lie in the same summand precisely when there is a
homotopy class $\phi\in\pi_2(\x,\y)$ with
\begin{eqnarray*}
n_w(\phi)=i-\ell&{\text{and}}&
n_z(\phi)=j-m.
\end{eqnarray*}

In the case studied here -- where the doubly-pointed Heegaard diagram
comes from from the marked Heegaard diagram of an oriented,
null-homologous knot $K$ -- this splitting can be interpreted in terms
of $\SpinC$ structures over $Y_0(K)$.  Specifically, fix a $\SpinC$
structure $\spinc$ over $Y$ and let
$\relspinct\in\RelSpinC(Y,K)=\SpinC(Y_0(K))$ be a $\SpinC$ structure
which extends it (c.f. Subsection~\ref{subsec:RelSpinC}).  consider
the subset $\CFKinf(Y,K,\relspinct)\subset
\CFinf(\Sigma,\alphas,\betas_0\cup\mu,\BasePt,\FiltPt)$ generated by
triples $[\x,i,j]$ with $\spinc_\BasePt(\x)=\spinc$ and 
\begin{equation}
\label{eq:DefOfHFK}
\relspinc_m(\x) + (i-j)PD[\mu]=\relspinct,
\end{equation}
where $[\mu]\in H_1(Y_0(K))$ is the homology class gotten by thinking
of the median $\mu$ as a closed curve in $Y_0(K)$.  According to
Lemma~\ref{lemma:RelSpinCThree}, this is a subcomplex.  Note that if
$\relspinct_1$ and $\relspinct_2$ represent the same $\SpinC$
structure over $Y$, then the complexes $\CFKinf(Y,\relspinct_1)$ and
$\CFKinf(Y,\relspinct_2)$ are isomorphic as chain complexes, and
indeed, the only difference is a shift in the $\Z\oplus\Z$-filtration.

Fix $\relspinct_0\in\RelSpinC(Y,K)$ and let $\spinc\in\SpinC(Y)$ be the compatible
element. We can view $\CFKinf(Y,K,\relspinct_0)$ as an extra $\Z$ filtration
on $\CFinf(Y,\spinc)$, by using the isomorphism
$$\Pi_1\colon \CFKinf(Y,K,\relspinct_0)\longrightarrow\CFinf(Y,\spinc)$$
given by
$$\Pi_1[\x,i,j]=[\x,i],$$
and declaring the extra $\Z$ filtration to be induced from the
projection 
of $[\x,i,j]$ to $j$.  

There is also a corresponding filtration structure on
$\CFa(Y,\spinc)$.  Specifically, let 
$$\CFK^{-,*}(Y,K,\relspinct_0)\subset \CFKinf(Y,K,\relspinct_0) $$
denote the subcomplex  corresponding to
$(i,j)$ with $i<0$, which in turn has a quotient complex
$\CFK^{+,*}(Y,K,\relspinct_0)$. Of course these complexes can be
thought of as induced filtrations on $\CFm(Y,\spinc)$ and
$\CFp(Y,\spinc)$ respectively. Consider the subcomplex
$\CFK^{0,*}(Y,K,\relspinct_0)\subset
\CFK^{+,*}(Y,K,\relspinct_0)$ generated by elements in the kernel of the induced
$U$-action; i.e. these generators all have $i=0$. This can be thought
of as a filtration on $\CFa(Y,\spinc)$.  The associated graded complex
for $\CFK^{0,*}(Y,K,\relspinct_0)$ admits an extra $\Z$-grading which,
of course, depended on our initial choice of $\relspinct_0$.
However, we can think of this object more invariantly as graded by
$\RelSpinC(Y,K)$. Specifically, for each 
$\relspinct\in \SpinC(Y_0(K))$, we have a complex
$\CFKa(Y,K,\relspinct)$ which is generated by intersection points $\x$
with $\relspinc_\Mark(\x)=\relspinct$, and its boundary operator counts only
homotopy classes $\phi\in\pi_2(\x,\y)$ with
$n_\BasePt(\phi)=n_\FiltPt(\phi)=0$. 
Then the associated graded graded complex for $\CFK^{0,*}(Y,K,\relspinct_0)$
is identified with
$$\bigoplus_{\{\relspinct\in\RelSpinC(Y,K)\big|\relspinct~{\text{extends}}~\spinc\}}
\CFKa(Y,K,\relspinct),$$
where the summand of the associated graded object for
$\CFK^{0,*}(Y,K,\relspinct_0)$ belonging to the integer $j$ corresponds
to the summand of the above complex belonging to 
$\relspinct=\relspinct_0-j\cm \PD[\mu]$.

When considering the complex $\CFKinf(Y,K,\relspinct)$ we always work
with a Heegaard diagram for the knot which is strongly
$\spinc$-admissible as a Heegaard diagram for $Y$ equipped with the 
$\SpinC$ structure $\spinc$ which is extended by $\relspinct$, while
when we consider its quotient complexes which are bounded below
(such as $\CFK^{+,*}$ and $\CFKa$), it suffices to consider only weakly
$\spinc$-admissible Heegaard diagrams.

\begin{theorem}
\label{thm:KnotInvariant}
Let $(Y,K)$ be an oriented knot, and fix a $\SpinC$
structure $\relspinct\in\RelSpinC(Y,K)$.
Then the filtered chain homotopy type of the chain
complex $\CFKinf(Y,K,\relspinct)$ is a topological invariant of the
oriented knot $K$ and the 
$\SpinC$ structure $\relspinct\in\RelSpinC(Y,K)$; i.e.  it is
independent of the choice of admissible, marked Heegaard diagram
$(\Sigma,\alphas,\betas_0,\mu,\Mark)$ used in its definition.
\end{theorem}

Indeed, we shall see that when the first Chern class of the
corresponding $\SpinC$ structure $\spinc$ over $Y$ is torsion, then
the above chain complex inherits in a natural way an absolute grading
from the absolute grading on $\CFinf(Y,\spinc)$
(c.f. Lemma~\ref{lemma:AbsGrade} below). We postpone the proof of the
above theorem till the next subsection, first stating some of its
immediate consequences.

\begin{cor}
\label{cor:KnotHomology}
The ``knot homology groups''
$\HFKa(Y,K,\relspinct)=H_*(\CFKa(Y,K,\relspinct))$ are topological
invariants of the knot $K\subset Y$ and $\relspinct\in\RelSpinC(Y,K)$.
\end{cor}

\vskip.2cm
\noindent{\bf{Theorem~\ref{thm:KnotInvariant}$\Rightarrow$ Corollary~\ref{cor:KnotHomology}.}}
It is clear from the above construction that the chain homotopy type
of $\CFKa(Y,K,\relspinct)$ is uniquely determined by the filtered chain
homotopy type of $\CFKinf(Y,K,\relspinct)$.
\qed
\vskip.2cm

In view of Proposition~\ref{prop:LinksToKnots}, passage from knot invariants
to link invariants is straightforward.

\begin{defn}
Let $L\subset Y$ be an $n$-component oriented link in $Y$, and let
$\kappa(L)$ denote the induced knot inside $\kappa(Y)\cong
Y\#^{n-1}(S^2\times S^1)$ as in Proposition~\ref{prop:LinksToKnots}. Then,
given $\relspinct\in\SpinC(\kappa(Y),\kappa(L))$ we define
$$\CFKinf(Y,L,\relspinct)=\CFKinf(\kappa(Y),\kappa(L),\relspinct).$$
\end{defn}

\begin{cor}
\label{cor:LinkInvariant}
Let $(Y,L)$ be an oriented link, and fix
$\relspinct\in\SpinC(\kappa(Y),\kappa(L)$. Then,
$\CFKinf(Y,L,\relspinct)$ is an invariant of the link $L\subset Y$.
\end{cor}

\vskip.2cm
\noindent{\bf{Theorem~\ref{thm:KnotInvariant}$\Rightarrow$ Corollary~\ref{cor:LinkInvariant}.}}
This implication is immediate, in view of
Proposition~\ref{prop:LinksToKnots}.
\qed
\vskip.2cm

As we see from above, the case of links is no more general than the
case of knots. Thus, we focus on the case of knots for most of the
present paper, returning to the disconnected case in
Section~\ref{sec:Links}.

\subsection{Proof of topological invariance of the knot homologies}
\label{subsec:TopInv}

The proof of Theorem~\ref{thm:KnotInvariant} is very closely modeled
on the proof of the topological invariance of $\HFc$ proved
in~\cite{HolDisk}.  Thus, we highlight here the difference between the
two results, leaving the reader to consult~\cite{HolDisk} for more
details.

Verification of the topological invariance rests on the following
description of how the invariant depends on the underlying Heegaard diagram:

\begin{prop}
\label{prop:HeegaardMoves}
If $(\Sigma,\alphas,\betas_0,\mu)$ and $(\Sigma',\alphas',\betas_0',\mu')$
represent the same knot complement, then we can pass from one to the other
by a sequence of the following types of moves (and their inverses):
\begin{itemize}
\item Handleslides and isotopies
amongst the $\alphas$ or the $\betas_0$
\item Isotopies of $\mu$
\item Handleslides of $\mu$ across some of the $\betas_0$
\item Stabilizations (introducing 
canceling pairs $\alpha_{g+1}$ and $\beta_{g+1}$ and increasing the
genus of $\Sigma$ by one).  \end{itemize} 
\end{prop}

\begin{proof}
This follows from standard Morse theory
(c.f. Lemma~\ref{HolDiskFour:lemma:LinkInvariance}
of~\cite{HolDiskFour}).
\end{proof}

For our proof, we adapt some of the Floer-homology constructions for
(singly-) pointed Heegaard diagrams described in~\cite{HolDisk} to the
doubly-pointed case.

Let $(\Sigma,\alphas,\betas,\BasePt,\FiltPt)$,
$(\Sigma,\betas,\gammas,\BasePt,\FiltPt)$ be a pair of doubly-pointed
Heegaard diagrams. Then, there is an induced (filtered) map:
$$F\colon \CFinf(\Sigma,\alphas,\betas,\BasePt,\FiltPt)\otimes 
\CFinf(\Sigma,\betas,\gammas,\BasePt,\FiltPt)\longrightarrow
\CFinf(\Sigma,\alphas,\gammas,\BasePt,\FiltPt)$$ defined as follows:
$$F([\x,i,j]\otimes [\y,\ell, m]) =\sum_{\{\psi\in
\pi_2(\x,\y,\w)\big|\Mas(\psi)=0\}}
\left(\#\ModFlow(\psi)\right)\cm [\w,i+\ell-n_\BasePt(\psi),j+m-n_\FiltPt(\psi)],$$
where $\#\ModFlow(\psi)$ denotes a signed count of pseudo-holomorphic
triangles.

\vskip.2cm
\noindent{\bf{Proof of Theorem~\ref{thm:KnotInvariant}}.}
First we observe that the filtered
chain homotopy type  $\CFKinf$  is independent
of the perturbation of complex structure $J$. This is a straightforward
modification of the corresponding discussion in~\cite{HolDisk}.
Similarly, independence of the complex under isotopies of the $\betas_0$ is
a straightforward modification of the corresponding invariance of $\HFc$:
for an isotopy consisting of a single pair creations, one considers the 
natural chain homotopy equivalence induced by an exact Hamiltonian isotopy
realizing the isotopy.

For independence of the groups under handleslides amongst the
$\betas_0$, follow the arguments of handleslide invariance of $\HFp$
as well. Now, we observe that if $\betas_0'$ is obtained from
$\betas_0$ by a handleslide, then the doubly-pointed Heegaard diagram
$(\Sigma,\mu\cup\betas_0,\mu'\cup\betas_0',\BasePt,\FiltPt)$ (where
$\mu'$ is a small exact Hamiltonian translate of $\mu$) has the
property that $\BasePt$ and $\FiltPt$ lie in the same connected
component of
$\Sigma-\mu-\beta_2-...-\beta_g-\mu'-\beta_2'-...\beta_g'$.  We
arrange for this Heegaard diagram to be strongly $\spinc_0$-admissible
as a pointed Heegaard diagram for $\#^{g}(S^2\times S^1)$ (using
either basepoint), where here $\spinc_0$ is the $\SpinC$ structure
with trivial first Chern class.  We let $\CFinf_{\delta}$ denote the
``diagonal'' summand of
$\CFinf(\Sigma,\mu\cup\betas_0,\mu'\cup\betas_0')$ generated by
$[\x,i,i]$ (where $[\x,i]$ represents $\spinc_0$).  It is easy to see
then that (c.f. Lemma~\ref{HolDiskOne:lemma:Tori1} of~\cite{HolDisk}),
$$\HFleq_\delta(\mu\cup\betas_0,\mu'\cup\betas_0',\BasePt,\FiltPt)
\cong \Z[U]\otimes_\Z \Wedge^* H_1(T^g).$$
Letting $\Theta$ represent the top-dimensional non-zero generator, we
use $$\x\mapsto F(\x\otimes\Theta)$$ to define the map associated to
handleslides. The proof that this is a filtered chain homotopy
equivalence now follows from a direct application of the methods from
Section~\ref{HolDiskOne:sec:HandleSlides} of~\cite{HolDisk}.  The same
remarks apply when one slides $\mu$ over one of the circles in
$\betas_0$. (Note that if one were to try to slide $\beta_2$ over
$\mu$, the arguments would break down in view of the fact that now
$\BasePt$ and $\FiltPt$ lie in different components of the complements
of the curves. This move, however, is not required for proving the
topological invariance of the knot homology, according to
Proposition~\ref{prop:HeegaardMoves} above.)

The above argument shows that the chain homotopy type of the
doubly-filtered chain complex is invariant under handleslides. It is
not difficult to see that the map also preserves the splitting
according to $\relspinct\in\RelSpinC(Y,K)$
(This is shown in a somewhat more general
context in Proposition~\ref{prop:RespectFiltration} below.)

Constructing such an isomorphism for isotopies and handleslides along
the $\alphas$ proceeds in an analogous manner: one shows that the
filtered chain homotopy type remains invariant under isotopies and
handleslides which do not cross the arc $\delta$ connecting $\BasePt$
and $\FiltPt$. Observe also that an isotopy of the $\alphas$ which does
cross this arc can be realized by a sequence of isotopies and
handleslides amongst the $\alphas$ which do not cross the arc
(c.f. Proposition~\ref{HolDiskOne:prop:PointedHeegaardMoves}
of~\cite{HolDisk}).

Stabilization invariance, too, follows directly as in~\cite{HolDisk}.
\qed
\vskip.2cm

\subsection{Absolute gradings}
\label{subsec:AbsGrading}

Fix a $\SpinC$ structure $\spinc\in\SpinC(Y)$, whose first Chern class
is torsion. Recall (c.f.~\cite{HolDiskFour}) that in this case, we defined
an absolute $\Q$-grading on $\CFinf(Y,\spinc)$ which
enjoys the following formula: if $W$ is a cobordism from $Y_1$ to $Y_2$,
$\spincf\in\SpinC(W)$ is a  $\SpinC$ structure
whose restrictions $\spinc_1$ and $\spinc_2$ to $Y_1$ and $Y_2$ respectively
are both
torsion, then the induced map $F_{W,\spincf}$ shifts degree by
\begin{equation}
\label{eq:DimensionFormula}
\frac{c_1(\spincf)^2 - 2\chi(W)-3\sigma(W)}{4},
\end{equation}
where $\chi(W)$ and $\sigma(W)$ denote the Euler characteristic and signature
of $W$ respectively.

The torsion hypothesis also gives us a canonical
$\relspinct_0\in\RelSpinC(Y,K)$ extending $\spinc$ which satisfies
\begin{equation}
\label{eq:DefCan}
\langle c_1(\relspinct_0),[{\widehat F}]\rangle = 0
\end{equation} for any choice
of Seifert surface $F$ for $K$ (indeed, the condition
that $\relspinct_0$ is
independent of the choice of Seifert surface is equivalent to the
condition that $c_1(\spinc)$ is torsion). We can endow
$\CFKinf(Y,K,\relspinct_0)$ with the absolute grading induced from the
map $$\Pi_1\colon \CFKinf(Y,K,\relspinct_0)\longrightarrow
\CFinf(Y,\spinc)$$ given by $\Pi_1[\x,i,j]=[\x,i]$.

Since $\relspinct_0\in\RelSpinC(Y,K)$ is uniquely determined by
$\spinc\in\SpinC(Y)$, we typically write
$\CFKinf(Y,K,\spinc)$ for $\CFKinf(Y,K,\relspinct_0)$.

\begin{lemma}
\label{lemma:AbsGrade}
Let $Y$ be an oriented three-manifold $K\subset Y$ be a knot, and fix
a $\SpinC$ structure $\spinc\in\SpinC(Y)$. Then, there is a convergent
spectral sequence of relatively graded groups whose $E^1$ term is
$$\bigoplus_{\{\relspinct\in\RelSpinC(Y,K)\big|
\relspinct~{\text{extends}}~\spinc\}}\HFKa(Y,K,\relspinct),$$
and whose
$E^\infty$-term is $\HFa(Y,\spinc)$. Moreover, when $c_1(\spinc)$
is torsion, the spectral sequence respects absolute gradings.
\end{lemma}

\begin{proof}
This is simply the Leray spectral sequence associated to the
filtration of $\CFa(Y,\spinc)$ given by $\CFK^{0,*}(Y,\relspinct_0)$
(where $\relspinct\in\RelSpinC(Y,K)$ is any choice of $\SpinC$ structure
which extends $\spinc\in\SpinC(Y)$). The statement in the absolutely graded
case follows from our definition of the absolute grading on $\HFKa$.
\end{proof}

Of course, we have analogous spectral sequences converging to
$\HFp(Y,\spinc)$, $\HFm(Y,\spinc)$, and $\HFinf(Y,\spinc)$ whose $E_1$
terms are derived from $\CFKinf(Y,K)$ in a straightforward way.

\subsection{Notational shorthand}
\label{subsec:Shorthand}
We make several further simplifications when considering
knots $K$ in $S^3$. We write simply $\CFKinf(S^3,K)$ for the
absolutely $\Z$-graded complex $\CFKinf(S^3,K,\relspinct_0)$, where here
$\relspinct_0\in\RelSpinC(S^3,K)$ is uniquely characterized by the property that
$c_1(\relspinct_0)=0$. Moreover, given an integer $n\in\Z$, write
$\HFKa(S^3,K,n)$ for $\HFKa(S^3,K,\relspinct_n)$, where $\relspinct_n\in\RelSpinC(S^3,K)$ is characterized by the property that $\langle
c_1(\relspinct_n),[{\widehat F}]\rangle = 2n$, where $F$ is a
compatible Seifert surface for the oriented knot $K$.

In this setting,  we also write
$\HFp(S^3_0(K),n)$ to denote $\HFp(S^3_0(K),\spincrel_n)$ (where we think of
$\relspinc_n$ as a $\SpinC$ structure over $Y_0$.

\subsection{Symmetries}
\label{subsec:Symm}

\begin{prop}
\label{prop:IndepOrient}
The dependence on the orientation of $Y$ is given by
$$\HFKa_*(Y,K,\relspinc)\cong \HFKa^*(-Y,K,\relspinc)$$ (where here
the left-hand-side denotes homology, and the right-hand-side denotes
cohomology). Moreover, if $\relspinc$ extends a torsion $\SpinC$ structure
over $Y$, then 
$$\HFKa_d(Y,K,\relspinc)\cong \HFKa^{-d}(-Y,K,\relspinc).$$
\end{prop}

\begin{proof}
This follows from the behaviour of $\HFa(Y)$ under
orientation reversal. Specifically, if $(\Sigma,\alphas,\betas,w,z)$
describes the knot $K$ in $Y$, then $(-\Sigma,\alphas,\betas,w,z)$
describes the knot $K$ in $-Y$. The identification of the 
complexes now follows as in~\cite{HolDiskTwo}.
\end{proof}

\begin{prop}
\label{prop:OrientKnot}
Fix a torsion $\SpinC$ structure $\spinc$ over $Y$.  For each
extension $\spincrel\in\RelSpinC(Y,K)$ of $\spinc$, we have that
$$\HFKa_d(Y,K,\spincrel)\cong \HFKa_{d-2m}(Y,-K,\spincrel),$$ where
$$m = \OneHalf \langle c_1(\spincrel),[{\widehat F}]\rangle$$
is calculated using 
the surface ${\widehat F}$ obtained by capping off a Seifert surface $F$
for the oriented knot $K$.
\end{prop}

\begin{proof}
Observe that for the
chain complex underlying $\CFKinf(Y,K,\relspinct_0)$, the roles of $i$ and $j$
can be interchanged; i.e. the same complex can be viewed at the same time
as $\CFKinf(Y,-K,\relspinct_0)$.  (Recall that
$\relspinct_0$ is chosen to satisfy Equation~\eqref{eq:DefCan}.) 
Correspondingly, we have a projection map
$$\Pi_2\colon \CFKinf(Y,K,\relspinct_0)
\cong \CFKinf(Y,-K,\relspinct_0)\longrightarrow \CFinfty(Y,\spinc)$$
given by $\Pi_2[\x,i,j]=[\x,j]$. (By all rights, the complex
$\CFinfty(Y,\spinc)$ ought to include the base point in its notation:
in this case we use $\FiltPt$ rather than $\BasePt$).  This
projection, too induces an absolute grading on
$\CFKinf(Y,K,\relspinct_0)$, but in fact it coincides with the absolute
grading induced from $\Pi_1$, since there is there is a
grading-preserving homotopy equivalence from $\CFinfty(Y,\spinc)$, as
defined using the base point
$\BasePt$, to the corresponding complex as defined using
$\FiltPt$.

Now, the stated isomorphism is induced by the map
$$U^m\colon \CFKa(Y,K,\spincrel)\cong
\CFK^{0,m}(Y,K,\relspinct_0)
\longrightarrow \CFK^{-m,0}(Y,K,\relspinct_0)\cong
\CFKa(Y,-K,\spincrel).$$
\end{proof}

There is a conjugation symmetry
$$J\colon \RelSpinC(Y,K)\longrightarrow \RelSpinC(Y,K)$$
induced from the usual conjugation symmetry on $\SpinC(Y_0(K))$.

\begin{prop}
\label{prop:JInvarianceGen}
For each $\relspinc\in\RelSpinC(Y,K)$, there is an identification
$$\CFKinf(Y,K,\relspinc)\cong \CFKinf(Y,-K,J\relspinc).$$
\end{prop}

\begin{proof}
First notice that in our definition the Heegaard diagram for a knot,
there is an asymmetry in the role of the $\alphas$ and the $\betas$;
specifically, our choice of meridian was chosen to be among the
$\betas$.  This was, of course, unnecessary. We could just as well
have chosen our meridian to be an $\alphas$-circle, repeated the above
definition, and obtained another knot filtration which is also a knot
invariant.  In fact, these two filtrations coincide, provided that
we switch the roles of $\BasePt$ and $\FiltPt$. Indeed, we can
always choose a Heegaard diagram for the knot $K$ for which there are
representatives for $K$ in both handlebodies, with both
$\alpha_1$ and $\beta_1$ representing meridians for the two
representatives, reversing the roles of  $\BasePt$ and $\FiltPt$. Such
a diagram can be obtained by  stabilizing a 
doubly-pointed Heegaard diagram for which the attaching circle  $\beta_1$
is a meridian for $K$. Specifically, we stabilize our original Heegaard surface
by attaching a handle with feet near $\BasePt$ and $\FiltPt$ respectively,
the circle $\alpha_1$ is supported inside the newly attached handle, and 
we introduce a canceling circle $\beta_2$ which runs along a longitude for $K$
except near $\beta_1$, which it avoids by running through the one-handle.
This is illustrated in Figure~\ref{fig:JInvariance}.

\begin{figure}
\mbox{\vbox{\epsfbox{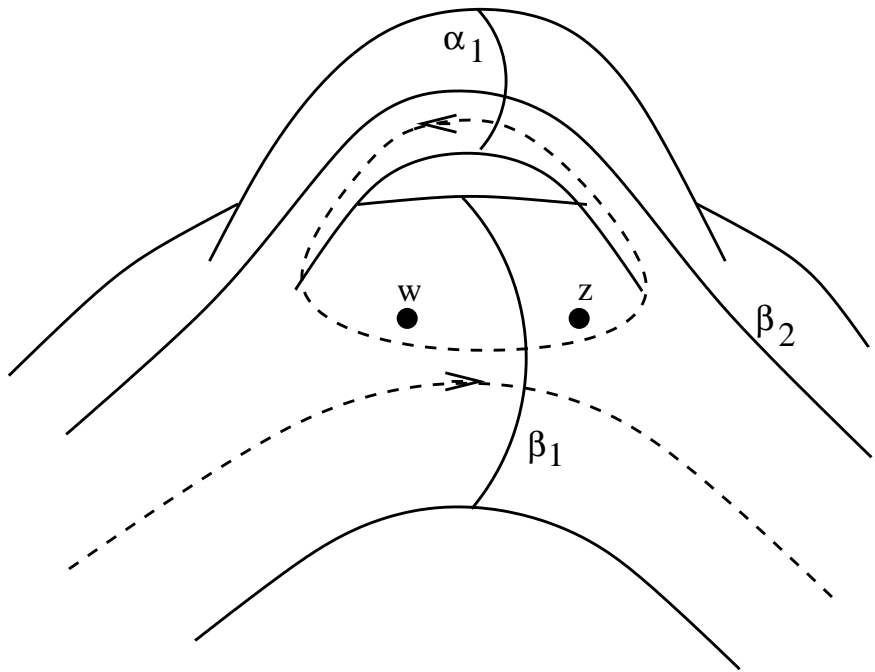}}}
\caption{\label{fig:JInvariance}
{\bf{Moving the meridian.}}  This illustrates the symmetry between the
roles of $\alpha_1$ and $\beta_1$ as meridians for a given knot
$K$. This knot $K$ is isotopic to either of the two (unlabeled)
dotted circles in $\Sigma$, and hence can be pushed into either
handlebody so that either $\alpha_1$ or $\beta_1$ is a meridian.  Note
that there may be additional $\alpha$ and $\beta$-circles entering the
picture, but it is important that they do not cross the dotted
circles, and do not separate $w$ and $z$ from these dotted circles.}
\end{figure}

With these remarks in place, we realize the conjugation symmetry
by the natural identification of the complexes
$(\Sigma,\alphas,\betas,\BasePt,\FiltPt)$ and
$(-\Sigma,\betas,\alphas,\FiltPt,\BasePt)$,
following~\cite{HolDisk}. It is easy to see that
this identification induces an isomorphism of $\CFKinf(Y,K,\spincrel)$ with
$\CFKinf(Y,K,J\spincrel)$.
\end{proof}

Fix a torsion $\SpinC$
structure $\spinc$ over $Y$, and let 
$\CFKinf(Y,K,\spinc)'$ denote the $\Z\oplus\Z$-filtered chain complex
$\CFKinf(Y,K,\spinc)$, endowed with the filtration 
$$\Filt[\x,i,j]=(j,i).$$
The above proposition gives  a degree-preserving
$\Z\oplus\Z$-filtered chain homotopy equivalence 
$$\Phi\colon \CFKinf(Y,K,\spinc)\longrightarrow  \CFKinf(Y,K,\spinc)'.$$

This gives the following corollary, which we state for the knot Floer homology:

\begin{prop}
\label{prop:JInvariance}
Let $Y$ be a three-manifold equipped with a torsion $\SpinC$ structure,
then
$$\HFKa_d(Y,K,\spincrel)\cong\HFKa_{d-2m}(Y,K,J\spincrel),$$
where $2m=\langle c_1(\spincrel),[{\widehat F}]\rangle$.
\end{prop}

\begin{proof}
Combine the results of Propositions~\ref{prop:OrientKnot} and
\ref{prop:JInvarianceGen}.
\end{proof}

\section{Relationship with three-manifold invariants}
\label{sec:Relationship}

We consider the following two $\Z\oplus\Z$-filtered, $\Z[U]$-subcomplexes 
of $\CFKinf(Y,K)$:
\begin{itemize}
\item $\CFKm(Y,K)$, which  denotes the subcomplex of  $\CFK(Y,K)$ generated by
$[\x,i,j]$ with $\max(i,j)<0$,
\item $\bCFKm(Y,K)$ (``big'' $\CFKm$), which 
denotes the subcomplex generated by $[\x,i,j]$ with $\min(i,j)<0$.
\end{itemize}
These have quotient complexes denoted 
\begin{eqnarray*}
\bCFKp(Y,K)&=&\CFKinf(Y,K)/\CFKm(Y,K) \\
\CFKp(Y,K)&=&\CFKinf(Y,K)/\bCFKm(Y,K).
\end{eqnarray*}
Our aim is to identify the homologies of $\CFKp(Y,K)$
(resp. $\bCFKp(Y,K)$) with $\HFp(Y_{-p}(K))$ (resp. $\HFp(Y_p(K))$)
for sufficiently large integer surgery coefficients
$p$, to be made precise shortly.

Fix a Seifert surface $F$ for $K\subset Y$, compatible with the
orientation of $K$. We can cap this off to obtain a surface ${\widehat
F}\subset Y_0(K)$. For each $m\in\Z$, let
$$\CFKinf(Y,K,F,m)=\bigoplus_{\{\relspinct\in\RelSpinC(Y,K)\big|
\langle c_1(\relspinct),[{\widehat F}]\rangle =2m\}}
\CFKinf(Y,K,\relspinct),$$
Similarly, for each integer $p$, we let $S$ denote the surface in $W_{-p}(K)$
obtained by closing off $F$, to obtain a surface of square $-p$.
Now, given $[m]\in\Zmod{p}$, let ${\mathfrak T}(F,p,[m])\subset \SpinC(Y_p(K))$ 
denote the set of $\SpinC$ structures which can be extended
over $\SpinC(W_{-p}(K))$ to give a $\SpinC$ structure
$\spincf$ with
$$\langle c_1(\spincf),[S]\rangle + p \equiv 2m\pmod{2p}.$$
Correspondingly, we let
$$\CFinf(Y_p(K),F,[m])=\bigoplus_{\spinct\in{\mathfrak T}(F,p,m)}
\CFinf(Y_p(K),F,\spinct).$$
When $F$ is understood from the context, then we drop it from the notation.

To relate these two complexes, we will use the following map.
Fix a positive integer $p$, and let
$(\Sigma,\alphas,\betas,\gammas,m)$ be a marked Heegaard triple for the
cobordism from $Y$ to $Y_{-p}$.  Consider the
chain map $$\Phi\colon \CFKinf(Y,K)\longrightarrow
\CFinf(Y_{-p})$$ defined by 
$$\Phi[\x,i,j]=\sum_{\y\in\Ta\cap\Tc}
\sum_{\{\psi\in\pi_2(\x,\Theta,\y)\big| 
n_\BasePt(\psi)-n_\FiltPt(\psi)=i-j,
\Mas(\psi)=0\}}
\left(\#\ModFlow(\psi)\right)\cm [\y,i-n_\BasePt(\psi)].$$ 

\begin{theorem}
\label{thm:LargeNegSurgeries}
Let $(Y,K)$ be an oriented knot, and fix a compatible Seifert surface $F$
for $K$. Then, for each integer $m\in\Z$, we have an integer
integer $N=N(m)$ so that for all integers $p\geq N$,
there is a Heegaard diagram for which the map $\Phi$ above 
induces isomorphisms
of chain complexes which fit into the following commutative diagram:
$$
\begin{CD}
0@>>>\bCFKm(Y,K,m)@>>> \CFKinf(Y,K,m) @>>>\CFKp(Y,K,m)@>>> 0 \\ &&
@V{\bPhim}VV @V{\Phi^\infty}VV @VV{\Phip}V \\
0@>>>\CFm(Y_{-p}(K),[m])@>>> \CFinf(Y_{-p}(K),[m])
@>>>\CFp(Y_{-p}(K),[m])@>>> 0.
\end{CD}
$$   Similarly, if we let $\CFKsa(Y,K,m)\subset 
\CFKp(Y,K,m)$ be the subset
generated by $[\x,i,j]$ with $\min(i,j)=0$, then $\Phi$ also induces
isomorphisms in the following diagram: $$
\begin{CD}
0@>>>\CFKsa(Y,K,m)@>>> \CFKp(Y,K,m) @>{U}>>\CFKp(Y,K,m)@>>> 0 \\ &&
@V{{\widehat\Phi}}VV @V{\Phip}VV @VV{\Phip}V \\ 0@>>>\CFa(Y_{-p},[m])
@>>> \CFp(Y_{-p},[m]) @>{U}>> \CFp(Y_{-p},[m]) @>>> 0.
\end{CD}
$$
\end{theorem}

\begin{figure}
\mbox{\vbox{\epsfbox{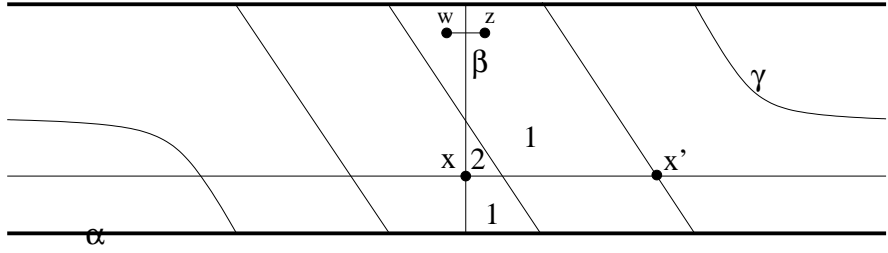}}}
\caption{\label{fig:LargeSurgeries}
{\bf{Illustration of Theorem~\ref{thm:LargeNegSurgeries}.}}
The integers denote (non-zero) local multiplicities of the ``small triangle'' $\Phi_0$
connecting $x$ and $x'$. This picture is taking place in a cylindrical 
region in $\Sigma$: the top dark line is identified with the bottom.}
\end{figure}

\begin{proof} 
First, note that $\Phi$ maps $\bCFKm(Y,K)$ into $\CFm(Y_{-p})$:
if $\psi$ has a holomorphic representative, then $n_\BasePt(\psi)\geq
0$ and $n_\FiltPt(\psi)\geq 0$, so if $[\x,i,j]\in\bCFKm(Y,K)$, then
either $i<0$, so clearly, $i-n_\BasePt(\psi)<0$; or $i\geq 0$ and
$j<0$, so $$i-n_\BasePt(\psi)=j-n_\FiltPt(\psi)<0.$$

Fix a marked Heegaard diagram $(\Sigma,\alphas,\betas,\Mark)$ for
$(Y,K)$, where $\Mark\in\beta_1$, so that $\beta_1$ is the
meridian. Fix an annular neighborhood of $\beta_1$ and a longitude
$\lambda$ for $K$. We will consider Heegaard triples
$(\Sigma,\alphas,\betas,\gammas,\BasePt)$ so that the $\gamma_i$ for
$i\geq 2$ are all exact Hamiltonian translates of the corresponding
$\beta_i$, and $\gamma_1$ is a small (embedded) perturbation of
the juxtaposition $p\beta_1+\lambda$. We will situate $\gamma_1$ so 
that $\beta_1$ intersects it in one point towards the middle of this winding
region (c.f. Figure~\ref{fig:LargeSurgeries}).

We say that an intersection $\x'\in\Ta\cap\Tb$ is supported in the winding
region if the coordinate on $\x'$ in $\gamma_1$ is supported there.
In this case, there is a uniquely determined intersection 
point $\x\in\Ta\cap\Tb$ which is closest to $\x'$. Moreover, there is 
a canonical ``small triangle'' $\psi_0\in\pi_2(\x,\Theta,\x')$
supported in the winding region. Conversely, given any $\x\in\Ta\cap\Tb$,
there are $p$ distinct intersection points for $\Ta\cap\Tc$, representing
the $p$ $\SpinC$ structures over $Y_{-p}(K)$ which are cobordant to
the $\SpinC$ structure represented by $\x$.

We claim that 
\begin{equation}
\label{eq:SpinCFormula}
\langle c_1(\spincrel(\x)),[{\widehat F}]\rangle +2(n_{\BasePt}(\psi)-n_{\FiltPt}(\psi))
= \langle c_1(\spinc_{\BasePt}(\psi)),[S]\rangle-p.
\end{equation}
Fix a $\SpinC$ structure $\spinc\in\SpinC(Y)$, and let 
$\x_1$ be some representative for $\spinc$. 
For $\psi=\psi_0(\x_1,\Theta,\x_1')\in\pi_2(\x,\Theta,\x')$
with $n_{\BasePt}(\psi)=n_{\FiltPt}(\psi)=0$, 
this is a straightforward application of
the first Chern class formula from~\cite{HolDiskFour}. 
Now, any $\psi$ which induces the same $\SpinC$ structure has the form
$\psi=\phi*\psi_0*\phi'$ with $\phi\in\pi_2(\x_2,\x_1)$ (for
$\Ta\cap\Tb$, $\phi'\in\pi_2(\x_1',\x_2')$. It follows from
Lemma~\ref{lemma:RelSpinCThree} Equation~\eqref{eq:SpinCFormula}
is unaffected by the juxtaposition
of $\phi$, while it is clear that it is unaffected by the juxtaposition
of $\phi'$. Finally, it is straightforward to see
that the equation is preserved by adding a  triply-periodic domain belonging
to $S$. This completes
the verification of Equation~\eqref{eq:SpinCFormula} for all Whitney
triangles $\psi$.

It follows from Equation~\eqref{eq:SpinCFormula}, that 
the restriction of $\Phi$ to $\CFKinf(Y,K,m)$
maps to $\CFp(Y_{-p}(K),[m])$.

Fix an integer $m\in\Z$. We claim that if $p$ is sufficiently large,
then that all the intersection points between $\Ta$ and $\Tc$ which
represent any $\SpinC$ structure $\spinc\in{\mathfrak T}(F,p,[m])$ are
supported in the winding region. This can be seen, for example, from
Equation~\eqref{eq:SpinCFormula}. In this case, consider the map
$$\Phi_0\colon \CFKinf(Y,K,m)\longrightarrow \CFinf(Y_{-p}(K),[m])$$
which carries $[\x,i,j]$ to $[\x',i-n_{\BasePt}(\psi_0)]$.  By the
Riemann mapping theorem, we have that $\#\ModFlow(\psi_0)=1$, so we
can think of $\Phi_0$ as the summand of $\Phi$ corresponding to the
homotopy class $\psi_0$.  Moreover, it is easy to see that $\Phi_0$
induces isomorphisms of groups
\begin{eqnarray*}
\CFKinf(Y,K,m)&\cong& \CFinf(Y_{-p}(K),[m]) \\
\bCFKm(Y,K,m)&\cong& \CFm(Y_{-p}(K),[m]) \\
\CFKp(Y,K,m)&\cong& \CFp(Y_{-p}(K),[m]) \\
\CFKsa(Y,K,m)&\cong& \CFa(Y_{-p}(K),[m]).
\end{eqnarray*}

Now, if the winding region has sufficiently small area relative to the
areas of the regions in $\Sigma-\alphas-\betas$ (divided by $n$), then
it is clear that $$\Phi=\Phi_0+{\text{lower order}},$$ with respect to
the energy filtrations on $\CFK(Y,K,s)$ and $\CF(Y_n(K),[s])$
(induced by areas of domains). Thus, by elementary algebra, the map
$\Phi$, which we know is a chain map, is also an isomorphism
(carrying $\bCFKm(Y,K,s)$ to $\CFm(Y,[s])$).
\end{proof}

For convenience, we state a version of the above result in the case
where we restrict to $\SpinC$ structures over $Y$ whose first Chern
class is torsion. In this case, recall that we defined
$\CFKinf(Y,K,\spinc)$ to be $\CFKinf(Y,K,\relspinct_0)$, where
$\relspinct_0\in\RelSpinC(Y,K)$ restricts $\spinc$ and also satisfies
$\langle c_1(\relspinct_0),[{\widehat F}]\rangle =0$, a notion which
now is independent of the choice of Seifert surface for $K$.  We let
$$\CFKinf(Y,K,\spinc)\{i\geq 0~\text{or}~j\geq -m\}\subset
\CFKinf(Y,K,\spinc)$$ denote the induced (quotient) complex generated
by $[\x,i,j]\in\CFKinf(Y,K,\spinc)$, subject to the stated
constraint $i\geq 0~\text{or}~j\geq -m$.
Moreover, let $[\spinc,m]\in\SpinC(Y_{-p})$
denote the $\SpinC$ structure over $Y_{-p}$ which is cobordant to $\spinc$,
via a $\SpinC$ structure $\spincf\in\SpinC(W_{-p})$ with 
$$\langle c_1(\spincf),[S]\rangle = 2m+p.$$

\begin{cor}
\label{cor:TorsionCase}
Let $K\subset Y$ be a knot, and fix $\spinc\in\SpinC(Y)$
with
whose first Chern class is torsion. Then, there is 
an integer $N$ with the property that for all $p\geq N$,
we have that
\begin{eqnarray*}
\HFp_{\ell}(Y_{-p}(K),[\spinc,m])\cong 
\left\{\begin{array}{ll}
H_k(\CFKinf(Y,K,\spinc)\{i\geq 0~\text{and}~j\geq -m\})
&{\text{if $|m|\leq g$}} \\
\HFp_k(Y) & {\text{otherwise,}}
\end{array}\right.
\end{eqnarray*}
where 
$$
\ell=k+\left(\frac{p-(2m+p)^2}{4p}\right).
$$ 
\end{cor}

\begin{proof}
If $|m|>g$, this follows readily from the adjunction inequality for
$\HFp(Y_0)$, together with long exact sequence for integral surgeries
(c.f.~\cite{HolDiskTwo}) (provided that $p>g$).

The cases where $|m|<g$, this follows from the above theorem.
Specifically, observe that if $\relspinct\in\RelSpinC(Y,K)$ extends $\spinc$,
then we have an isomorphism of chain complexes
$$\CFKinf(Y,K,\relspinct)\cong
\CFKinf(Y,K,\relspinct_0)$$ defined by $$[\x,i,j]\mapsto [\x,
i,j+\OneHalf \langle c_1(\relspinct),[{\widehat F}]\rangle].$$ This
preserves the absolute $\Z$ grading as defined in
Subsection~\ref{subsec:AbsGrading}, shifting the filtration in the
obvious way. Thus, the isomorphism in its relatively graded form is a
direct consequence of Theorem~\ref{thm:LargeNegSurgeries}. For the
graded statement, we recall that the map $\Phi$ inducing the
isomorphism is a filtered version of the map from $\CFinf(Y)$ to
$\CFinf(Y_{-p})$ induced by counting holomorphic triangles
(c.f.~\cite{HolDiskFour}), associated to the $\SpinC$ structure over
$W_{-p}(K)$ with $$\langle c_1(\spincf),[S]\rangle = 2m+p,$$ according
to Equation~\eqref{eq:SpinCFormula}. Thus, the shift in dimension
follows from Equation~\eqref{eq:DimensionFormula}.
\end{proof}

\begin{remark}
It follows from the adjunction inequality for $\HFp(Y_0)$, together
with properties of the integer surgeries long exact sequence, that we
can take $N=2g-1$, where $g$ denotes the genus of the knot.
\end{remark}

In the same vein, for each $m\in\Z$,  we define a map
$$\Psi\colon \CFinf(Y_p,[m])\longrightarrow \CFKinf(Y,K,m)$$
by
$$\Psi[\x,i]=\sum_{\y\in\Ta\cap\Tc}
\sum_{\{\psi\in\pi_2(\x,\Theta,\y)\big|
\Mas(\psi)=0, n_{\BasePt}(\psi)-n_{\FiltPt}(\psi)=m\}}
\#\ModFlow(\psi)[\y,i-n_\BasePt(\psi),
i-n_\FiltPt(\psi)].$$
Note that the analogue of Equation~\eqref{eq:SpinCFormula} in 
this case guarantees
that the right-hand-side is contained in $\CFKinf(Y,K,m)\subset \CFKinf(Y,K)$.
We now have the following analogue of
Theorem~\ref{thm:LargeNegSurgeries}:

\begin{theorem}
\label{thm:LargePosSurgeries}
Let $(Y,K)$ be an oriented knot, and fix a compatible Seifert surface
$F$. For each integer $m$, we have an integer $N$ so that for all
integers $p\geq N(m)$, there is a Heegaard diagram for which the map
$\Psi$ above induces isomorphisms of chain complexes which fit into
the following commutative diagram: $$
\begin{CD}
0@>>>\CFm(Y_{p}(K),[m])@>>> \CFinf(Y_{p}(K),[m])
@>>>\CFp(Y_{p}(K),[m])@>>> 0, \\
&& @V{\Psim}VV @V{\Psi^\infty}VV @VV{\bPsip}V \\
0@>>>\CFKm(Y,K,m)@>>> \CFKinf(Y,K,m) @>>>\bCFKp(Y,K,m)@>>> 0.
\end{CD}
$$ 
Similarly, if we let $\bCFKsa(Y,K)\subset \CFKp(Y,K)$ be the subset
generated by $[\x,i,j]$ with $\max(i,j)=0$, then $\Psi$ also induces
isomorphisms in: $$
\begin{CD}
0  @>>>\CFa(Y_{p},[m])@>>>\CFp(Y_{p},[m]) @>{U}>> \CFp(Y_{p},[m]) @>>> 0 \\
&& @V{\bPsia}VV @V{\bPsip}VV @VV{\bPsip}V \\ 
0 @>>>\bCFKsa(Y,K,m)@>>> \bCFKp(Y,K,m) @>{U}>>\bCFKp(Y,K,m)@>>> 0.
\end{CD}
$$
\end{theorem}

\begin{proof}
This is a straightforward modification of the proof of
Theorem~\ref{thm:LargeNegSurgeries}.
\end{proof}

It is worth pointing out the following result:

\begin{cor}
\label{cor:HighestDegree}
Suppose that $K$ is a knot in $S^3$, and let 
$d$ be the largest integer for which $\HFKa(Y,K,d)\neq 0$, 
and suppose that $d>1$.
Then, 
$$\HFp(S^3_0(K),\spinc)\cong \HFKa(S^3,K,d)$$
as $\Zmod{2}$ graded $\Z$-modules,
where here $\spinc\in\SpinC(S^3_0(K))$ is the $\SpinC$ structure with
$\langle c_1(\spinc),[{\widehat F}]\rangle =2d-2$.
\end{cor}

\begin{proof}
On the one hand, we have the short exact sequence
$$
\begin{CD}
0 @>>> C\{i<0~{\text{and}}~j\geq d-1\} @>>> C\{i\geq 0~{\text{or}}~j\geq d-1\}
@>{\Psi}>> C\{i\geq 0\}@>>> 0.
\end{CD}$$
Now, it is easy to see by taking filtrations that
$$H_*(C\{i<0~{\text{and}}~j\geq d-1\})\cong \HFKa(S^3,K,d).$$
Of course, according to Theorem~\ref{thm:LargePosSurgeries}, for sufficiently
large integers $n$, 
$$H_*(C\{i\geq 0~{\text{or}}~j\geq d\})\cong \HFp(S^3_n(K),[d-1]),$$
and indeed the map $\Psi$ is modeled on the map induced by the cobordism
from $S^3_n(K)$ to $S^3$, endowed with the $\SpinC$ structure $\spincf$ 
with 
$$\langle c_1(\spincf),[{\widehat F}]\rangle = 2d-2+n.$$ 
Note that the map $\psi$ induced by $\Psi$ on homology is surjective, since it is a 
$\Z[U]$ module map and
$$H_*(C\{i\geq 0\})\cong \HFp(S^3),$$
and $\psi$ induces an isomorphism in all sufficiently large degrees.

We compare this with the integral surgeries long exact sequence,
according to which we have $$
\begin{CD}
... @>>>\HFp(S^3_0(K),\spinc) @>>> \HFp(S^3_p(K),[d-1]) @>{F}>> \HFp(S^3) @>>> ...
\end{CD}
$$
Now, if we endow $\HFp(S^3_p(K),[d-1])$ and $\HFp(S^3)$ with the filtrations
given by the absolute grading, then $F$ has the form $\Psi+L$, where $L$
is a sum of homogeneous maps which map to lower order than $\Psi$ (by at least $2d-2$,
c.f. Equation~\eqref{eq:DimensionFormula}). It then follows that the
induced map on homology $\psi+\ell$ is surjective, and also that
this kernel is identified with the kernel of $\psi$. 
\end{proof}

\section{Adjunction inequalities}
\label{sec:Adj}

In this section, we prove the following adjunction inequality for $\HFKa$.

\begin{theorem}
\label{thm:Adjunction}
Let $K\subset Y$ be an oriented knot, and 
suppose that 
$\HFKa(Y,K,\spincrel)\neq 0$. 
Then, for each Seifert surface $K$
for $F$ of genus $g>0$, we have that
$$\Big|\langle c_1(\spincrel),[{\widehat F}] \rangle \Big| \leq 2g(F).$$
\end{theorem}

\begin{proof}
Following Section~\ref{HolDiskTwo:sec:Adjunction} of~\cite{HolDiskTwo} (see especially
Lemma~\ref{HolDiskTwo:lemma:SpecialHeegaard}), we can construct a
special Heegaard diagram
$(\Sigma,\alphas,\betas_0\cup\mu,\BasePt,\FiltPt)$ for the knot $K$
with the following properties.  In the corresponding diagram for
$Y_0(K)$ (where the meridian $\mu$ is exchanged for the longitude
$\lambda$, but we have not yet wound the longitude along $\mu$, as
pictured in Figure~\ref{fig:WindConvention}), there is a periodic
domain $\PerDom$ which bounds $\lambda\cup\alpha_1$, all its
multiplicities are one and zero and its Euler characteristic is
$-2g$. Moreover, the meridian $\mu$ meets the support of $\PerDom$ in
an arc which meets none of the other attaching circles. 

Let $\x$, now, be any intersection point generating $\CFKa(Y,K)$,
and let $\x'$ be the nearby intersection for the $Y_0$ Heegaard diagram,
obtained now after winding. Observe that for this Heegaard diagram, there is 
a variant of $\PerDom$, denoted $\PerDom'$, where there are multiplicities $-1$,
as well, now along some neighborhood of the meridian $\mu$. Indeed, if
$x_1\in \x$ denotes the intersection point on $\mu$, then $x_1'$ is supported
in this region. However, none of the other $x_i\in\x'$ are supported in the
winding region. Note that there must be another $x_2\in\x'$ which is supported
on $\alpha_1$ (i.e. it lies on the boundary of $\PerDom'$, where the multiplicity
on one side is $+1$). See Figure~\ref{fig:Adjunction} for an illustration.

\begin{figure}
\mbox{\vbox{\epsfbox{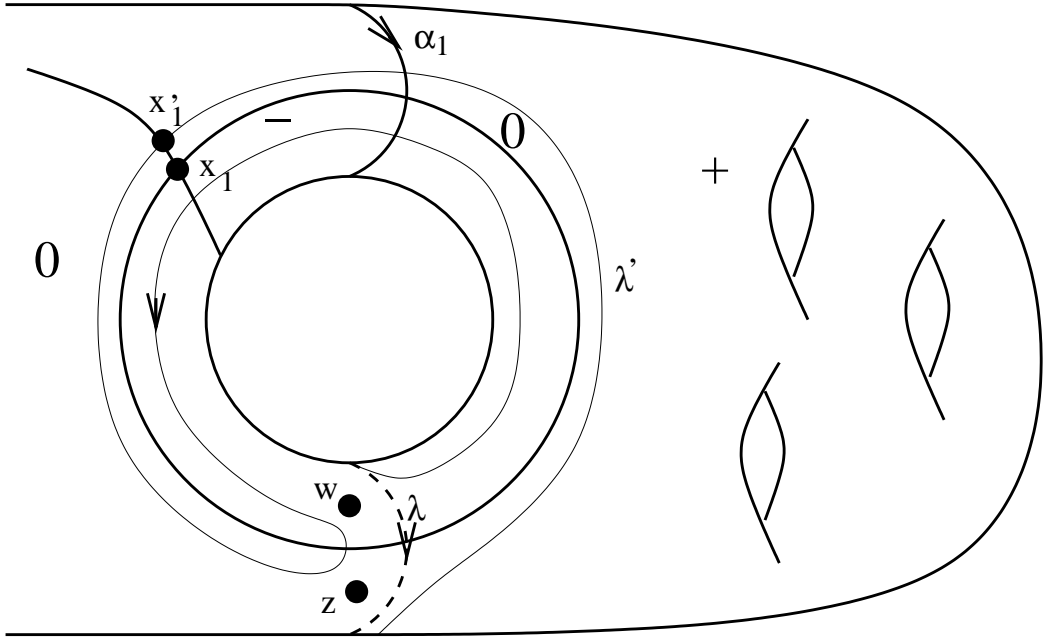}}}
\caption{\label{fig:Adjunction}
{\bf{Adjunction inequality.}}
An illustration for the proof of Theorem~\ref{thm:Adjunction}. Here,
the dashed circle is the original longitude, while $\lambda'$ is the one
obtained after winding. We have indicated some of the local multiplicities
of the periodic domain $P'$ obtained after winding, by writing $+$ for the
regions where its local multiplicity is $+1$, $-$ where it is $-1$, and $0$
where it vanishes. (These can be obtained by realizing that the boundary of
$P$ is $\alpha_1\cup \lambda'$, with the orientations indicated by the arrows.)}
\end{figure}

Now, according to Equation~\eqref{eq:MeasureCalc}, we see that
$$\langle c_1(\spinc'(\x')),[{\widehat F}]\rangle
= -2g + \#(x_i~\text{in the interior of $\PerDom$}) \geq -2g.$$
The result now follows by conjugation 
invariance (Proposition~\ref{prop:JInvariance}).
\end{proof}

\section{Examples}
\label{sec:Examples}

Before turning to some general calculational devices for $\CFKinf$, we
give here a few Heegaard genus two examples where $\CFKinf(Y,K)$ can be
explicitly determined through more direct means.

The key trick which facilitates these calculations is the following:

\begin{prop}
\label{prop:Destabilization}
Let $(\Sigma,\alphas,\betas,\StartPt,\EndPt)$ be a doubly-pointed
Heegaard diagram of genus $g$. Suppose that the curves $\alpha_1$ and
$\beta_1$ intersect in a single point, and suppose that $\beta_1$
meets none of the other $\alpha_i$ for $i>1$.  Consider
the doubly-pointed Heegaard diagram
$(\Sigma',\alphas',\betas',\StartPt',\EndPt')$
of genus $g-1$, where $\Sigma'$ is obtained from $\Sigma$ by surgering
out $\beta_1$, and $\alphas'$, $\betas'$, $\StartPt'$, and $\EndPt'$
are obtained by viewing $\{\alpha_2,...,\alpha_g\}$,
$\{\beta_2,...,\beta_g\}$, $\StartPt$, and $\EndPt$ as 
supported in
$\Sigma'$. Then,
$$\CFinf(\Sigma,\alphas,\betas,\StartPt,\EndPt)\cong 
\CFinf(\Sigma',\alphas',\betas',\StartPt',\EndPt').$$
\end{prop}

\begin{proof}
In the case where $\alpha_1$ 
meets none of the 
curves in $\{\beta_2,...,\beta_g\}$,
this is a
statement of the stabilization invariance of the doubly-pointed complex.
It is easy to reduce to this case by performing a series of handleslides
of the $\beta_i$ (for $i>1$) which meet $\alpha_1$, over $\beta_1$
(along a subarc in $\alpha_1$ which connects $\beta_i$ with $\beta_1$).
Note that these handleslides do not cross either basepoint.
\end{proof}

\subsection{The trefoil}
\label{subsec:Trefoil}

\begin{figure}
\mbox{\vbox{\epsfbox{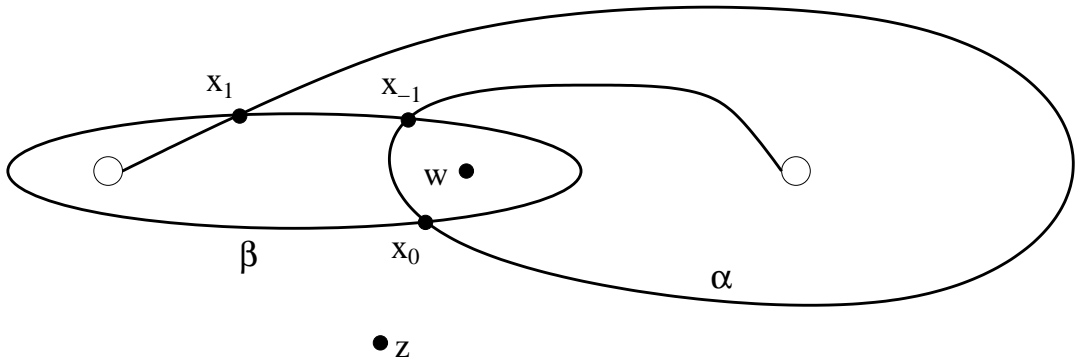}}}
\caption{\label{fig:Trefoil}
{\bf{Destabilized doubly-pointed Heegaard diagram for the trefoil.}}
This picture is meant to take place in a torus, where the two hollow
circles are identified.}
\end{figure}

Start with a genus two Heegaard diagram for the left-handed
trefoil. This can be
destabilized once (Proposition~\ref{prop:Destabilization}) to obtain
the following picture shown in Figure~\ref{fig:Trefoil}, which takes
place inside a torus. Note that we have destabilized out the meridian
of the knot. 


\begin{prop}
\label{prop:Trefoil}
The filtered chain complex $\CFKinf(Y,K)$ is freely generated as a
$\Z$-module by generators $[x_1,i,i+1]$, $[x_{-1},i+1,i]$, and
$[x_0,i,i]$, where $i$ is an arbitrary integer. The differential given
by:
\begin{eqnarray*}
\partial [x_0,i,i]&=&0 \\
\partial [x_1,i,i+1]&=& [x_0,i,i] \\
\partial [x_{-1},i+1,i]&=& [x_0,i,i]. 
\end{eqnarray*}
Moreover, the absolute grading of the generator $[\x_0,0,0]$ is $1$.
\end{prop}

\begin{proof}
That the chain complex has the stated form follows readily from the
fact there are two orientation-preserving embedded Whitney disks in
the torus: $\phi_{-1}$, which connects $x_{-1}$ to $x_0$, and
$\phi_1$, which connects $x_1$ to $x_0$. So, by the Riemann mapping
theorem, both admit unique holomorphic representatives (up to
translation). Moreover, $\phi_{-1}$ contains the basepoint $w$ (and
not $z$), while $\phi_1$ contains $z$ (and not $w$).

Strictly speaking, this information suffices to identify the chain complex $\CFinf$ for
the doubly-pointed Heegaard diagram (generated by $[x,i,j]$ for
intersection points $x\in\{x_{-1},x_0,x_1\}$, and $i,j\in\Z$), which
splits as a direct sum of chain complexes isomorphic to $\CFKinf(S^3,K)$.
It is immediate from the symmetry property of $\HFKa$
(c.f. Proposition~\ref{prop:IndepOrient}) that the summand
$\CFKinf(S^3,K)$ is generated by $[x_0,i,i]$, $[x_1,i,i+1]$ and
$[x_{-1},i+1,i]$.

To calculate the absolute grading, we use Lemma~\ref{lemma:AbsGrade},
and the observation that the complex generated by $[x_1,0,1]$,
$[x_0,0,0]$, and $[x_{-1},0,-1]$ (thought of as a quotient of a
subcomplex of $\CFKinf(Y,K)$) calculates $\HFa(S^3)$. In this complex,
the generator $[x_{-1},0,-1]$ is the only one which persists in homology,
so it must have degree zero.
\end{proof}

As an easy application of the results from
Section~\ref{sec:Relationship}, one could give yet another calculation
of $\HFp$ for three-manifolds which are surgeries along the trefoil.
Moreover, the technique employed above can be readily generalized to
give a combinatorial description of $\CFK$ for an arbitrary two-bridge
knot. However, since $\HFp$ of three-manifolds obtained as integral
surgeries along such knots has already been calculated by Rasmussen
in~\cite{Rasmussen}, we do not pursue this any further here,
turning instead to a related calculation.

\subsection{Examples from two-bridge links}

There is a class of knots admitting genus two Heegaard diagrams, but
which are not two-bridge knots in the usual sense. Our aim here is to
show how to reduce the calculation of the complex $\CFK$ for this
class of knots to a purely combinatorial problem
(c.f. Propositions~\ref{prop:HeegaardDiagram} and
\ref{prop:Differential} below).

Recall that a two-bridge link in $S^3$ is a link which is contained in
a Euclidean neighborhood $B$ on which there is function $f\colon
B\longrightarrow \R$ whose restriction to the link has exactly two
local maxima. We consider links of this type which have two
components. Such links can be described by certain (four-stranded)
braids $\sigma$. Specifically, think of a braid as a mapping class
$\sigma$ of the plane with four marked points $\{a_1,b_1, a_2,b_2\}$
which permutes the points. Let $I_1$ and $I_2$ be a pair of disjoint arcs
which connect $a_1$ to $b_1$ and $a_2$ to $b_2$ respectively.  We
construct (the projection of) a link by joining $I_1$ and $I_2$ to $\sigma(I_1)$
and $\sigma(I_2)$ respectively along their boundaries, so that $I_1$ and $I_2$ pass
``over'' their images under $\sigma$. When the braid $\sigma$ fixes
the subsets $\{a_1,b_1\}$ and $\{a_2,b_2\}$, the link has two
components. In fact, by introducing twists which leave the isotopy
class of the link unchanged, we can assume the braid fixes the four
points (and indeed that it fixes a small tubular neighborhood of these
points).

\begin{figure}
\mbox{\vbox{\epsfbox{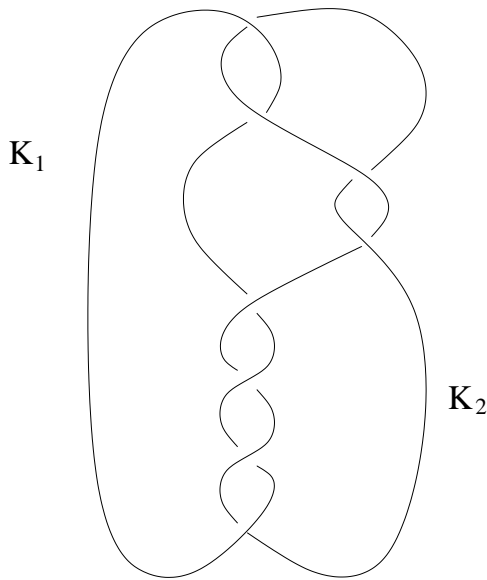}}}
\caption{\label{fig:TwoBridgeLink}
{\bf{A two-bridge link.}} }
\end{figure}

\begin{figure}
\mbox{\vbox{\epsfbox{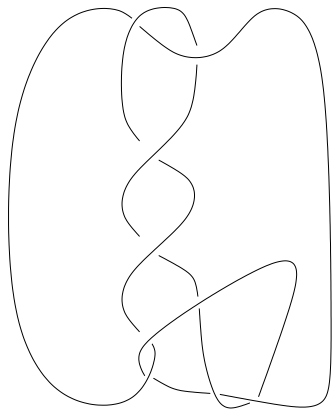}}}
\caption{\label{fig:BigKnot}
{\bf{A knot with Heegaard genus two.}} This knot is obtained by
``blowing down'' (with sign $-1$) the component $K_1$ of the
link from Figure~\ref{fig:TwoBridgeLink}.}
\end{figure}

Let $K_1\cup K_2$ be a two-component two-bridge link.  It is easy to
see that the two components of the link are each individually
unknotted. Thus, if we perform $\pm 1$ surgery on $K_1$ (or, more
generally, $1/n$ where $n$ is an arbitrary integer), then the ambient
three-manifold is still diffeomorphic to $S^3$, and we let $K$ denote
the image of $K_2$ under this diffeomorphism. Thus, $K$ is obtained
from $K_2$ by introducing a ``Rolfsen twist'' along the linking circle
$K_1$ (see for example~\cite{GompfStipsicz} for a detailed description
of this operation; see Figure~\ref{fig:BigKnot} for a picture of the
result of performing the operation on the link shown in
Figure~\ref{fig:TwoBridgeLink}).

Let $\sigma$ denote the braid associated to a two-component link
$K_1\cup K_2$. We construct a
doubly-pointed genus one Heegaard diagram for $K$ from $\sigma$ as
follows.  Delete a pair of disks $A_1$ and $B_1$ centered at $a_1$ and
$b_1$. The Heegaard surface is formed by attaching a cylinder $C_1$ to
$S^2-A_1-B_1$ along its two boundary circles.  We let $a_2$ and $b_2$
serve as our two basepoints $\StartPt$ and $\EndPt$ respectively. We
consider the curve $\beta$ obtained from closing up the arc $I_1\cap
(S^2-A_1-B_1)$ by an arc inside the attached cylinder $C_1$. Similarly, we
consider the curve $\alpha$ obtained by closing up $\sigma(I_1)$ by an
arc inside $C_1$. Note that there is some freedom in
how we close off $\alpha$ and $\beta$ (inside $C_1$). In
particular, by introducing Dehn twists supported in the cylinder 
on the $\alpha$, 
we can arrange for the algebraic intersection number
$\#\alpha\cap\beta$ to be an arbitrary integer.

\begin{prop}
\label{prop:HeegaardDiagram}
Consider the doubly-pointed, genus one Heegaard diagram obtained as
above, where we have arranged for the curves $\alpha$ and $\beta$ to
have $\#\alpha\cap\beta= \pm 1$. This Heegaard diagram is a
destabilization of the Heegaard diagram for the knot $K$ obtained
from the two-bridge link $K_1\cup K_2$ by ``blowing down'' $K_1$
(given framing $\pm 1$).
\end{prop}

\begin{proof}
In the standard way (see, for example,~\cite{GompfStipsicz}), a two-bridge
link gives a genus two Heegaard diagram of $S^3$. This is related to
the automorphism $\sigma$ as follows. We consider the sphere with four
distinguished disks $A_1$, $B_1$, $A_2$, and $B_2$ (neighborhoods of
the $a_1$, $b_1$, $a_2$, and $b_2$ respectively) each fixed by
$\sigma$. The Heegaard surface $\Sigma_2$ is obtained by attaching
cylinders $C_1$ and $C_2$ to $S^2-A_1-B_1-A_2-B_2$ which connect
$\partial A_1$ to $\partial B_1$ and $\partial A_2$ to $\partial B_2$
respectively.  Meridians $\mu_1$ and $\mu_2$ for the two components
$K_1$ and $K_2$ are appropriately oriented, embedded, homologically
non-trivial circles supported in $C_1$ and $C_2$ respectively. The
remaining circles, $\gamma_1$ and $\gamma_2$, are obtained by closing
off the arcs $\sigma(I_1)$ and $\sigma(I_2)$ inside $C_1$ and $C_2$
respectively. Longitudes $\lambda_1$ and $\lambda_2$ for $K_1$ and
$K_2$ are obtained by closing off $I_1$ and $I_2$ in $C_1$ and $C_2$
respectively. Thus, $(\Sigma_2,\{\gamma_1,\gamma_2\},\{\mu_1,\mu_2\})$
describes $S^3$.  The Heegaard diagram for a three-manifold obtained
by integral surgery along $K_1$ is obtained from the above by replacing $\mu_1$
by the curve $\delta$ obtaining by juxtaposing $\lambda_1$ and some
number of copies of $\mu_1$. Thus,  when $\delta$ is chosen so that
$\#(\gamma_1\cap \delta)=\pm 1$, then $(\Sigma_2,\{\gamma_1,\gamma_2\},
\{\delta\},\mu_2)$ is a Heegaard diagram for the knot induced from $K_2$
inside $S^3_{\pm 1}(K_1)$. The corresponding doubly-pointed Heegaard
diagram
$(\Sigma_2,\{\gamma_1,\gamma_2\},\{\delta,\mu_2\},\StartPt,\EndPt)$
clearly has the property that $\mu_2$ meets $\gamma_2$ in a single
point (inside $C_2$). Thus, we can destablizing as in
Proposition~\ref{prop:Destabilization}, to obtain the doubly-pointed
Heegaard diagram described in the statement of the present proposition.
\end{proof}

Thus, $\CFK$ for knots obtained from two-bridge links is reduced to a
calculation of the doubly-pointed chain complex for the torus. Such a
calculation is purely combinatorial, i.e. independent of the choice of
complex structure on the torus, as we shall see in the following results
(see especially Proposition~\ref{prop:Differential} below).

When analyzing Floer homology taking place in the torus $T$, 
we find it convenient to pass to the  universal covering
space $\pi\colon \C\longrightarrow T$. To this end, we will assume
that ${\alpha}$ and ${\beta}$ are a pair of homotopically non-trivial,
embedded curves in the torus $T$.  In this case, 
$\pi^{-1}(\alpha)$
(and also $\pi^{-1}(\beta)$) is a properly
embedded submanifold of $\C$ homeomorphic to $\R\times \Z$.
Let ${\widetilde \alpha}$ be one connected
component of $\pi^{-1}(\alpha)$ and ${\widetilde\beta}$ be a connected
component of $\pi^{-1}(\beta)$. Given reference points $w$ and
$z$ in the torus, $\pi^{-1}(w)$ and $\pi^{-1}(z)$
lift to lattices $W$ and $Z$ respectively. 
Moreover, $\CFKinf(S^3,K)$ is
generated by triples $[{\widetilde x},i,j]$ subject to the relation
in Equation~\eqref{eq:DefOfHFK} 
where here the ${\widetilde x}$ are
intersection points of ${\widetilde
\alpha}$ and ${\widetilde
\beta}$. Indeed, by path lifting,
given $x,y\in \alpha\cap \beta$, the space of homotopy classes of
disks $\pi_2(x,y)$ is identified with the space of homotopy classes of
Whitney disks in $\C$ interpolating connecting the corresponding
lifts ${\widetilde x}, {\widetilde y}\in
{\widetilde\alpha}\cap{\widetilde \beta}$ (a space denoted
${\widetilde \phi}\in\pi_2({\widetilde x},{\widetilde
y})$). Clearly, if $\phi\in\pi_2(x,y)$ and $\phiCov\in\pi_2(\xCov,\yCov)$
is its lift, then 
$n_w(\phi)$ is obtained as a sum over all $\wCov\in W$
of $n_{\wCov}(\phiCov)$.
It is also not difficult to see an identification of moduli
spaces $\ModFlow(\phi)\cong
\ModFlow({\widetilde \phi})$.

\begin{prop}
\label{prop:Differential}
Suppose that $\alpha$ and $\beta$ are two homotopically non-trivial,
embedded curves in the torus $T$, and let $\alphaCov$ and $\betaCov$
be corresponding curves in the universal cover $\C$.  Let
${\widetilde x}, {\widetilde y}\in {\widetilde \alpha}\cap{\widetilde
\beta}$ be a pair of intersection points, and let ${\widetilde
\phi}\in\pi_2({\widetilde x},{\widetilde y})$ be the homotopy class of
Whitney disk connecting ${\widetilde x}$ and ${\widetilde y}$. Then if
$\Mas(\phiCov)=1$, we have that $\#\UnparModFlow(\phiCov)=1$ if
and only if $\cald(\phiCov)\geq 0$.
\end{prop}

Before proving this proposition, it is useful to have the following
formula for the Maslov index.  To state it, we define the local
multiplicity of $\phiCov$ 
at $\xCov$ as follows. Given
$\phiCov\in\pi_2(\xCov,\yCov)$, 
there are four regions in $\C-\alphaCov-\betaCov$ which contain
$\xCov$ as a corner point. Choosing interior points $z_1,...,z_4$ in
these four regions, define
\begin{equation}
\label{eq:DefMult}
{\overline n}_{\xCov}(\phiCov)=
\frac{n_{z_1}(\phiCov)+n_{z_2}(\phiCov)+n_{z_3}(\phi)+n_{z_4}(\phiCov)}{4}.
\end{equation}
We define ${\overline n}_{\yCov}(\phiCov)$ similarly, only now choosing
the four points in the region containing $\yCov$ as a boundary.

\begin{lemma}
\label{lemma:MasFormula}
Suppose that $\alpha$ and $\beta$ are two homotopically non-trivial,
embedded curves in the torus $T$, and let $\alphaCov$ and $\betaCov$
be corresponding curves in the universal cover $\C$.  
Given $\xCov, \yCov\in\piCov(\alphaCov,\betaCov)$,
$$\Mas(\phiCov)=2({\overline n}_{\xCov}(\phiCov)+{\overline n}_{\yCov}(\phiCov)).$$
\end{lemma}

\begin{proof}
As we go from $\xCov$ to $\yCov$ in $\alphaCov$, we encounter finitely many
intersection points $\alphaCov\cap\betaCov$.  Let
$\{\xCov_1,...\xCov_n\}$ denote this sequence of points, arranged in
the order they are encountered, so that $\xCov_1=\xCov$ and
$\xCov_n=\yCov$.  There is a canonical Whitney disks ${\widetilde
\phi}_i\in{\widetilde \pi}_2({\widetilde x}_i,{\widetilde x}_{i+1})$: 
the domain is bounded by the Jordan curve obtained by juxtaposing the
two subarcs of ${\widetilde \alpha}$ and ${\widetilde \beta}$ which
connect ${\widetilde x}_i$ and ${\widetilde x}_{i+1}$ (oriented so
that we go from $\xCov_i$ to $\xCov_{i+1}$ in the
$\alphaCov$-arc). Thus, the support of $\cald(\phi_i)$ is an embedded
disk whose local multiplicity is either $+1$ or $-1$.  Clearly, the,
homotopy class $\phiCov$ can be decomposed as juxtapositions of these
canonical Whitney disks, i.e.  $\phiCov= \phiCov_1*...*\phiCov_{n-1}$.

We prove the lemma by induction on the length
$n$. Clearly, when $n=1$, both quantities
are zero.

For the inductive step decompose $\phiCov=\phiCov_1*\phiCov'$, where
$\phiCov'\in \piCov(\xCov_{2},\xCov_n)$ has length $n-1$.  Let
$\epsilon$ denote the local multiplicity inside $\cald(\phi_1)$.  A
case-by-case analysis shows that
\begin{equation}
\label{eq:DesiredEq}
2(n_{\xCov_1}(\phiCov)+n_{\xCov_2}(\phiCov))=
2(n_{\xCov_{2}}(\phiCov')+n_{\xCov_n}(\phiCov'))+\epsilon;
\end{equation}
while it is easy to see that $\Mas(\phiCov_1)=\epsilon$. Thus, the
result follows from Equation~\eqref{eq:DesiredEq},
using the inductive hypothesis, and the additivity
of the Maslov index under juxtaposition.

The cases required for establishing Equation~\eqref{eq:DesiredEq}
are divided according to the order in which $\xCov_1$, $\xCov_2$,
and $\xCov_n$ are encountered as we traverse $\betaCov$. Observe that for
this ordering, $\xCov_n<\xCov_1$.

For example, suppose that $\xCov_n< \xCov_{1}< \xCov_2$.  In this
case, the local multiplicities of the four regions containing
$\xCov_1$ have the form $m+\epsilon$, occurring three times, and $m$,
occurring once, for some integer $m$. Now,
traversing the arc in $\alphaCov$ from $\xCov_1$ to $\xCov_{2}$, we
see that the four neighboring multiplicities are: $m$ (with multiplicity two)
and $m+\epsilon$ (also with multiplicity two). It follows readily that
\begin{eqnarray*}
n_{\xCov_{2}}(\phiCov')=n_{\xCov_1}(\phiCov)-\OneHalf\epsilon, &&
n_{\xCov_n}(\phiCov')=n_{\xCov_n}(\phiCov),
\end{eqnarray*}
verifying Equation~\eqref{eq:DesiredEq}.
The case where $\xCov_2 < \xCov_n < \xCov_{1}$ works the same way.

Finally, when $\xCov_n< \xCov_2 <\xCov_{1}$, we have instead that
\begin{eqnarray*}
n_{\xCov_{2}}(\phiCov')=n_{\xCov_1}(\phiCov), &&
n_{\xCov_n}(\phiCov')=n_{\xCov_2}(\phiCov)-\OneHalf\epsilon,
\end{eqnarray*}
also giving Equation~\eqref{eq:DesiredEq}.
\end{proof}

We now turn to a proof of Proposition~\ref{prop:Differential}.

\vskip.2cm 
\noindent{\bf{Proof of  Proposition~\ref{prop:Differential}.}}
Clearly, if $\phiCov$ has a holomorphic representative, then
$\cald(\phiCov)\geq 0$. 

In the more interesting direction, suppose that
$\phiCov\in\pi_2(\xCov,\yCov)$ is a homotopy class with
$\Mas(\phiCov)=1$ and $\cald(\phiCov)\geq 0$. Clearly, $\phiCov$ is
determined by the restriction of $\phiCov$ to its boundary.  Moreover,
its boundary curve
$C$ is formed from an arc $a$ in $\alphaCov$ from $\xCov$ to $\yCov$ 
and an arc 
$b$ in $\betaCov$ from $\yCov$ to $\xCov$. Of course, $C$ is
null-homotopic. Our goal is to construct a covering space of the plane
with two disks removed, in which the lift $C'$ of $C$ is a
null-homotopic, embedded circle.

\begin{figure}
\mbox{\vbox{\epsfbox{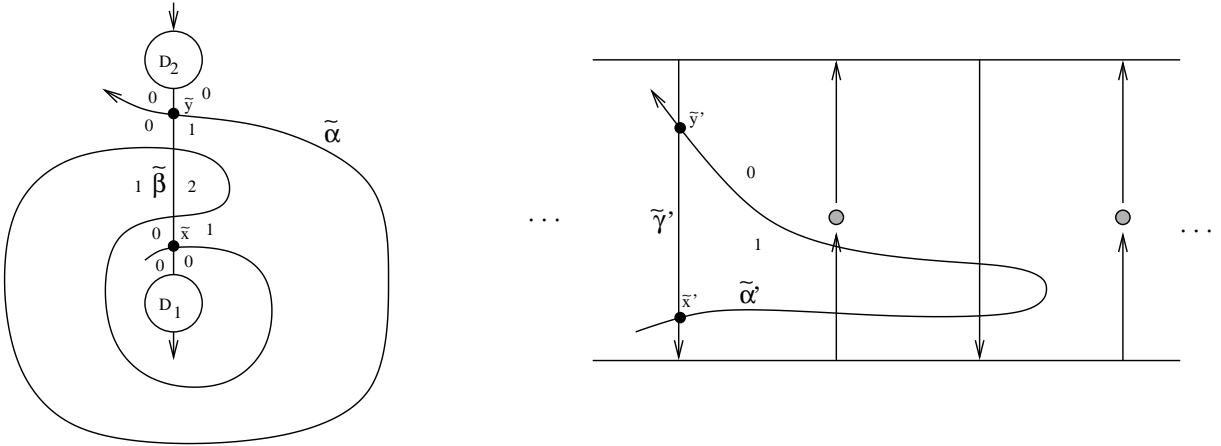}}}
\caption{\label{fig:PosDomain}
{\bf{Lifting a domain.}} On the left-hand-side, we have illustrated a 
domain connecting $\xCov$ and $\yCov$, corresponding to a map $\phiCov$
with
$\Mas(\phiCov)=1$. The ambient space is $\C$, with the disks
$D_1$ and $D_2$ deleted. 
On the right, we have illustrated an induced domain in
$X$. Here, the picture is an infinite strip with an infinite sequence
of puncture points, indicated here by the grey circles.
The top boundary
corresponding to the universal cover of $\partial D_2$, 
and the bottom corresponding to the unversal cover of $\partial D_1$.}
\end{figure}

To this end, Lemma~\ref{lemma:MasFormula}, together with the
hypothesis that $\cald(\phiCov)\geq 0$, shows that three out of the
four local multiplicities around each of $\xCov$ and $\yCov$ vanish,
and the remaining ones are $+1$. We choose a pair of disks, $D_1$ and
$D_2$, centered at points on $\betaCov$ with $D_1$ just after $\xCov$,
and the other just before $\yCov$ in $\betaCov$ (with respect to the
orientation induced on $\betaCov$ from $\partial\cald(\phiCov)$),
choosing the disks small enough to be disjoint from the support of
$\cald(\phiCov)$, and so that they are disjoint from
$\alphaCov$. Clearly, $(D_1\cup D_2)\cap \betaCov$ separates
$\betaCov$ into three components, one of which is compact (containing
both $\xCov$ and $\yCov$) and two of which are non-compact. Let
$\gamma$ denote the compact subarc.  Intersection number with $\gamma$
induces a homomorphism of $\pi_1(\C-D_1-D_2)$ to $\Z$, and hence a
covering space $X$ of $\C-D_1-D_2$. This covering space is easily seen
to be an infinite strip with a discrete, countable set of points
$\Lambda$ removed (corresponding to the point at infinity in
$\C-D_1-D_2$). In $X$, the pre-image of $\gamma$ is an infinite
collection of separating arcs (connecting the two boundaries), of
which we choose one, denoted $\gamma'$. Similarly, the pre-image of
$\alphaCov$ is an infinite collection of arcs (which begin and end in
the puncture points), of which we choose one, $\alphaCov'$. Let
$\xCov'$ and $\yCov'$ be the lifts in $\alphaCov'\cap \gamma'$ of
$\xCov$ and $\yCov$ respectively.  (c.f. Figure~\ref{fig:PosDomain}
for an illustration).

\begin{figure}
\mbox{\vbox{\epsfbox{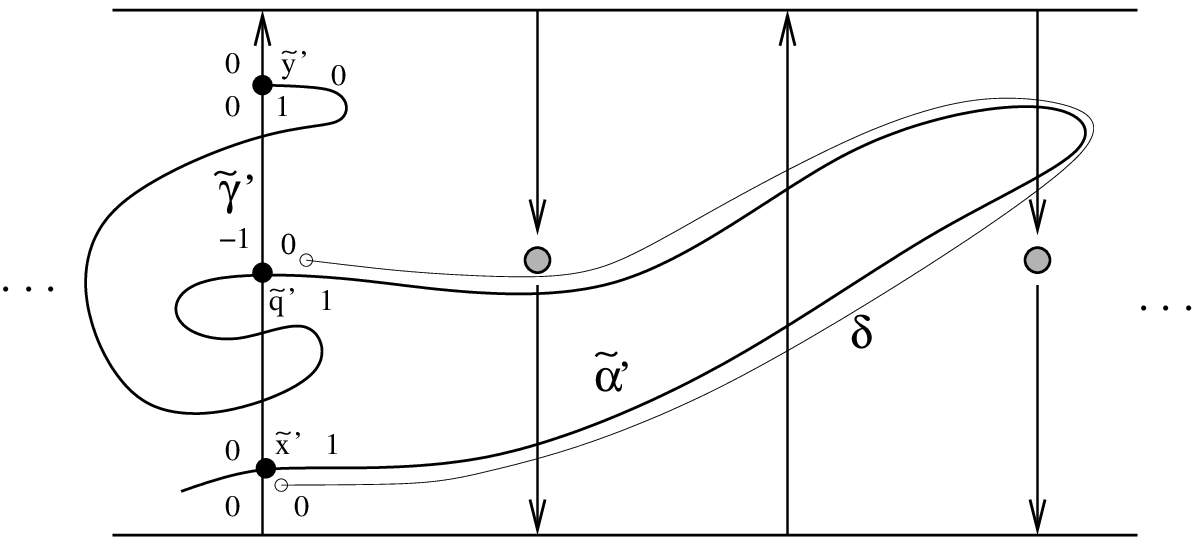}}}
\caption{\label{fig:NegDomain}
{\bf{Proof of Proposition~\ref{prop:Differential}.}}  In the case
where the arc in $a'$ meets the arc $b'$ in than two points, as we
have illustrated, $\cald(\phiCov)$ is forced to be non-positive,
provided $\Mas(\phiCov)=1$. The light curve $\delta$ illustrated here
has projection to $\R^2-D_1-D_2$ which is disjoint from $\alphaCov$,
and intersection number zero with $\betaCov$. Thus, the multiplicities
of $\cald(\phiCov)$ at the two endpoints must coincide (indeed, they must
vanish). It follows that just on the other side of $\gammaCov'$ from
the endpoint near $\qCov'$, the local multiplicity of the domain
$\cald(\phiCov)$ is $-1$.}
\end{figure}

The lift $C'$ of $C$ to $X$, of course, is obtained by adjoining the
subarc $a'$ in $\alphaCov'$ from $\xCov'$ to $\yCov'$ to the subarc
$b'$ in $\betaCov'$ from $\yCov'$ to $\xCov'$. We claim that $a'$ and
$b'$ meet only along their boundary, and hence that $C'$ is an
embedded curve.  To this end, observe that $\gamma'$ separates $X$
into two halves. Call the ``right'' half the one which the curve $a'$
points towards, at $\xCov'$. Let $\qCov'$ be the first crossing (with
respect to the $a'$ orientation) after $\xCov$ encountered as we
traverse $a'$.  With respect to the induced orientation on $\gamma$,
$\xCov'$ and $\yCov'$ are extremal. It is an easy consequence now that
one of the four local multipilcities near $\qCov'$ is $-1$.  Indeed,
by following parallel to $a'$, we can find a curve $\delta$ supported
on the ``right'', connecting one of the regions near $\xCov'$ with
zero multiplicity to one of the four neighboring regions near
$\qCov'$, so that $\delta$ is supported near $\alphaCov'$, its
projection to $\C$ is disjoint from $\alphaCov$, and, since its
endpoints are both near $\gamma'$ but to the right, the intersection
number of the projection of $\delta$ with $\betaCov$ is zero. This
shows that in $\R^2-D_1-D_2$ there is a region containing $\qCov$ (the
projection to $\R^2-D_1-D_2$ of $\qCov'$) to the right of $\gamma$,
which has multiplicity zero. Since $\qCov'$ falls between $\xCov'$ and
$\yCov'$, it follows that the facing region on the other side of
$\gamma$ has multiplicity $-1$.  (See the illustration in
Figure~\ref{fig:NegDomain}.)  This contradicts the hypothesis that
$\cald(\phiCov)\geq 0$.

We have established that $C'$ is an embedded curve in $X$. Thus, by
the Jordan curve theorem, $C'$ separates the strip into two regions, a
disk $D$, and another region which contains the boundary of the strip.
Since $\cald(\phiCov)$ is compact it follows readily that none of the
points in $\Lambda$ is contained in the disk; i.e. $C'$ bounds a disk
in $X$. (More precisely, the disk has two corner points
at the two points in $C'$ where $a'$ and $b'$ meet.) 
It follows now by
the Riemann mapping theorem that there is a unique holomorphic Whitney
disk representing this disk; and hence by path-lifting,
it follows that
$\#\ModFlow(\phiCov)=1$, as well.
\qed
\vskip.2cm

\begin{figure}
\mbox{\vbox{\epsfbox{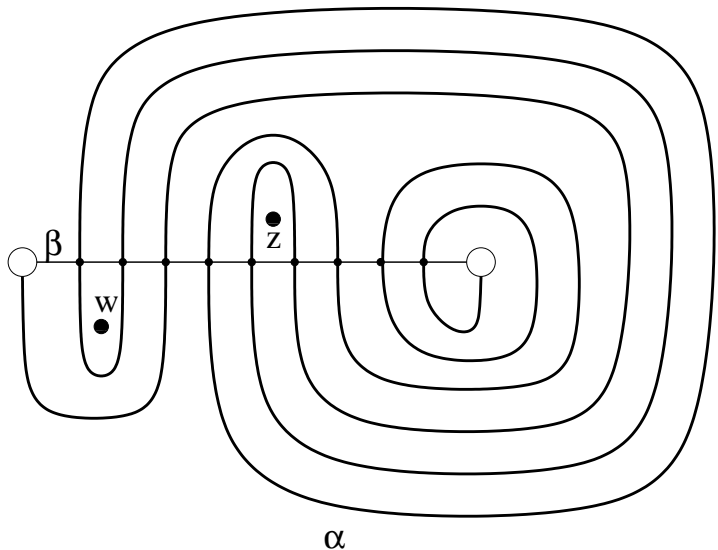}}}
\caption{\label{fig:HeegaardDiagram}
{\bf{A doubly-pointed Heegaard diagram.}} This is the Heegaard diagram
for the knot from Figure~\ref{fig:BigKnot}, as constructed in
Proposition~\ref{prop:HeegaardDiagram}. The two empty circles are glued
along a cylinder, so that no new intersection points are introduced
between the curve $\alpha$ (the darker curve) and $\beta$ (the
lighter, horizontal curve).  Note that there are nine intersection
points between $\alpha$ and $\beta$, which, when counted with sign,
give an intersection number of $+1$.}
\end{figure}

\begin{figure}
\mbox{\vbox{\epsfbox{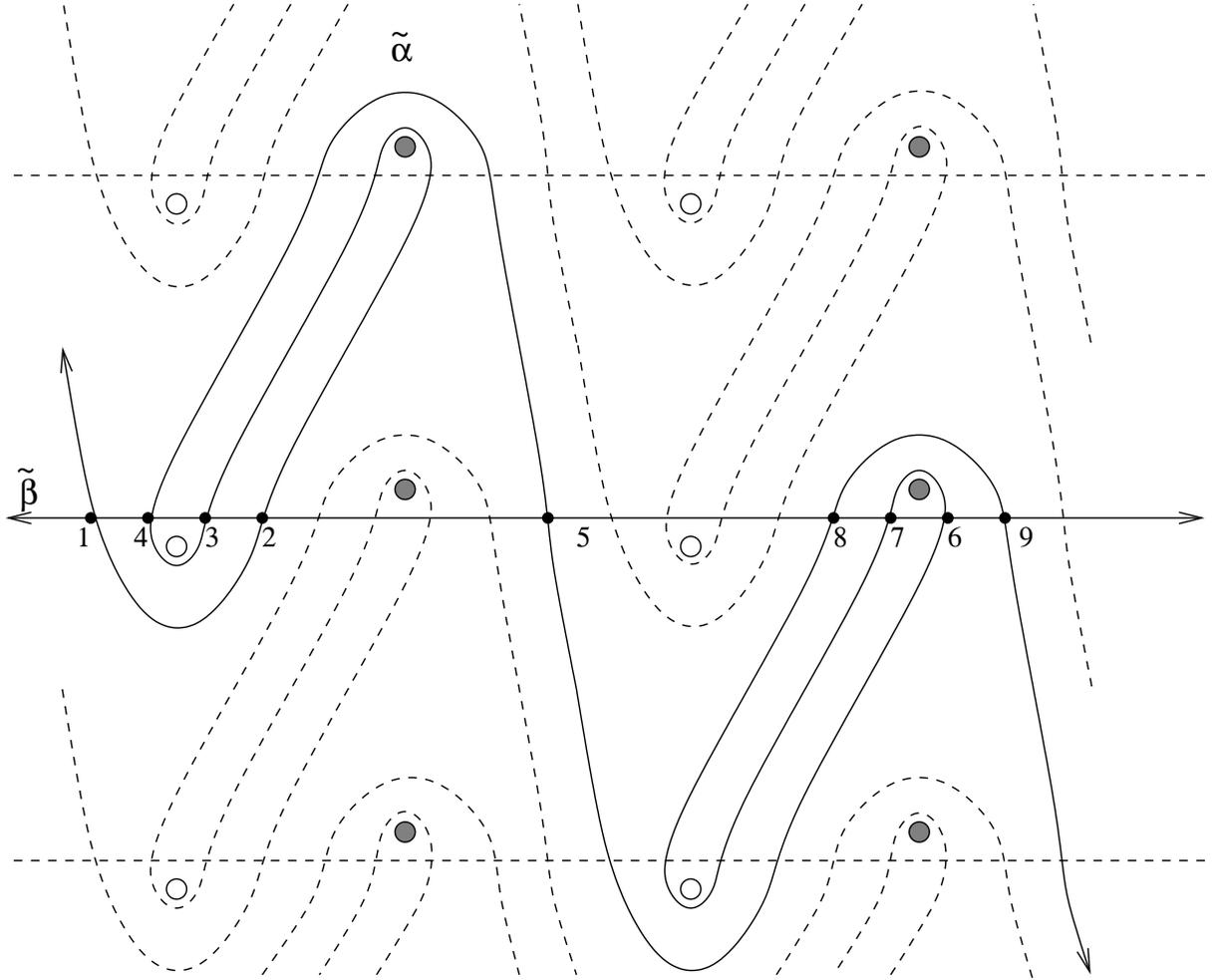}}}
\caption{\label{fig:CoverHeeg}
{\bf{Lift of Figure~\ref{fig:HeegaardDiagram} to the universal
cover.}}  Two of the lifts of $\alpha$ and $\beta$, ${\widetilde
\alpha}$ and ${\widetilde \beta}$, are pictured here as dark curves,
while other lifts are shown with dashed lines.  The lattice
$\pi^{-1}(\{w\})$ is represented by the lightly-filled circles, while the
lattice of $\pi^{-1}(\{z\})$ is represented by the darkly-filled
circles, while.  The nine intersection points of ${\widetilde
\alpha}$ and ${\widetilde \beta}$ are labeled in the order in which
they are encountered along ${\widetilde\alpha}$.}
\end{figure}

As an illustration of these techniques, consider the link pictured in
Figure~\ref{fig:TwoBridgeLink}. Label the generators of the braid
group on four elements by $\tau_1$, $\tau_2$, $\tau_3$ (so that
$\tau_i$ transposes the $i^{th}$ and $(i+1)^{st}$ strands).  The braid
element $\sigma$ corresponding to the link, then can be written as a
product $\tau_2^2 \cm \tau_3^2 \cm
\tau_2^{-4}$. The Heegaard diagram, then, is obtained by corresponding
Dehn twists in the sphere with four marked points. Blowing down the other
component  (with sign $-1$), we obtain the  diagram 
illustrated in Figure~\ref{fig:HeegaardDiagram}. Note that this
gives us the knot pictured in Figure~\ref{fig:BigKnot} (which
appears in knot tables as the knot $9_{42}$). Moreover, lifting this
doubly-pointed Heegaard diagram to the universal covering space $\C$
for the torus, we obtain the picture shown in
Figure~\ref{fig:CoverHeeg}.

In this case, ${\widetilde \alpha}\cap {\widetilde \beta}$ consists of
$n=9$ points.  Thus, $\C-{\widetilde \alpha}-{\widetilde\beta}$ has
eight compact, connected components (and, more generally, $n-1$
compact connected components), which we label $D_1,...,D_8$, under the
numbering convention as follows. Let $C_i$ denote the Jordan curve
obtained by juxtaposing the segment of ${\widetilde \alpha}$
connecting ${\widetilde x}_i$ to ${\widetilde x}_{i+1}$ with the
corresponding segment in ${\widetilde \beta}$. Then, $D_i$ is the
component of $\C-{\widetilde \alpha}-{\widetilde \beta}$ which is
inside $C_i$, and whose boundary contains ${\widetilde x}_i$.  For
example, as we see from the figure, $C_1$, whose interior represents a
flow from $\xCov_1$ to $\xCov_2$, can be decomposed as a sum
$D_1+D_3$; $C_2$, which gives a flow from $\xCov_3$ to $\xCov_2$,
bounds simply $D_2$; while $C_3$ bounds $D_3$, giving a flow from
$\xCov_3$ to $\xCov_4$. It follows that the domain of the Whitney disk
from $\xCov_4$ to $\xCov_1$ is given by $D_1+D_2$. This disk has
$\Mas=1$, and it clearly supports a single flowline (indeed, it is
represented by an embedded disk in the plane). Moreover, the domain
contains no points point from the $W$ lattice and one from the $Z$ lattice,
so we have determined that the ${\widetilde x}_4$ coefficient of
$\partial [{\widetilde x}_1,i,j]$ is $\pm [{\widetilde x}_4,i,j-1]$.

\begin{figure}
\mbox{\vbox{\epsfbox{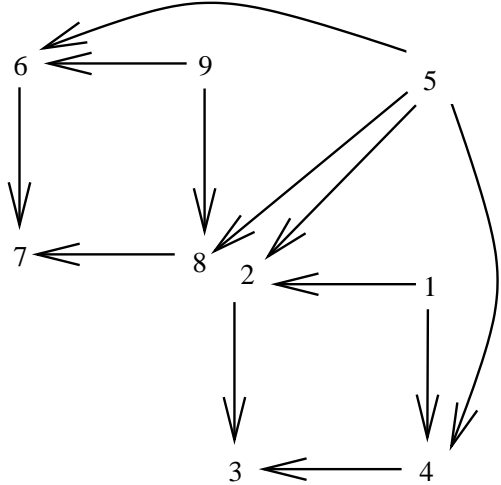}}}
\caption{\label{fig:PictureDiff}
{\bf{The chain complex.}}  This depicts the differential for
$\CFKinf(S^3,K)$. The integers here correspond to intersection points
with the labeling conventions from the text: i.e. the integer $k$
corresponds to ${\widetilde x}_k$. An arrow denotes a non-trivial
differential, and (rough) coordinates in the plane correspond to the
indices $i,j$. Of course, the complex $\CFKinf(S^3,K)$ consists of
infinitely many copies of this complex; more precisely 
it is the tensor product of this complex 
with $\Z[U,U^{-1}]$.}
\end{figure}

Proceeding in the like manner, one 
verifies that $\CFKinf(S^3,K)$ 
(with the notational shorthand from Subsection~\ref{subsec:Shorthand})
is generated as a $\Z$-module by generators
indexed by $i\in\Z$ of the form:
\begin{equation}
\label{eq:TheGenerators}
\begin{array}{lll}
[{\widetilde x}_1, i+1, i], &
[{\widetilde x}_2, i, i], &
[{\widetilde x}_3, i, i-1],\\ 
{[{\widetilde x}_4, i+1, i-1]}, &
[{\widetilde x}_5, i+1, i+1], &
[{\widetilde x}_6, i-1, i+1],\\ 
{[{\widetilde x}_7, i-1, i]}, &
[{\widetilde x}_8, i, i], &
[{\widetilde x}_9, i, i+1];
\end{array}
\end{equation}
and the boundary operator on $\CFKinf(S^3,K)$ (up to some signs) is given by:
\begin{eqnarray*}
\partial[{\widetilde x}_1,i+1,j]&=& [\xCov_2,i,i]+[\xCov_4,i+1,i-1] \\
\partial[\xCov_2,i,i]&=&[\xCov_3,i,i-1] \\
\partial[\xCov_3,i,i-1]&=& 0 \\
\partial[\xCov_4,i+1,i-1]&=& [\xCov_3,i,i-1] \\
\partial[\xCov_5,i+1,i+1]&=&
[\xCov_2,i,i]+[\xCov_4,i+1,i-1]+[\xCov_6,i-1,i+1]+[\xCov_8,i,i] \\
\partial[\xCov_6,i-1,i+1]&=&[\xCov_7,i-1,i] \\
\partial[\xCov_7,i-1,i]&=& 0 \\
\partial[\xCov_8,i,i]&=&[\xCov_7,i-1,i] \\
\partial[\xCov_9,i,i+1] &=& [\xCov_6,i-1,i+1]+[\xCov_8,i,i]
\end{eqnarray*}
In calculations, it is often useful to depict this schematically, 
as pictured in Figure~\ref{fig:PictureDiff}.

Consider the complex $\CFK^{0,*}(Y,K)$, which is generated by the nine
generators as enumerated in Equation~\eqref{eq:TheGenerators}, but
with $i=0$. We obtain a complex whose homology is $\Z$, supported in a
single dimension. In fact, this homology is represented by the cycle
$[\xCov_1,0,-1]\pm [\xCov_5,0,0]$. Indeed (as in to
Lemma~\ref{lemma:AbsGrade}) this generator represents the Floer
homology class of ${\HFa}(S^3)\cong\Z$, so its absolute grading is
zero.

In particular, we have the following:

\begin{prop}
\label{prop:NineFortyTwo}
Let $K$ denote the knot $9_{42}$ (illustrated in Figure~\ref{fig:BigKnot}).
Various invariants for the knot $K$ are given by :
\begin{eqnarray*}
\HFKa(S^3,K,i)&\cong&
\left\{\begin{array}{ll}
\Z_{(3)} & {\text{if $i=2$}}  \\ \\
\Z^2_{(2)} & {\text{if $i=1$}} \\ \\
\Z^2_{(1)} \oplus \Z_{(0)} & {\text{if $i=0$}} \\ \\
\Z^2_{(0)} & {\text{if $i=-1$}} \\ \\
\Z_{(-1)} & {\text{if $i=-2$}} \\ \\
0 & {\text{otherwise.}}  
\end{array}
\right. \\ \\
\HFp(S^3_0(K),i)&\cong&
\left\{\begin{array}{ll}
{\mathcal T}_{-1/2}\oplus{\mathcal T}_{1/2}
& {\text{if $i=0$}} \\ \\
\Z & {\text{if $i=\pm 1$}} \\ \\
0 & {\text{otherwise}}
\end{array}
\right.
\end{eqnarray*}
In particular, $d_{\pm 1/2}(S^3_0(K))=\pm 1/2$.
\end{prop}

\begin{proof}
The calculation of $\HFK(S^3,K)$
follows from inspection of the chain complex.

For instance, observe that $\CFKa(S^3,K,2)$ is generated by a single
element $[{\widetilde x}_6,0,2],$, while $\CFa(S^3,K,0)$ is generated by
three elements $[{\widetilde x}_2,0,0]$, $[{\widetilde x}_5,0,0]$, and
$[{\widetilde x}_8,0,0]$.

To calculate $\HFp(S^3_0(K))$, we first use
Theorem~\ref{thm:LargeNegSurgeries} to identify
$\HFp(S^3_{-p}(K),[s])$ for large $p$ with $H_*(\bCFKp(Y,K,s))$ (with
the grading shift). The graded long exact squence for integral
surgeries then shows that $\HFp(S^3_0(K))$ has the claimed form.
\end{proof}

\section{Connected sums of knots}
\label{sec:ConnectedSums}

Recall that there are multiplication maps
$$
\begin{array}{lll}
\CFa(Y_1)\otimes_{\Z} \CFa(Y_2)&\longrightarrow&  \CFa(Y_1\# Y_2), \\
\CFm(Y_1)\otimes_{\Z[U]} \CFm(Y_2)&\longrightarrow& \CFm(Y_1\# Y_2), \\
\CFinf(Y_1)\otimes_{\Z[U,U^{-1}]} \CFinf(Y_2)&\longrightarrow& \CFinf(Y_1\# Y_2), 
\end{array}
$$
which induce homotopy equivalences. 

To state the generalization for knot invariants, we use the following
convention: the tensor product of two $\Z$-filtered
(resp. $\Z\oplus\Z$-filtered) filtered complexes $(C_1,\Filt_1)$ and
$(C_2,\Filt_2)$ is naturally a $\Z$-filtered (resp. $\Z\oplus
\Z$-filtered) complex $(C_1\otimes C_2,\Filt_{\otimes})$, where the
tensor product filtration is defined by $$\Filt_{\otimes}(x_1\otimes
x_2)=\Filt_1(x_1)+\Filt_2(x_2)$$ for arbitrary homogeneous generators
$x_1\in C_1$ and $x_2 \in C_2$.

Given $\relspinct_1\in\RelSpinC(Y_1,K_1)$
and $\relspinct_2\in\RelSpinC(Y_2,K_2)$,
let $\relspinct_3=\relspinct_1\#\relspinct_2\in\RelSpinC(Y_1\#Y_2,K_1\#K_2)$
be characterized by the properties that 
\begin{itemize}
\item if $\spinc_i\in\SpinC(Y_i)$ denote the $\SpinC$ structures underlying
$\relspinct_i$ for $i=1,...,3$, where $Y_3=Y_1\#Y_2$, then $\spinc_3$ is
cobordant to $\spinc_1\cup\spinc_2$ under the natural cobordism from $Y_3$ to
$Y_1\cup Y_2$
\item if $F_1$ and $F_2$ are Seifert surfaces for $K_1$ and $K_2$, and $F_3$
is the corresponding Seifert surface for $K_1\#K_2$ obtained by boundary
connected sum,
then
$$\langle c_1(\relspinct_3),[F_3]\rangle =
\langle c_1(\relspinct_1),[F_1]\rangle +
\langle c_1(\relspinct_2),[F_2]\rangle.
$$
\end{itemize}

\begin{theorem}
\label{thm:ConnectedSumsOfKnots}
There are filtered chain homotopy equivalences 
$$
\begin{array}{lll}
\CFK^{0,*}(Y_1,K_1)\otimes_\Z \CFK^{0,*}(Y_1,K_1) 
&{\longrightarrow}& 
\CFK^{0,*}(Y_1\# Y_2,K_1\#K_2) \\
\CFK^{<,*}(Y_1,K_1)\otimes_{\Z[U]}\CFK^{<,*}(Y_1,K_1) 
&{\longrightarrow}& 
\CFK^{<,*}(Y_1\# Y_2,K_1\#K_2) \\
\CFKinf(Y_1,K_1)\otimes_{\Z[U,U^{-1}]} \CFKinf(Y_1,K_1) 
&{\longrightarrow}& 
\CFKinf(Y_1\# Y_2,K_1\#K_2)
\end{array}
$$
which respect the splittings of the complexes according to their summands;
in particular,
$$\CFKinf(Y_1\#Y_2,\relspinct_3)
\cong
\bigoplus_{\{\relspinct_1,\relspinct_2\in\RelSpinC(Y_1,K_1)\times\RelSpinC(Y_2,K_2)\big|\relspinct_1\#\relspinct_2=\relspinct_3\}}\CFKinf(Y_1,\relspinct_1)\otimes
\CFKinf(Y_2,\relspinct_2).$$
\end{theorem}

\begin{figure}
\mbox{\vbox{\epsfbox{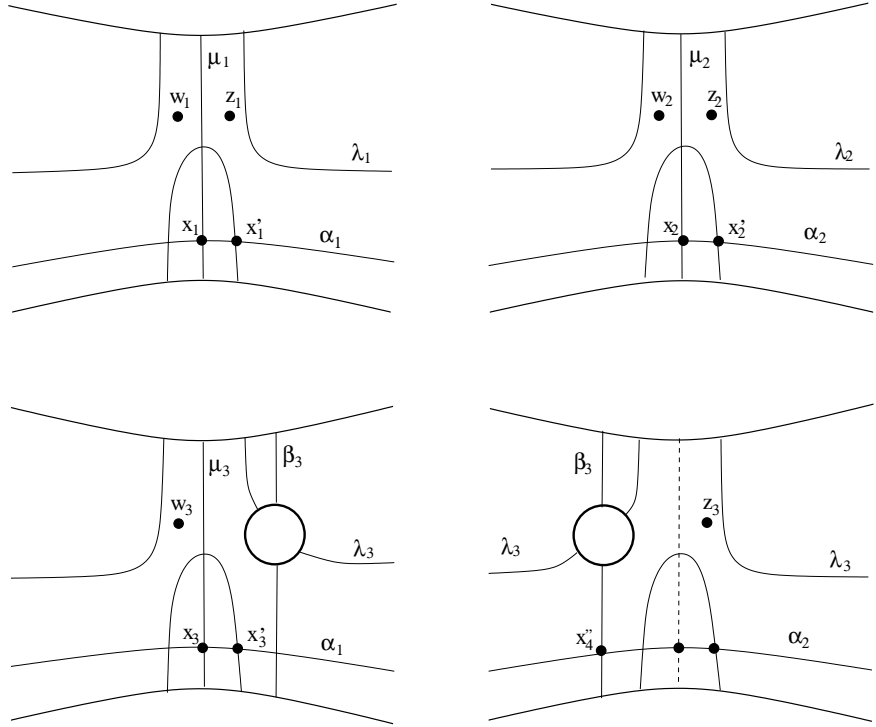}}}
\caption{\label{fig:ConnectedSumKnots}
{\bf{Connected sum of knots: before and after.}} The top picture
represents the (disconnected) union of two original manifolds
$(Y_1,K_1)$ and $(Y_2,K_2)$, while the second represents their
connected sum (in particular, in the second picture, the two open
circles are identified). The curve $\mu_3$ (corresponding to $\mu_1$
in the previous diagram) represents a meridian for $K_1\# K_2$ while
the curve $\beta_2$ which is obtained by handlesliding $\mu_2$ over
$\mu_1$ (after forming the connected sum) is a new attaching circle for
$Y_3=Y_1\# Y_2$. The circle $\lambda_3$ is the longitude for $K_1\#
K_2$. Note that $\beta_3$ and $\lambda_3$ are the only curves which go
through the connected sum neck. The curve $\beta_3$ is easily seen not
to appear in the boundary of the periodic domain for $K_1\# K_2$,
while the curve $\lambda_3$ appears in it with multiplicity one.  The
pair $x_3'\times x_4''$ corresponds under the tensor product map to
$x_1\times x_2$.}
\end{figure}

\begin{proof}
We sketch here the definition of the usual connected sum map (see
Section~\ref{HolDiskTwo:sec:ConnectedSums} of~\cite{HolDiskTwo}).  Let
$(\Sigma_1,\alphas_1,\betas_1,\StartPt_1)$ and
$(\Sigma_2,\alphas_2,\betas_2,\StartPt_2)$ be pointed Heegaard diagrams
representing $Y_1$ and $Y_2$.  There are maps
$$\CFinf(Y_1)\longrightarrow 
\CFinf(Y_1\#^{g_2}(S^2\times S^1))$$
$$\CFinf(Y_2)\longrightarrow \CFa(\#^{g_1}(S^2\times S^1)\# Y_2)$$
sending $\x_1$ and $\x_2$ to $\x_1\times\Theta_1$ and $\Theta_2\times
\x_2$ respectively, where $\Theta_1$ and $\Theta_2$ are intersection
points representing the canonical (top-dimensional) $\HFleq$ class for
the connected sum of $S^2\times S^1$.  We then compose this with a
count $f$ of holomorphic triangles in the Heegaard triple
$(\Sigma_1\#\Sigma_2,\alphas_1\times \alphas_2,
\betas_1\times\alphas_2, \betas_1\times\betas_2)$.
In the above expression, any time some $g$-tuple of circles is repeated,
it is to be understood that one takes a small exact Hamiltonian
translate to make the intersections transverse. 

The refinement
$$f^\infty\colon  \CFKinf(Y_1,K_1)\otimes_{\Z[U,U^{-1}]}\CFKinf(Y_2,K_2)
\longrightarrow
\CFKinf(Y_1\#Y_2,K_1\#K_2)
$$
is defined as follows. Consider the map
\begin{eqnarray*}
\lefteqn{f([\x_1,i_1,j_1]\otimes [\x_2,i_2,j_2])}\\
&=& 
\sum_{\y\in{\mathbb T}_{\alpha_1\times\alpha_2}\cap
{\mathbb T}_{\beta_1\times\beta_2}}
\sum_{\{\psi\in\pi_2(\x_1\times\Theta_2,\Theta_1\times\x_2,\y)\}} 
\#\ModFlow(\psi)\cm [\y,i_1+i_2-n_{\BasePt_1}(\psi),
j_1+j_2-n_{\FiltPt_2}(\psi)],
\end{eqnarray*}
where we continue with the notation from above, with the understanding
that the meridians $\mu_1$ and $\mu_2$ appear in the $g$-tuples
$\betas_1$ and $\betas_2$ respectively. Of course, the above map is a
perturbation of the map sending $[\x_1,i_1,j_1]\otimes [\x_2,i_2,j_2]$
to $[\x_1\times \x_2,i_1+i_2,j_1+j_2]$. The map clearly respects
filtrations, since $n_{\EndPt_2}(\psi)\geq 0$ whenever $\psi$ has a
holomorphic representative. We go from this Heegaard diagram to a
Heegaard diagram ro the following Heegaard diagram for the knot
$(Y_1\# Y_2,K_1\#K_2)$, 
$$(\Sigma,\alphas_1\cup\alphas_2,\beta_3\cup
(\betas_1-\mu_1)'\cup(\betas_2-\mu_2)', \mu_3),$$
where here the primes
denote small exact Hamiltonian translates, $\mu_3$ is a curve obtained
from $\mu_1$ by a small exact Hamiltonian translate, and $\beta_3$ is
obtained as a handleslide of $\mu_2$ over $\mu_1$ -- in particular,
the new Heegaard diagram for $Y_1\# Y_2$ differs from the 
$(\Sigma_1\#\Sigma_2,\alphas_1\times\alphas_2,\betas_1\times\betas_2)$
by a single handelside. 
Moreover, it is easy to see that the composite map here respects relative filtrations.

We verify that this map respects the splitting according to $\RelSpinC(Y_i,K_i)$.
This is clear on the level of $\SpinC(Y_i)$; thus, it remains to show
that if $F_1$ and $F_2$ are compatible
Seifert surfaces for $K_1$ and $K_2$, then
\begin{equation}
\label{eq:SumOfSpinC}
\langle c_1({\relspinct}(\x_1)),[{\widehat F_1}]\rangle 
+ \langle c_1({\relspinct}(\x_2)),[{\widehat F_2}]\rangle 
=
\langle c_1({\widehat \relspinct}(\x_1\times \x_2)),
[{\widehat {F_1\#_b F_2}}]\rangle, 
\end{equation}
where $F_1\#_b F_2$ is the Seifert surface for $K_1\#K_2$ inside
$Y_1\#Y_2$ obtained by boundary connected sum.  To this end, observe
that if $\lambda_1$ and $\lambda_2$ are curves in the Heegaard
diagrams for $Y_1$ and $Y_2$ representing longitudes for the knots
$K_1$ and $K_2$ respectively, then $\lambda_1\#\lambda_2$ represents a
longitude for $K_1\# K_2$. Correspondingly, if $P_1$ and $P_2$ are the
periodic domains (corresponding to ${\widehat F}_1$ and ${\widehat
F}_2$) in $(Y_1)_0(K_1)$ and $(Y_2)_0(K_2)$ respectively, then the
periodic domain for $F_1\#_b F_2$ in $(Y_1\# Y_2)_0(K_1\#K_2)$ is
obtained from $P_1$ and $P_2$ by a boundary connected sum (see
Figure~\ref{fig:ConnectedSumKnots}). Under this operation, the Euler
measure satisfies $$\chi(P_1\#_b P_2)=\chi(P_1)+\chi(P_2)-1,$$ and it
is easy to see that if $\x_3$ denotes the intersection point
corresponding to $\x_1\times \x_2$, then $$\Mass_{\x_3}(P_1\#_b P_2)=
\Mass_{\x_1}(P_1)+\Mass_{\x_2}(P_2)+\OneHalf.$$ (Again, the reader is asked to consult
Figure~\ref{fig:ConnectedSumKnots}.)  Equation~\eqref{eq:SumOfSpinC}
now follows immediately from these observations and
Equation~\eqref{eq:MeasureCalc}.

Having defined $f^{\infty}$, it is now easy to adapt the arguments 
for three-manifolds to verify that $f^{\infty}$ is a chain homotopy equivalence
of filtered complexes. (Moreover, the corresponding constructions for
$\CFK^{<,*}$ and
$\CFK^{0,*}$ proceed in the same way.)
\end{proof}

Note that the above theorem has the following immediate corollary
(which can be thought of as a refinement of the multiplicativity of
the Alexander polynomial under connected sums):

\begin{cor}
\label{cor:Kunneth}
Let $\relspinct_1\in\RelSpinC(Y_1,K_1)$ and $\relspinct_2\in\RelSpinC(Y_2,K_2)$ 
be $\SpinC$ structures,  and let $\relspinct_1\#\relspinct_2\in\RelSpinC(Y_1\# Y_2,K_1\#K_2)$ denote
the $\SpinC$ structure whose
underlying $\SpinC$ structure is cobordant to $\spinc_1\cup\spinc_2$
over $Y_1\cup Y_2$ (under the natural cobordism to $Y_1\#Y_2$), and
$$\langle c_1({\relspinct_1\#\relspinct_2}),[{\widehat{F_1 \#_b F_2}}]\rangle = 
\langle c_1(\relspinct_1),[{\widehat{F_1}}]\rangle +
\langle c_1(\relspinct_2),[{\widehat{F_2}}]\rangle.$$
Then, for each $\relspinct_3\in\RelSpinC(Y_1\#Y_2,K_1\#K_2)$, we have
$$\HFKa(Y_1\# Y_2,K_1\#K_2,\relspinct_3)\cong 
\bigoplus_{\relspinct_1\#\relspinct_2=\relspinct_3}
H_*(\CFKa(Y_1,K_1,\relspinct_1)\otimes 
\CFKa(Y_2,K_2,\relspinct_2)).$$
\end{cor}

\begin{proof}
This follows immediately from the theorem.
\end{proof}
\section{Long exact sequences}
\label{sec:Sequences}

One can readily adapt the surgery long exact sequences
from~\cite{HolDiskTwo} to the context of knot invariants, in several ways.
We give here a discussion whose level of generality is suitable for the
applications to follow. 

Recall that if $\gamma$ is a framed knot in a three-manifold $Y$, and
$\spinc$ is a $\SpinC$ structure on the four manifold $W=W_\gamma(Y)$,
then there are maps induced by counting holomorphic triangles
$$\begin{array}{lll}
{\widehat f}_{W,\spinc}\colon\CFa(Y) &\longrightarrow& \CFa(Y_\gamma), \\
f^-_{W,\spinc}\colon\CFm(Y) &\longrightarrow& \CFm(Y_\gamma), \\
f^\infty_{W,\spinc}\colon\CFinf(Y) &\longrightarrow& \CFinf(Y_\gamma).
\end{array}
$$ 

In the presence of a Seifert surface $F$ for a knot, 
for each $\spinc\in\SpinC(Y)$, there is a
$\relspinct_0\in\RelSpinC(Y,K)$ extending $\spinc$, which is
uniquely characterized by
the property that 
$$\langle c_1(\relspinct_0),[{\widehat F}]\rangle =0.$$
Thus, we can identify,
$$\CFinf(Y,\spinc)\cong\CFKinf(Y,\relspinct_0),$$ 
as $\Z$-filtered complexes; and this isomorphism endows
$\CFinf(Y,\spinc)$ with an extra $\Z$ grading
(i.e. our representatives for $\CFKinf$ are triples $[\x,i,j]$,
and the $\Z$ filtration is given by the value of $j$).
We sometimes suppress $\relspinct_0$ from the notation, writing
$\CFKinf(Y,K,\spinc)$ when $F$ is understood from the context.
(Note that this shorthand coincides with the corresponding shorthand
we used in the case where $c_1(\spinc)$ is torsion.)

Suppose now that $K$ is an oriented knot in $Y$ which is disjoint from
$\gamma$, and whose linking number with $\gamma$ is zero. In this
case, we can choose a Seifert surface $F$ for $Y$ which is disjoint
from $\gamma$; thus, this can be closed up in $Y$ or $Y_\gamma$ to
obtain surfaces ${\widehat F}$ and ${\widehat F}'$.

\begin{prop}
\label{prop:RespectFiltration}
Let $Y$ be a three-manifold with a null-homologous knot $K$, and let
$\gamma$ be a framed knot which is disjoint from $K$, and whose
linking number with $K$ is zero, and let $Y'$ be the
three-manifold induced from $Y$ by surgery along $\gamma$ and
$K'$ be the induced knot. Then, for each $\SpinC$ structure
$\spinc$, the map $f_{W,\spinc}$ induced by counting holomorphic
triangles respects the above induced $\Z$-filtrations.
\end{prop}

\begin{proof}
We describe here the case of $\CFinf(Y)$; the other two cases follow similarly.
Let $(\Sigma,\alphas,\betas,\gammas,\StartPt,\EndPt)$ be a
Heegaard triple, where for $i=2,...,g$, $\beta_i$ is a small
exact Hamiltonian translate of $\gamma_i$, with $\mu=\beta_2$, 
while the $\beta_1$ and $\gamma_1$ are different (the $\gamma_1$ here
corresponds to the longitude of $\gamma$, and $\beta_1$ to its meridian)
Following~\cite{HolDiskFour}, the map 
$$f^\infty_{W,\spinc}\colon \CFKinf(Y,K) \longrightarrow \CFKinf(Y',K')$$
induced by the cobordism is defined by
$$f^\infty_{W,\spinc}([\x,i,j])
=
\sum_{\{\y\in\Ta\cap\Tc\}}
\sum_{\{\psi\in\pi_2(\x,\Theta_{\beta,\gamma},\y)\big|
\Mas(\psi)=0\}}
\#\ModFlow(\psi)\cm
[\y,i-n_{\BasePt}(\psi),j-n_{\FiltPt}(\psi)].  $$ 

We must verify that the range of this map is $\CFKinf(Y',K')$ as
claimed, i.e. that if 
\begin{equation}
\label{eq:WhatWeWant}
\langle c_1(\relspinc_\Mark(\x),[F]\rangle
+(i-j) =0,
\end{equation} then if $\psi\in\pi_2(\x,\Theta,\y)$, we also have that
$$\langle c_1(\relspinc_\Mark(\y)),[{\widehat F}']\rangle +
(i-n_\BasePt(\psi)-j+n_\FiltPt(\psi))=0,$$ as well. 

As a preliminary remark, we claim that if $P$ is a triply-periodic
domain for the Heegard triple describing the cobordism from $Y$ to
$Y'$, then $n_\BasePt(P)=n_\FiltPt(P)$. The reason for this is
that $P$ can be viewed as a null-homology of $\gamma$ in $Y$ (as the
only curve which appears in $\partial P$ which is not isotopic to an
attaching circle for $Y$ is a longitude for $\gamma$). Thus,
the multiplicity with which $\beta_2$ appears in $\partial P$ 
can be interpreted as the linking number of $\gamma$ with $\mu$,
which we assumed to vanish. But this multiplicity is given by
$n_\BasePt(P)-n_\FiltPt(P)$.

In view of these remarks, it suffices to verify
that 
\begin{equation}
\label{eq:SpecialCase}
\langle c_1(\spincrel_\Mark(\x)),[F]\rangle 
= \langle c_1(\spincrel_\Mark(\y)),[F']\rangle, 
\end{equation}
for two particular intersection
points $\x$, $\y$, which can be connected by a triangle
$\psi\in\pi_2(\x,\Theta,\y)$ with $n_\BasePt(\psi)=n_\FiltPt(\psi)=0$
(by appealing to Lemma~\ref{lemma:RelSpinCThree}).
Indeed, it suffices to consider the case where the the point in $\x$
lying on $\beta_2$ is a small perturbation of the point on $\y$ lying
on $\gamma_2$ (recall here that $\gamma_2$ is a small perturbation of
$\beta_2$). It is easy to see that $\psi$ induces a triangle
$\psi'\in\pi_2(\x',\Theta,\y')$, and hence a $\SpinC$ structure on the
cobordism from $Y_0(K)$ to $Y_0'(K')$ whose restrictions to $Y_0(K)$
and $Y_0(K')$ are $\spinc_\BasePt(\x')$ and $\spinc_\BasePt(\y')$
respectively. Since $[{\widehat F}]$ and $[{\widehat F'}]$ are
homologous in this cobordism (again, we are using the fact that the
linking number of $\gamma$ with $K$ vanishes),
Equation~\eqref{eq:SpecialCase} follows at once, and hence also
Equation~\eqref{eq:WhatWeWant}.

As defined, the map clearly respects filtrations.
\end{proof}

We have now the following easy consequence of these observations, as
combined with the long exact sequences of~\cite{HolDiskTwo}. Note that
in the following statement, 
we let
$$\HFKa(Y,K,m)=
\bigoplus_
{\{\relspinct\in\RelSpinC(Y,K)\big|\langle c_1(\relspinct),[{\widehat F}]\rangle = 2m\}}\HFKa(Y,K,\relspinct).$$

\begin{theorem}
\label{thm:LongExactSeq}
Let $(Y,K)$ be a knot, and $\gamma$ be a framed knot in $Y-K$.  
Then, by choosing all Seifert surfaces compatibly, 
the counts of holomorphic triangles give a long exact sequence
for each integer $m$:
$$
\begin{CD}
... @>>> \HFKa(Y_{-1}(\gamma), K,m)
@>{f_1}>> \HFKa(Y_{0}(\gamma),K,m) @>{f_2}>>  \HFKa(Y,K,m)@>{f_3}>>...
\end{CD}
$$
In cases where the groups are graded (e.g. when $Y$ is a homology sphere),
the maps $f_1$ and $f_2$ are homogeneous maps which
lower the absolute gradings by $\OneHalf$,
while $f_3$ is a sum of homogeneous maps, which are non-increasing
in the absolute grading.
\end{theorem}

\begin{proof}
This follows easily from inspecting the proof of
the surgery long exact sequence for $\CFKa$ (see
Theorem~\ref{HolDiskTwo:thm:GeneralSurgery} of~\cite{HolDiskTwo}),
and observing that the homotopy equivalences also respect the filtrations, 
following arguments from Proposition~\ref{prop:RespectFiltration} above.

Recall that the maps $f_1$, $f_2$, and $f_3$ are sums of maps induced
by $\SpinC$ structures over the
corresponding cobordisms. The statements about gradings now follow
immediately from the dimension formula of~\cite{HolDiskFour},
stated in Equation~\eqref{eq:DimensionFormula} above.
\end{proof}

It is not difficult to see that the above theorem generalizes the
skein exact sequence stated in the introduction. We return to this
point in Section~\ref{sec:Links}.

\section{Some fibered examples}
\label{sec:FiberedExamples}

As an illustration of some of the above techniques, we give some
calculations of knot homology groups for some simple (genus one)
fibered knots.  These calculations, together with basic properties of
the knot invariants, allow us to calculate $\HFp(S^1\times
\Sigma_g,\spinct)$ (and also some other three-manifolds which fiber
over the circle), where $\spinct$ is some non-torsion $\SpinC$
structure.  It is interesting to compare these with Seiberg-Witten
results, c.f.~\cite{MSzT} and~\cite{MunozWang}, 
with the quantum cohomology calculations of~\cite{ThadBert},
and also with the
``periodic Floer homology'' of Hutchings,
see~\cite{HutchingsSullivan}.

Consider the three-component Borromean rings in $S^3$, and fix $m,n\in
\Z\cup\{\infty\}$ (indeed, we will typically consider
$m,n\in\{0,1,\infty\}$).  Perform $m$-surgery on one component,
$n$-surgery on another, and let $B(m,n)$ denote the remaining
component, thought of as knot in the surgered three-manifold. Thus,
$B(m,n)$ is a null-homologous knot inside $L(m,1)\# L(n,1)$ (with the understanding
that $L(0,1)\cong S^2\times S^1$ and $L(1,1)\cong S^3$).  By Kirby
calculus, one can see that $B(1,1)\subset S^3$ is diffeomorphic to the
right-handed trefoil knot, and $B(m,\infty)\subset L(m,1)$ is an unknot. 

In general, $L(m,1)\#L(n,1)-B(m,n)$ is diffeomorphic to the mapping
torus of an automorphism of the torus with a disk removed
$T^2-\CDisk$. In particular, the complement $(S^2\times S^1)-B(0,0)$ is
diffeomorphic to the product of a circle with $T^2-\CDisk$, and
$(S^2\times S^1)-B(1,0)$ is the mapping torus of a non-separating Dehn
twist acting on $T^2-\CDisk$.

\begin{figure}
\mbox{\vbox{\epsfbox{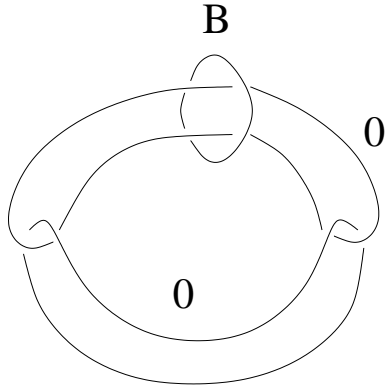}}}
\caption{\label{fig:BorromeanKnot}
A Kirby calculus picture of the ``Borromean knot'' $B\subset\#^2(S^2\times S^1)$.}
\end{figure}

\begin{remark}
In the following three-manifolds $Y$, there is a unique
$\spinct\in\SpinC(Y)$ with trivial first Chern class (and indeed, 
it is clear that $\HFa(Y,\spinc)=0$ for all other $\SpinC$
structures).  Thus, there is no ambiguity in writing $\HFKa(Y,K,j)$
to denote $\HFKa(Y,K,\relspinct_j)$, where
$\relspinct_j\in\SpinC(Y,K)$ is characterized by the fact that it
extends $\spinct$, and satisfies $\langle
c_1(\relspinct_j),[{\widehat F}]\rangle=2j$.
\end{remark}

\begin{prop}
\label{prop:KnotHomology}
The knot homology groups $\HFKa$ of Borromean knots are given by:
\begin{eqnarray*}
\HFKa(S^3,B(1,1),j)&\cong& \left\{
\begin{array}{ll}
\Z_{(0)} & {\text{if $j=1$}} \\
\Z_{(-1)} & {\text{if $j=0$}} \\
\Z_{(-2)} & {\text{if $j=-1$}} \\
0 & {\text{otherwise}} 
\end{array}\right. \\
\HFKa(S^2\times S^1,B(0,1),j)&\cong& \left\{
\begin{array}{ll}
\Z_{(\OneHalf)} & {\text{if $j=1$}} \\
\Z^2_{(-\OneHalf)} & {\text{if $j=0$}} \\
\Z_{(-3/2)} & {\text{if $j=-1$}} \\
0 & {\text{otherwise}} 
\end{array}\right. \\
\HFKa(\#^2(S^2\times S^1),B(0,0),j)&\cong& \left\{
\begin{array}{ll}
\Z_{(1)} & {\text{if $j=1$}} \\
\Z_{(0)}^2 & {\text{if $j=0$}} \\
\Z_{(-1)} & {\text{if $j=-1$}} \\
0 & {\text{otherwise}} 
\end{array}\right.
\end{eqnarray*}
\end{prop}

\begin{proof}
The first isomorphism follows immediately from the trefoil calculation
from Subsection~\ref{subsec:Trefoil}. 

The calculation of $\HFKa$ for $B(0,1)$ could be made in the same spirit -- by reducing
to a genus one Heegaard diagram for $S^2\times S^1$. Alternatively, one
can use the long exact sequence of Theorem~\ref{thm:LongExactSeq}:
$$
\tiny{\begin{CD}
\HFKa(S^3,B(1,\infty),j)@>{f_1}>>\HFKa(S^2\times S^1, B(1,0),j)@>{f_2}>>\HFKa(S^3,
B(1,1),j)@>{f_3}>>...
\end{CD}}
$$
Since $B(1,\infty)\subset S^3$ is an unknot,
$\HFKa(S^3,B(1,\infty),j)=0$ when $j\neq 0$. Moreover the map $f_2$
lowers degree by $1/2$ (c.f. Equation~\eqref{eq:DimensionFormula}). 
Thus, the calculation of $\HFKa(S^2\times S^1,B(1,0),j)$ for $j\neq 0$ is complete.

In the case where $j=0$, of course,
\begin{eqnarray*}
\HFKa(S^3,B(1,\infty),0)\cong \Z_{(0)}
&{\text{and}}&
\HFKa(S^3,B(1,1),0)\cong \Z_{(-1)},
\end{eqnarray*} while the map $f_3$ is non-increasing in grading.
Thus, it must vanish identically, and, since $f_1$ also lowers degree
by $\OneHalf$, the computation of $\HFKa(S^3,B(1,0),j)$ for all $j$ is
complete.

For $B(0,0)$, we now have
\begin{equation}
\label{eq:B00}
\tiny{
\begin{CD}
\HFKa(S^2\times S^1,B(0,\infty),j)@>{f_1}>>
\HFKa(\#^2(S^2\times S^1), B(0,0),j)@>{f_2}>>\HFKa(S^2\times S^1,
B(0,1),j)@>{f_3}>>...
\end{CD}}
\end{equation}
Thus, the case where $j\neq 0$ follows exactly as before. 

Before continuing the calculation when $j=0$, it is useful to make some
observations. Specifically, recall that 
$$\bigoplus_{j}\HFKa(\#^2(S^2\times S^1),B(0,0),j)$$
is the homology of the graded object associated to a filtration of
$\HFa(\#^2(S^2\times S^1))$. Thus, since 
$\HFKa(\#^2(S^2\times S^1),B(0,0),j)=0$ for $j<-2$, we obtain a natural map
$$\Z_{(-1)}\cong \HFKa(\#^2(S^2\times S^1),B(0,0),-1)\longrightarrow \HFa(\#^2(S^2\times S^1))
\cong \Z_{(-1)}\oplus\Z^2_{(0)}\oplus \Z_{(1)}.$$
We claim that this map actually induces an isomorphism onto the $\Z_{(-1)}$
factor on the right-hand-side. This follows immediately from the following
commutative diagram:
$$
\begin{CD}
\HFKa(S^3,B(-1,-1),-1) @>{\cong}>> \HFa(S^3)\cong \Z_{(0)} \\
@V{f}VV @V{g}VV \\
\HFKa(\#^2(S^2\times S^1),B(0,0),-1)@>{i}>> \HFa(\#^2(S^2\times S^1)),
\end{CD}
$$ where here the two horizontal maps are induced by inclusions, the
fact that the first is an isomorphism follows readily from the
calculation of $B(1,1)$, while the vertical maps $f$ and $g$ are
induced by the obvious cobordisms (here $g$ is induced by the
two-handle additions, and $f$ is the induced map on the knot complex).
In particular, it is a straightforward calculation that
(c.f.~\cite{HolDiskFour}) that $g$ induces an isomorphism to
$\HFa_{(-1)}(\#^2(S^2\times S^1))\cong \Z$. It follows that $i$ is an
isomorphism for degree $-1$, as claimed.

We return now to the calculation of $\HFKa(\#^2(S^2\times S^1),j)$
when $j=0$.  Recall that $\HFKa(S^2\times S^1,B(0,\infty),0)\cong
\Z_{(-1/2)}\oplus \Z_{(1/2)}$, and of course $\HFKa(S^2\times
S^1,B(0,1),j)\cong \Z^2_{(-1/2)}$. Verifying the calculation thus
reduces to showing that $f_3$ (from Long Exact
Sequence~\eqref{eq:B00}) surjects onto the $\Z_{(-1/2)}$ summand of
$\HFKa(S^2\times S^1,B(0,\infty),0)$.  If it did not, the map would
have to have cokernel, and hence the rank (over possibly some finite
field) of $$G=\bigoplus_{j}\HFKa_{-1}(\#^2(S^2\times S^1), B(0,0),j)$$
would be two.  But this is impossible: we have already shown that
$\HFKa(\#(S^2\times S^1),B(0,0),-1)$ represents a generator for
$\HFa_{-1}(\#^2(S^2\times S^1))$. It also follows from the previous
calculation that $\HFKa(\#^2(S^2\times S^1))$ is supported in
dimension $1$, so the remaining generator must come from
$\HFKa_{-1}(\#^2(S^2\times S^1),0)$. Indeed, by dimension reasons,
this generator cannot be killed by any higher differential (in the
spectral sequence onnecting $\HFKa(\#^2(S^2\times S^1),B(0,0))$ to
$\HFa(\#^2(S^2\times S^1))$),
and hence contradicting the fact that
$\HFa_{-1}(\#^2(S^2\times S^1))\cong \Z$.
\end{proof}

We would like to use the above calculations to calculate $\HFp$ of a
number of three-manifolds. Indeed, the manifolds we will consider have
positive first Betti number, and in this case, it is interesting to
obtain information of Floer homology as a module over the ring
$\Z[U]\otimes_\Z \Wedge^* H_1(Y;\Z)/\Tors$
(c.f. Subsection~\ref{HolDiskOne:subsec:DefAct} of~\cite{HolDisk}). We pause now
to recall the construction of this action.

Given a curve $\gamma$ in a Heegaard surface for $Y$ which
is disjoint from all the $\alpha_i\cap\beta_j$, the action of the corresponding
class $[\gamma]\in H_1(Y;\Z)$ is induced by the chain map
$$A_\gamma\colon 
\CFinf(Y)\longrightarrow \CFinf(Y)$$
which decreases absolute degree by one, and is 
defined by
$$A_\gamma[\x,i]=\sum_{\y}\sum_{\{\phi\in\pi_2(\x,\y)\big|\Mas(\phi)=1\}}
\#\UnparModFlow(\phi)\cm \langle \gamma,\phi\rangle \cm [\y,i-n_\BasePt(\phi)],$$
where here $\langle \gamma,\phi\rangle$ is the intersection number of
the
restriction of 
$$\phi\colon 
[0,1]
\times\R
\longrightarrow\Sym^g(\Sigma)$$ to $\{1\}\times \R$ (where it maps
into $\Ta$) with the subset of $\left(\gamma\times \Sym^{g-1}(\Sigma)\right)\cap
\Ta$ (c.f. Section~\ref{HolDiskOne:subsec:DefAct} of~\cite{HolDisk}).
We can obviously extend this to $\CFKinf(Y,K)$ by defining
$$A_\gamma[\x,i,j]=\sum_{\y}\sum_{\{\phi\in\pi_2(\x,\y)\big|\Mas(\phi)=1\}}
\#\UnparModFlow(\phi) \cm \langle  \gamma,\phi \rangle \cm [\y,j-n_\BasePt(\phi),
j-n_\FiltPt(\phi)],$$ so that it respects the filtration
(and indeed, it is not hard to see that the filtered chain homotopy
type of the induced map $A_\gamma$ is an invariant of the knot $K$, as well). In
particular, this induces also an action, which we denote by $A_\gamma^0$
on $\HFKa(Y,K)$. In both examples $(Y,K)=(S^2\times S^1,B(0,1))$ and $(\#^2(S^2\times S^1),B(0,0))$
calculated above, for each fixed integer $j$, $\HFKa(Y,K,j)$ 
is concentrated in a single degree, so it follows at once that the action 
$A_\gamma^0$
is trivial. In these cases, there are higher maps, for example
$$A_\gamma^1\colon \HFKa(K,j)\longrightarrow \HFKa(K,j-1),$$
which count $\phi$ with $n_\BasePt(\phi)=0$ and $n_\FiltPt(\phi)=1$
(weighted, once again, by $\langle \gamma,\phi\rangle$).

We can use Proposition~\ref{prop:KnotHomology} to obtain information
about $\HFp$ for some simple three-manifolds which fiber over the
circle.  We begin with the case where the monodromy is trivial.  The
key point for these calculations is that if we perform zero-surgery on
the knot $\#^g B(0,0)\subset \#^{2g}(S^2\times S^1)$, then we get the
three-manifold $\Sigma_g\times S^1$.

We now explain the notation required to state the following result.
Fix an integer $k$, and let $\HFp(\Sigma_g\times
S^1,k)$ denote the summand of $\HFp(\Sigma\times S^1)$ corresponding
to the $\SpinC$ structure $\spinc$ with 
\begin{eqnarray*}
\langle c_1(\spinc),[\Sigma_g]\rangle = 2k, &{\text{and}}&
\langle c_1(\spinc),\gamma \times S^1 \rangle=0
\end{eqnarray*} for all curves $\gamma\subset
\Sigma_g$ (of course, $\HFp$ is trivial on $\SpinC$
structures which do not satisfy this latter condition, by the
adjunction inequality).

Next, for integers $g$ and $d$, consider the group
\begin{equation}
\label{eq:DefX}
X(g,d)=\bigoplus_{i=0}^{d} \Wedge^{2g-i} H^1(\Sigma_g)\otimes_\Z
\left(\Z[U]/U^{d-i+1}\right).
\end{equation} This can be viewed as a module over the ring
$\Z[U]\otimes_\Z \Wedge^* H_1(\Sigma_g)$, where 
the action of $\gamma\in H_1(\Sigma_g)$ 
is given by
\begin{equation}
\label{eq:DefDifferential}
D_\gamma (\omega\otimes U^j)=(\iota_\gamma \omega)\otimes U^j +
\PD(\gamma)\wedge \omega\otimes U^{j+1},
\end{equation}
where here 
$$\iota_\gamma\colon \Wedge^i H^1(\Sigma_g)\longrightarrow
\Wedge^{i-1} H^1(\Sigma_g)$$ denotes the contraction homomorphism
(and $U$ acts in the obvious way). In fact, we can
endow $X(g,d)$ with a relative grading so that multiplication by $U$
lowers degree by two, and $D_\gamma$ lowers
degree by one.  It is interesting to note that if $\Sym^d(\Sigma)$
denotes the $d$-fold symmetric product of $\Sigma$, then 
\begin{equation}
\label{eq:MacDonald}
X(g,d)\cong
H_*(\Sym^d(\Sigma_g))
\end{equation}
as relatively graded groups, according to
MacDonald, c.f.~\cite{MacDonald}. Indeed, $H_*(\Sym^d(\Sigma_g))$ is
naturally a module over $\Z[U]\otimes_\Z\Wedge^* H_1(\Sigma_g)$, where
$U$ acts as cap product with the Poincar\'e dual of $\{x\}\times
\Sym^{d-1}(\Sigma)$, and also $H_1(\Sigma_g)\cong
H^1(\Sym^d(\Sigma_g))$ acts by cap product with one-dimensional
cohomology. With this action, the isomorphism of
Equation~\eqref{eq:MacDonald} is given as $\Z[U]\otimes_\Z\Wedge^*
H_1(\Sigma_g)$-modules, again following~\cite{MacDonald}.

Of course, $\HFp(\Sigma\times S^1)$ is also a module over
$\Z[U]\otimes_\Z \Wedge^* H_1(\Sigma\times S^1)$ and hence, under the
natural inclusion $H_1(\Sigma)$ in $H_1(\Sigma\times S^1)$, it can be
viewed as a module over $\Z[U]\otimes \Wedge^* H_1(\Sigma_g)$ (in
fact, the following proof also shows that the circle factor in
$\Sigma\times S^1$ annihilates $\HFp(\Sigma\times S^1,\spinc)$,
provided $c_1(\spinc)\neq 0$).

It is very suggestive to compare the results on the module structure
with the quantum cohomology ring calculations of the $d$-fold
symmetric product of $\Sigma$ obtained by Bertram and
Thaddeus~\cite{ThadBert}.

\begin{theorem}
\label{thm:SigmaTimesSOne}
Fix an integer $k\neq 0$. Then, there is an identification of $\Z$-modules
\begin{equation}
\label{eq:IdentifyRings}
\HFp(\Sigma\times S^1,k)\cong
X(g,d),
\end{equation}
where $$d=g-1-|k|.$$ (Of course,
the group is trivial if $|k|>g-1$). Indeed, there is a
filtration on $\HFp(\Sigma\times S^1,k)$ as a 
$\Z[U]\otimes_{\Z}\Wedge^*H_1(\Sigma)$-module with the property that 
the identification of Equation~\eqref{eq:IdentifyRings} 
is an isomorphism between
the associated graded modules. Indeed, when $3d<2g-1$, there is
an isomorphism as in Equation~\eqref{eq:IdentifyRings} on the level of 
$\Z[U]\otimes_\Z\Wedge^* H_1(\Sigma)$-modules.
\end{theorem}

\begin{proof}
In view of the conjugation symmetry, it suffices to consider the case
where $k>0$.

First, note that it follows from the calculation of $\HFKa$ for
$B(0,0)$ given in Proposition~\ref{prop:KnotHomology} and
Corollary~\ref{cor:Kunneth} that $\HFKa(\#^{2g}(S^2\times
S^1),\#^{g}B(0,0),j)$ is a free $\Z$-module of rank
$\left(\begin{array}{c} 2g \\ g+j
\end{array}\right)$ supported entirely in grading $j$. More invariantly, we can write
\begin{equation}
\label{eq:MultiBorr}
\HFKa(\#^{2g}(S^2\times S^1),\#^{g}B(0,0),j)\cong 
\Wedge^{g+j}
H^1(\Sigma_g;\Z)_{(j)}
\end{equation}
(where as usual the subscript $(j)$ denotes the absolute grading of
the corresponding group).  

Now we can view $\CFKinf(Y)$ as a $\Z$-filtered subcomplex, by
defining the filtration level of $[\x,i,j]=i+j$. And hence, we obtain
a spectral sequence whose $E_1$ term is $\Z[U,U^{-1}]\otimes_\Z
\HFKa(Y,K)$, converging to $\HFKinf(Y)$. Moreover, as the above calculations show, the 
filtration (on the $E_1$ term) is proportional to the absolute degree
of each group, so the spectral sequence collapses after the $E_2$
stage. The same remarks apply to the complexes $\CFK^{0,*}$ and
$\CFK^{*,0}$, as well. Indeed, we argue that the $d_2$ differentials
vanish, as well. 
This follows immediately from the observations that
$\CFK^{0,*}$ and $\CFK^{*,0}$ must both calculate
$\HFa(\#^{2g}(S^2\times S^1))$, so for dimension reasons, all higher
differentials must vanish.
It also follows readily that
$$H_*(C\{i\geq 0~{\text{or}}~j\geq k\})\cong H\{i\geq 0~{\text{or}}~
j\geq k\},$$
where we have adopted the convention that $H\{i\geq 0~{\text{or}}~
j\geq k\}$ denotes the quotient module of $H_*(C)$ by the submodule of
$H_*(C\{i~\leq 0~{\text{and}}j\leq k\})$. (We use this notation later
in this proof, as well.)

We claim  that with respect to the
induced identification $$\HFinf(\#^{2g}(S^2\times
S^1),\#^{g}B(0,0))\cong \Wedge^{*} H^1(\Sigma_g;\Z)\otimes_\Z
\Z[U,U^{-1}],$$ the action of $$\Z[U]\otimes_\Z \Wedge^*H_1(\#^{2g}(S^2\times
S^1))\cong \Z[U]\otimes_\Z\Wedge^* H_1(\Sigma)$$ is the action
$D_\gamma$ defined above. (Note that the induced action $A_\gamma^0$
is trivial, $A_\gamma^1$ is given by the map $D_\gamma$ above, and all
higher actions vanish for dimension reasons.)

To see why this is true, we proceed as follows. We restrict to
a ``vertical slice''
$$\HFKa(\#^{2g}(S^2\times S^1),\#^g B(0,0))\cong \Wedge^* H^1(\Sigma_g;\Z)
\cong \HFa(\#^{2g}(S^2\times S^1)),$$
where here as usual we view $\CFK^{0,*}(\#^{2g}(S^2\times S^1),\#^g B(0,0))$
as a filtration for 
$\CFa(\#^{2g}(S^2\times S^1))$, by forgetting about
$\FiltPt$.
In the model for $\HFa(\#^{2g}(S^2\times S^1))\cong \Wedge^*
H^1(\Sigma_g;\Z)$, the action by $\gamma$ is given by $\iota_\gamma$
(compare Section~\ref{HolDiskTwo:sec:Examples} of~\cite{HolDiskTwo}).
This verifies that if 
$$\omega\otimes U^j\in \Wedge^{2g-i}
H^1(\Sigma_g)\otimes_\Z\left(\Z[U]/U^{d-i+1}\right),$$ then the component of
$A_\gamma(\omega\otimes U^j)$ in $\Wedge^{2g-i-1} H^1(\Sigma_g)\otimes
_\Z\left(\Z[U]/U^{d-i+1}\right)$ -- the ``vertical component'' -- 
is given by $D_\gamma(\omega\otimes U^j)$. For
the component in $\Wedge^{2g-i+1}\otimes_\Z\left(\Z[U]/U^{d-i+1}\right)$, we
consider now $\CFKa^{*,0}(\#^{2g}(S^2\times S^1),\#^g B(0,0))$ as a filtration
of $\CFa(\#^{2g}(S^2\times S^1))$, by forgetting about $\BasePt$.
In this model, 
$$\HFa(\#^{2g}(S^2\times S^1))\cong \bigoplus_j
\Wedge^{g+j} H^1(\Sigma_g;\Z)\otimes U^j.$$ Now, the
the action of $\gamma$ is forced to be given by the map
$$(\omega\otimes U^j) \mapsto \PD(\gamma)\wedge \omega \otimes U^{j+1},$$
verifying that the 
$D_\gamma$ describes the horizontal component of $A_\gamma$. Since $A_\gamma$
decreases dimension by one, there are no other components to be verified.

Let $Y_n$ denote the three-manifold obtained by performing $+n$
surgery along $B$. This manifold is, in fact, the total space of a
circle bundle over $\Sigma$ with Euler number $n$. It follows from the
above remarks, together with Theorem~\ref{thm:LargePosSurgeries} that
for large enough $n$,
$\HFpRed(Y_n)=0$. Indeed, if we consider the cobordism $W$
obtained by the natural two-handle addition, and let
$f_1\colon \HFp(Y_n,\spinc_k)\longrightarrow
\HFp(\#^{2g}(S^2\times S^1))$ denote the map associated to the $\SpinC$ structure
$\spincfour$ with 
$$\langle c_1(\spincfour), [S]\rangle = 2k+n$$
(where $S$ is is obtained by completing the Seifert surface for $B$ inside
$W$), then we see that the kernel of $f_1$
is given by
\begin{eqnarray*}
H\{i<0,j\geq k\}&\cong& X(g,d),
\end{eqnarray*}
where $d=g-1+|k|$. 
We claim in fact that the homology of this kernel of $f_1$ 
is
identified with $\HFp(S^1\times \Sigma,k)$. 

To this end, recall the integer surgeries long exact sequence,
according to which we have
$$
\begin{CD}
... @>{F}>>\HFp(\#^{2g}(S^1\times S^2))
@>>>\HFp(S^1\times \Sigma,k)
@>{G}>>
\HFp(Y_n,[k])
@>>> ..., 
\end{CD}
$$ where here the map $F$ is the sum of the maps associated to all the
$\SpinC$ structures over the cobordism $W$.  Note also that since
$\HFpRed(\#^{2g}(S^1\times S^2))=0$, and $\HFpRed(S^1\times
\Sigma,k)\cong \HFp(S^1\times \Sigma,k)$ (since $k\neq 0$), it follows
that $\HFp(S^1\times \Sigma_g)\cong \Ker F$.  Now, we can write
$$F=f_1+f_2,$$ where $f_1$ is the homogeneous map described earlier,
and $f_2$ has lower order, in the sense that if $\xi$ is a homogeneous
element of dimension $d$, then $f_2(\xi)$ is a sum of elements of
dimension less than $f_1(\xi)$. More concretely,
according to Equation~\eqref{eq:DimensionFormula},
if $\xi$ supported in a single dimension, then so is $f_1(\xi)$, 
and $f_2(\xi)$ can be written as a sum of terms, all of which are
supported in dimension at least $2k$ smaller than the dimension of $f_1(\xi)$.

In this setting, there is an identification $A=\Ker f_1$ with
$B=\Ker F$, as $\Z$-modules.
To see this, let $$R_1\colon H\{i\geq 0\} \longrightarrow
H\{i\geq 0~{\text{or}}~j\geq k\}$$ denote the natural
inclusion. Since $f_2$ has lower order than $f_1$, it is easy to see
that there is an injective $\Z$-linear map 
$$\Pi\colon \Ker(f_1+f_2)\longrightarrow \Ker(f_1),$$
given by taking the an element $\xi\in A$ to
its homogeneous term in the highest degree, and hence, 
$\rk(B)\geq \rk(A)$. In the other direction, consider the map
\begin{equation}
\label{eq:DefR}
R=\sum_{i=0}^{\infty}(-R_1 \circ f_2)^i.
\end{equation}
Since $R$ is non-increasing in the filtration induced by the absolute
grading, while $f_2$ is strictly decreasing in this filtration (and,
of course, this filtration is bounded below), it follows that $R(\xi)$
is well-defined for any $\xi$. Indeed, it is easily seen that $R$
defines an injection from $A$ to $B$, and hence, their ranks agree.

To analyze the $\Z[U]\otimes_\Z\Wedge^* H_1(\Sigma)$-module structure,
we proceed as follows. Of course, the filtration mentioned in the
statement of the theorem is induced by identifying $\HFp(S^1\times
\Sigma)=\Ker F\subset \HFp(Y_n)$, and noting that the latter group has
a filtration induced by absolute grading. The first assertion about
the module structure is clear, now. For the second, observe that in
general neither the maps $\Pi$ nor $R$ are module
homomorphisms. However, the restriction of $R_1$ to $X=H\{i\geq
0~\text{and}~j< k\}$ does respect the module structure (as does $f_2$,
since it is induced by a cobordism). In fact, since in our case, the
absolute degree of $H\{i,j\}$ is given by $i+j$, it follows that any
element of $\HFp(\#^{2g}(S^2\times S^1))$ with degree $<g-k$ is
contained in $X$.  Since the maximal dimension of any element of
$H\{i< 0~\text{and}~ j\geq k\}$ is $g-2$, and $f_2$ decreases absolute
degree by at least $2k$, the inequality $2g-1>3d$ implies that $f_2$
maps $H\{i< 0~\text{and}~j\geq k\}$ into $X$, and hence $R$ is a
module homomorphism and hence also an isomorphism.
\end{proof}

More generally, let $\phi$ be an automorphism of $\Sigma_g$, and let
$M_\phi$ denote the mapping torus of $\phi$, $$M_{\phi}\cong
\frac{[0,1]\times \Sigma_g}{(0,x)\sim (1,\phi(x))}.$$ When $\phi$ can
be written as a product of disjoint non-separating Dehn twists, one
can use the above techniques to calculate $\HFp(M_\phi,k)$ in many cases.  We
content ourselves here with a calculation in the case where $\phi$ is
a single negative Dehn twist along a non-separating simple, closed
curve $\gamma\subset \Sigma$. 

To this end, we define another differential on $X(g,d)$. Write $\Sigma=T\# \Sigma'$,
where $T$ is a torus and $\gamma\subset T$ is a non-separating simple, closed curve in
$T$. Let
\begin{eqnarray*}
\Wedge^*_+ H^1(\Sigma)&=&\left(\Wedge^0 H^1(T)\otimes_\Z \Wedge^*(\Sigma')\right) \oplus
\left(\Wedge^2 H^1(T)\otimes_\Z \Wedge^*(\Sigma')\right) \\
\Wedge^*_- H^1(\Sigma)&=&\Wedge^1 H^1(T)\otimes_\Z \Wedge^*(\Sigma'), 
\end{eqnarray*}
so that $\Wedge^* H^1(\Sigma)=\Wedge^*_+H^1(\Sigma)\oplus\Wedge^*_-H^1(\Sigma)$.
Let $X(g,d)=X_+(g,d)\oplus X_-(g,d)$ be the corresponding splitting.
We define $D'_\gamma$ so that 
\begin{eqnarray}
\label{eq:DefNewDifferential}
D'_{\gamma}|X_+(g,d)\equiv 0 
&{\text{and}}&
D'_{\gamma}|X_-(g,d)\equiv D_\gamma|X_-(g,d)
\end{eqnarray}

\begin{theorem}
\label{thm:DehnTwist}
Let $\phi$ denote a negative Dehn twist along a
non-separating simple, closed curve $\gamma\subset \Sigma$, and let
$M_\phi$ denote its mapping torus. For each integer $k$,  
let $\HFp(M_\phi,k)$ denote
$\HFp$ of $M_\phi$ evaluated on the $\SpinC$ structure $\spinc$
characterized by the properties that $\langle
c_1(\spinc),[F]\rangle=2k$ and $c_1(\spinc)$ vanishes on all homology
classes represented by tori. 
Suppose that $0\neq k$ and also suppose that
\begin{equation}
\label{eq:DegreeInequality}
3d<2g-1,
\end{equation}
where $d=g-1-|k|$.  Then, there is an isomorphism of $\Z[U]\otimes_\Z
\Wedge^*H_1(\Sigma)$-modules 
\begin{equation}
\label{eq:DehnTwist}
\HFp(M_\phi,k)\cong H_*(X(g,d),D'_\gamma),
\end{equation}
where $D'_\gamma$ is the differential defined by
Equation~\eqref{eq:DefNewDifferential}.
When, Equation~\eqref{eq:DegreeInequality} does not hold, we still obtain 
an isomorphism of the form stated in Equation~\eqref{eq:DehnTwist},
as $\Z$-modules.
\end{theorem}

\begin{proof}
If $\phi$ is a negative Dehn twist along $\gamma$, then we can realize
$M_\phi$ as zero-surgery along the connected sum
$B(0,1)\#\left(\#^{g-1}B(0,0)\right)$. 
We adopt the notation from the proof of Theorem~\ref{thm:SigmaTimesSOne}
before. Note that $\HFKa$ of $B(0,0)$ and $B(0,1)$ differ
by only a shift in dimension, so
it is still the case that the knot filtration is proportional to 
the absolute degree, so the spectral sequence collapses after the $E_2$
stage. 
There is one crucial difference now:  it is no longer the case that
the $d_2$ differential is trivial. Indeed, since $\HFKa(S^2\times
S^1,B(0,1))$ is the homology of the graded object associated to a
filtration of $\CFa(S^2\times S^1)$, and $$\HFa(S^2\times S^1)\cong
\Z_{(-\OneHalf)}\oplus \Z_{(\OneHalf)},$$ it follows that the $d_2$
differential is non-trivial; indeed, the map $$d_2\colon
\Z^2_{(-\OneHalf)}\cong \HFKa(S^2\times S^1,B(0,1),0)\longrightarrow
\HFKa(S^2\times S^1,B(0,1),-1)\cong
\Z_{(-3/2)}$$ 
surjects. In the notation from the theorem, we can write this as the
map $H^1(T^2)\longrightarrow H^0(T^2)$ given by evaluation against
$[\gamma]$ (i.e. this is the ``horizontal component'' of $D_{\gamma}'$.

We now proceed as before, now letting $H$ denote the homology of $C$
with respect to the $d_1$ differential. Thus, the group $H\{i<0,j\geq
k\}\cong X(g,d)$, which is the kernel of $f_1$, is now endowed with
the $d_2$ differential induced by $D'_\gamma$. (The fact that this
differential is given by $D'_\gamma$ follows the proof of
Theorem~\ref{thm:SigmaTimesSOne}, in the calculation of the module
$\Z[U]\otimes_\Z \Wedge^* H_1(\Sigma_g)$-module structure of
$\HFp(\Sigma\times S^1,k)$. In particular, it follows from dimension
reasons, followed by a consideration of the filtration
$\CFK^{*,0}(S^2\times S^1,B(0,1))$ of $\CFa(S^2\times S^1)$.)

We once again consider  the integral surgeries long exact sequence,
and let $Y_n$ denote the three-manifold obtained by large $n$ surgery
on the link $B(1,0)\#^{g-1}B(0,0)$, and 
identify $\HFp(M_\phi,k)\cong H_*(\Ker(f_1+f_2),d_2)$.
Proceeding as before, we construct the identification of $\Ker(f_1+f_2)$ with
$\Ker(f_1)$, by considering the map $R$ defined as in Equation~\eqref{eq:DefR}.
The inequality $3d<2g-1$ now ensures that $R$ is a chain map,
using the $d_2$ differential (and hence
also an isomorphism of chain complexes).

For the final remark, when Equation~\eqref{eq:DegreeInequality} is
violated, we still observe that the map induced by $f_1$ is surjective
on homology.  We then appeal to the same argument as in
Corollary~\ref{cor:HighestDegree}.
\end{proof}

It follows from the above result that when $g>2$, $\HFp(M_\phi,g-2)$ is
isomorphic to the relative singular cohomology $H^*(\Sigma,\gamma)$,
while $\HFp(M_{\phi^{-1}},g-2)\cong H_*(\Sigma-\gamma)$.  This should be
compared with a result of Seidel~\cite{Seidel}, see
also~\cite{HutchingsSullivan}.

\section{Links}
\label{sec:Links}

Throughout most of the paper, we assumed that our oriented knot $K$ is
connected. In view of Proposition~\ref{prop:LinksToKnots}, passing to
the case of oriented links is quite straightforward. Our aim here is
to highlight some issues in this generalization, with the aim of
showing how the properties of the link invariant stated in the
introduction follow from the other results proved above.

As a first point, recall that in some applications, especially to
exact sequences (c.f. Section~\ref{sec:Sequences}), we found it
convenient to fix a Seifert surface $F$ for the knot
$K$. Correspondingly, when considering a link $L$, we also fix a
Seifert surface $F$ for $L$. 

\begin{lemma}
\label{lemma:SplittingWellDefined}
Fix a Seifert surface $F$ for $L$ in $Y$, and let $F'$ be any
extension of $F$ to the zero-surgery of $\kappa(Y)$ along
$\kappa(L)$. Then, the group $$\HFKa(Y,L,m)=
\bigoplus_{\{\relspinct\in\RelSpinC(\kappa(Y),\kappa(L))\big|
\langle c_1(\relspinct),[F']\rangle = 2m\}}\HFKa(Y,L,\relspinct)$$
depends on $F'$ only through the 
(relative homology class of the) induced Seifert surface $F$ in $Y$
(i.e. it is independent of the extension $F'$).
\end{lemma}

\begin{proof}
The indeterminacy of $F'$ comes from the attached one-handles --
more precisely
writing $\kappa(Y)=Y\# \#^{n-1}(S^2\times S^1)$, the
indeterminacy amounts to adding homology classes classes represented
by spheres coming from $S^2\times S^1$, which in turn
can be represented by embedded tori in the complement of $\kappa(K)$
inside $\kappa(Y)$. However, it is
straightforward to see that if the first Chern class of $\relspinct\in
\RelSpinC(\kappa(Y))$ evaluates non-trivially on such a homology
class, then we can find an admissible Heegaard diagram for the knot
$\kappa(L)$ in $\kappa(Y)$ with the property that no intersection
point represents $\relspinct$. This follows from a straightforward
adaptation of the proof of the adjunction inequality for $\HFp$
(Theorem~\ref{HolDiskTwo:sec:Adjunction} of~\cite{HolDiskTwo}, see also
Theorem~\ref{thm:Adjunction} above)
that $\HFKa(\kappa(Y),\kappa(L),\relspinct)=0$.
\end{proof}

Now, consider the case where $Y$ is an oriented three-manifold,
equipped with an oriented link $L_+$, and let $\gamma\subset Y$ be an
unknot which spans a disk $D$ which meets $L$ in two algebraically
cancelling transverse points. Let $L_-$ denote the new link induced
from $L_-$ after introducing a full twist, as pictured in
Figure~\ref{fig:Skein}. In comparing with this picture, the disk $D$
can be thought of as a horizontal disk.  Let $L_0$ denote the new link
obtained from $L$ by resolving the link, so that it no longer meets
$D$.

Recall that in Section~\ref{sec:Sequences}, we chose a Seifert surface
$F$ for $L$, and use it to define the integral splitting of
$\HFK(Y,L)$ as in Lemma~\ref{lemma:SplittingWellDefined}.

\begin{theorem}
\label{thm:SkeinExactSequence}
Let $L_+$, $L_0$, and $L_-$ be the links related by a skein move as above.
Then we have a long exact sequence of the form:
\begin{equation}
\label{eq:SkeinOne}
\begin{CD}
...@>>>\HFKa(Y,L_-)@>>> \HFKa(Y,L_0) @>>> \HFKa(Y,L_+)@>>>...,
\end{CD}
\end{equation}
while if they belong to different components of $L$, we have a long
exact sequence of the form
\begin{equation}
\label{eq:SkeinTwo}
\begin{CD}
...@>>>\HFKa(Y,L_-)@>>> \HFKa(Y',L_0')
@>>> \HFKa(Y,L_+)@>>>...,
\end{CD}
\end{equation}
where here $(Y',L_0')$ denotes the link obtained by connected sum of
$(Y,L_0)$ with the ``Borromean knot'' $(\#^2(S^2\times S^1),B(0,0))$
from Figure~\ref{fig:BorromeanKnot}.  Moreover, the maps in the
exact sequence respect the $\Z$-splittings of the homology groups
induced from a Seifert
surface $F$ for $L_+$ in $Y$ (and corresponding Seifert surface for
$L_-$ and $L_0$).
\end{theorem}

\begin{figure}
\mbox{\vbox{\epsfbox{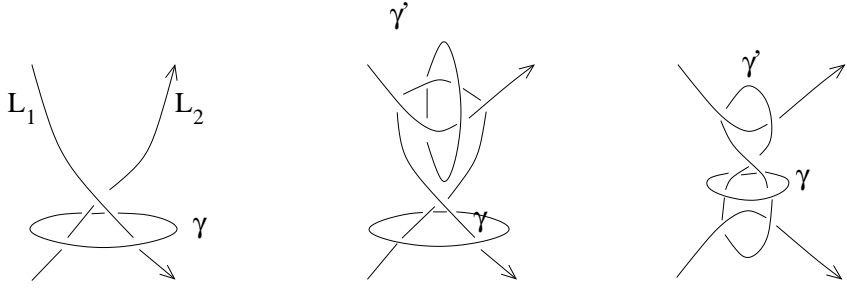}}}
\caption{\label{fig:TradeLinkCrossings}
{\bf{Trading link crossings.}}  The identification of 
$(\kappa(Y)_0(\gamma),\kappa(L))$ with
$(Y,L_0)\#(S^2\times S^1,B(0,0))$.}
\end{figure}

\begin{proof}
This result is a special case of Theorem~\ref{thm:LongExactSeq}, after some
remarks.

First, observe that under the diffeomorphism $Y_{-1}(\gamma)\cong Y$,
the image of the link obtained from viewing $L_+$ as a link in
$Y-\gamma\subset Y_{-1}(\gamma)$ is $L_-$.  

Suppose next that both strands of $L$ meeting $D$ belong to the same
component of $L$. In this case, we can trade the zero-framed
two-handle attached along $\gamma$ for a one-handle, without changing
the underlying three-manifold $Y_0(\gamma)\cong Y\# (S^2\times
S^1)$. In the case where $L$ had two components, in fact, we see that
the knot induced from $K$ then coincides with the knot
$\kappa(L_0)$. It is easy to see that in this case,
Theorem~\ref{thm:LongExactSeq} translates into Long Exact
Sequence~\eqref{eq:SkeinOne}.

Suppose next that two strands of $L$ meeting in $D$ belong to
different components of $L$. In fact, for simplicity, we assume that
$L$ has two components $L_1$ and $L_2$ (indeed, the case where $L$
consists of more than two components can be reduced to this case,
after attaching sufficiently many one-handles.) In this case, by
definition, the link invariant of $L$ inside $Y_0(\gamma)$ agrees with
a knot invariant of the knot induced from $\kappa(L)$ inside
$\kappa(Y)_0(\gamma)$. Recall that $\kappa(Y)$ is the three-manifold
obtained from $Y$ by attaching a one-handle. We can replace the
one-handle in $\kappa(Y)$ by a zero-framed two-handle $\gamma'$, to
obtain $Y'=Y\#^2(S^2\times S^1)$ with a new knot $K'$. It is easy
to see $$(Y',K')=(Y,K)\#(\#^2(S^2\times S^1),B(0,0))$$ as claimed,
with the knots $\gamma$ and $\gamma'$ playing the role of the two
zero-framed circles in the Borromean knot, as illustrated in
Figure~\ref{fig:TradeLinkCrossings}.  Strictly speaking, to see that
the operation is local on the knot $K$, it is useful to trade the
two-handle specified by $\gamma$ for a three-handle (with the curve
$\gamma'$ running through it), and slide the one-handle freely along
$K$.

With these remarks in hand, we see Long Exact
Sequence~\eqref{eq:SkeinTwo} as a special case of
Theorem~\ref{thm:LongExactSeq}.
\end{proof}

In view of the calculation of $\HFKa(\#^2(S^2\times S^1),B(0,0))$ in
Proposition~\ref{prop:KnotHomology}, together with the K\"unneth
principle for connected sums from Corollary~\ref{cor:Kunneth}, we see
that the versions of the skein exact sequence given in the
introduction coincide with those stated in
Theorem~\ref{thm:SkeinExactSequence}.

The Euler characteristic calculation of $\HFKa$ stated in
Equation~\eqref{eq:EulerChar} is an easy application of the skein
exact sequence (though an alternate proof could also be given in the
spirit of Theorem~\ref{HolDiskTwo:thm:EulerOne} of~\cite{HolDiskTwo}).
More precisely, note that $\HFKa(S^3,L)$ inherits an absolute
$\Zmod{2}$ grading from $\HFa(\#^{n-1}(S^2\times S^1))$ (where here
$n$ is the number of components of $L$).  As in~\cite{AbsGraded}, we
see that this grading is preserved by both the map from
$\HFKa(S^3,L_-)$ to $\HFKa(S^3,L_0)$ and the map from $\HFKa(S^3,L_+)$
to $\HFKa(S^3,L_-)$, while it is reversed by the remaining map. It
follows from the two skein exact sequences then that if we write
$$\chi\left(\HFKa(S^3,L)\right)=\sum_{i\in\Z}
\chi\left(\HFKa(S^3,L,i)\right)\cm T^i$$ then
$$\chi\left(\HFKa(S^3,L_-)\right)-\chi\left(\HFKa(S^3,L_0)\right)-\chi\left(\HFKa(S^3,L_+)\right)=0$$
if $L_0$ has more components than $L_+$; otherwise,
$$\chi\left(\HFKa(S^3,L_-)\right)-(T^{1/2}-T^{-1/2})^2\cm
\chi\left(\HFKa(S^3,L_0)\right)-\chi\left(\HFKa(S^3,L_+)\right)= 0.$$
Since $\chi(\HFKa(S^3,u))=1$ (when $u$ is the unknot),
Equation~\eqref{eq:EulerChar} now follows from the usual skein
relation relation characterization of the Alexander-Conway polynomial.
The fact that this absolute $\Zmod{2}$ grading is related to the
absolute $\Q$-grading in the classical case discussed in the
introduction follows easily from~\cite{AbsGraded}.

Equation~\eqref{eq:Reflection} follows immediately from
Proposition~\ref{prop:IndepOrient},
Equation~\eqref{eq:KnotConjInv} follows from
Proposition~\ref{prop:JInvariance}, 
Equation~\eqref{eq:LinkOrientationReversal} follows
from Proposition~\ref{prop:JInvarianceGen}
(noting that when we pass from $K$ to $-K$, we also must reverse the 
orientation of the Seifert surface, so
$\langle c_1(\spinc),[{\widehat F}]\rangle = 
\langle c_1(J\spinc),[-{\widehat F}]\rangle$), 
and 
Equation~\eqref{eq:KunnethFormula} follows from
Corollary~\ref{cor:Kunneth}.

For the behavior of $\HFKa$ under disjoint union, note that by definition
of the link invariant,
$$\HFKa(S^3,L_1\cup L_2)=\HFKa(S^2\times S^1, L_1\# L_2)$$
(where here we think of $L_1\# L_2$ as supported in a ball in $S^2\times S^1$).
Now, Equation~\eqref{eq:DisjointUnion} from a straightforward adaptation
of the proof that $\HFa(S^3\#(S^2\times S^1))\cong \HFa(Y)\otimes H_*(S^1)$
given in~\cite{HolDiskTwo}.

\commentable{
\bibliographystyle{plain}
\bibliography{biblio}
}

\end{document}